\documentclass[12pt,a4paper,subeqn,reqno]{article}

\usepackage{graphicx}
\usepackage{amssymb}
\usepackage{epstopdf}

\usepackage{latexsym}

\usepackage{amsmath,amsthm}

\usepackage[mathscr]{eucal}
\usepackage[all]{xy}
\usepackage{mathrsfs}
\usepackage{authblk}

\usepackage{amsfonts}
\usepackage{subfigure}
\usepackage{caption}
\usepackage[dvipsnames]{xcolor}
\usepackage{shuffle}
\usepackage{makeidx}
\usepackage[all]{xy}
\usepackage{shuffle}
\usepackage{float}
\usepackage[draft]{fixme}

\allowdisplaybreaks
\usepackage[text={169mm,245mm},centering]{geometry}

\usepackage{mathtools}
\usepackage{dsfont}
\usepackage{url}
\usepackage{enumerate}
\usepackage[figuresright]{rotating}

\usepackage{datetime2}

\usepackage[hyperindex=true,hidelinks]{hyperref}

\usepackage{scalerel,stackengine}
\stackMath
\newcommand\rwh[1]{%
\savestack{\tmpbox}{\stretchto{%
  \scaleto{%
    \scalerel*[\widthof{\ensuremath{#1}}]{\kern.1pt\mathchar"0362\kern.1pt}%
    {\rule{0ex}{\textheight}}
  }{\textheight}%
}{2.4ex}}%
\stackon[-6.9pt]{#1}{\tmpbox}%
}
\newtheorem{theorem}{Theorem}[section]

\newtheorem{lemma}{Lemma}[section]
\newtheorem{conjecture}[lemma]{Conjecture}
\newtheorem{corollary}[lemma]{Corollary}
\newtheorem{proposition}[lemma]{Proposition}

\theoremstyle{definition}

\newtheorem{example}[lemma]{Example}
\newtheorem{definition}[lemma]{Definition}
\newtheorem{definition-lemma}[lemma]{Definition-Lemma}
\newtheorem{definition-theorem}[lemma]{Definition-Theorem}

\newtheorem{remark}[lemma]{Remark}

\newcommand{\ka}{\kappa}
\newcommand{\om}{\omega}
\newcommand{\ph}{\varphi}
\newcommand{\tht}{\theta}
\newcommand{\eps}{\varepsilon}
\newcommand{\ue}{\underline\eps}
\newcommand{\uvar}{\ue}

\newcommand{\uK}{\underline{K}}
\newcommand{\ul}{{\underline{\ell}}}
\newcommand{\uh}{{\underline{h}}}
\newcommand{\uo}{{\underline{1}}}
\newcommand{\up}{{\underline{p}}}
\newcommand{\us}{{\underline{s}}}
\newcommand{\uS}{{\underline{S}}}
\newcommand{\uT}{{\underline{T}}}

\newcommand{\one}{{\mathds 1}}
\newcommand{\uone}{{\underline{\one}}}

\newcommand{\whp}{\wh{p}}

\definecolor{rossocorsa}{rgb}{0.83, 0.0, 0.0}

\newcommand{\bl}[1]{\textcolor{blue}{#1}}
\newcommand{\bbl}[1]{\textcolor{blue}{#1}}
\newcommand{\rred}{\color{red}}
\newcommand{\red}[1]{\textcolor{red}{#1}}
\newcommand{\newc}[1]{\textcolor{rossocorsa}{#1}} 
\newcommand{\bblack}{\color{black}}

\newcommand{\marge}[1]{\marginpar{{\rred\scriptsize\em #1}}}

\newcommand{\datestamp}{{\small{File:\;\hbox{\tt\jobname.tex}
\; \DTMnow}}}
\newcommand\extrafootertext[1]{%
    \bgroup
    \renewcommand\thefootnote{\fnsymbol{footnote}}%
    \renewcommand\thempfootnote{\fnsymbol{mpfootnote}}%
    \footnotetext[0]{#1}%
    \egroup
  }
  \newcommand\blfootnote[1]{
    \begingroup
    \renewcommand\thefootnote{}\footnote{#1}
    \addtocounter{footnote}{-1}
    \endgroup
}

\DeclareMathOperator{\den}{den}
\DeclareMathOperator{\SU}{SU}
\DeclareMathOperator{\SL}{SL}
\DeclareMathOperator{\GL}{GL}
\DeclareMathOperator{\C}{\mathbb{C}}

\DeclareMathOperator{\Span}{Span}

\DeclareMathOperator{\cB}{\mathcal{B}}

\DeclareMathOperator{\Q}{\mathbb{Q}}
\DeclareMathOperator{\sQ}{\mathscr{Q}}
\DeclareMathOperator{\cQ}{\mathcal{Q}}
\newcommand{\ocQ}{\overrightarrow{\mathcal{Q}}}
\newcommand{\oQs}{{\sQ_{\!s}}} 

\newcommand{\cV}{\mathscr{V}}
\newcommand{\cVs}{\cV_{\!s}}
\newcommand{\cW}{\mathscr{W}}

\newcommand{\cO}{\mathcal{O}}

\newcommand{\gP}[1]{\mathscr{P}_{\hspace{-.2ex}#1}}
\DeclareMathOperator{\gR}{\mathscr{R}}

\newcommand{\cT}{\mathcal{T}_r}

\DeclareMathOperator{\Z}{\mathbb{Z}}
\newcommand{\JX}{\mathscr{J}_{\!X}}

\DeclareMathOperator{\cE}{\mathcal{E}}

\DeclareMathOperator{\cN}{\mathcal{N}}

\DeclareMathOperator{\cL}{\mathcal{L}}

\DeclareMathOperator{\fH}{\mathfrak{H}}
\DeclareMathOperator{\fL}{\mathfrak{L}}
\DeclareMathOperator{\fR}{\mathfrak{R}}
\DeclareMathOperator{\fS}{\mathfrak{S}}

\DeclareMathOperator{\Tr}{Tr}
\DeclareMathOperator{\M}{\mathcal{M}}

\DeclareMathOperator{\CS}{CS}
\DeclareMathOperator{\Irr}{Irr}

\DeclareMathOperator{\spin}{spin}
\DeclareMathOperator{\ev}{ev}
\DeclareMathOperator{\Gal}{Gal}
\DeclareMathOperator{\rank}{Rank}
\DeclareMathOperator{\Res}{Res}
\DeclareMathOperator{\interior}{int}
\DeclareMathOperator{\PSL}{PSL}
\DeclareMathOperator{\ID}{id}

\numberwithin{equation}{section}

\newcommand{\col}{\colon\thinspace}          

\DeclareMathOperator{\bS}{\mathcal{S}}
\DeclareMathOperator{\cS}{{\mathscr{S}_{\!\CS}}}

\DeclareMathOperator{\R}{\mathbb{R}}
\DeclareMathOperator{\HH}{\mathbb{H}}
\DeclareMathOperator{\DD}{\mathbb{D}}

\DeclareMathOperator{\IM}{\Im\hspace{.04em}\textrm{m}}
\DeclareMathOperator{\RE}{\Re\hspace{.04em}\textrm{e}}

\DeclareMathOperator{\WRT}{WRT}
\newcommand{\WRTk}{\WRT_k}
\newcommand{\WRTone}{\WRT_1}

\newcommand{\ii}{^{-1}}
\newcommand{\al}{\alpha}
\newcommand{\wh}{\widehat}
\newcommand{\hh}{\widehat}
\newcommand{\ti}{\tilde}
\newcommand{\wt}{\widetilde}

\newcommand{\ga}{\gamma}
\newcommand{\uga}{{\underline{\ga}}}

\newcommand{\Ga}{\Gamma}

\def\doublestroke#1{\pdfliteral{1 Tr .5 w}#1\pdfliteral{0 Tr 0 w}}
\newcommand{\GA}{\doublestroke{\Ga}}
\newcommand{\TGA}{\wt\GA}
\newcommand{\gTGA}{\big\langle g\bulM\wt\GA \big\rangle}
\newcommand{\bull}{\raisebox{.05ex}{\scalebox{.7}{$\bullet$}\!}}
\newcommand{\bulM}{\,\bull\,\,}

\newcommand{\De}{\Delta}
\newcommand{\La}{\Lambda}
\newcommand{\sig}{\sigma}

\newcommand{\Om}{\Omega}
\newcommand{\Omb}{{  \hspace{.1em} \raisebox{.6ex}{\scalebox{1}[.55]{$|$}} \hspace{-.41em}\scalebox{1.15}[1]{$\Omega$}}}
\newcommand{\Ombb}{{  \hspace{.15em} \raisebox{.36ex}{\scalebox{1}[.4]{$|$}} \hspace{-.355em}\scalebox{1.2}[1]{$\scriptstyle\Omega$}}}

\newcommand{\la}{\lambda}
\newcommand{\ze}{\zeta}

\newcommand{\ens}{\enspace}
\newcommand{\iimp}{\;\Rightarrow\;}
\newcommand{\imp}{\ens\Rightarrow\ens}
\newcommand{\Imp}{\quad\Rightarrow\quad}

\newcommand{\oone}[1]{ \left[{#1}\right]_{\raisebox{-.05ex}{\!$\scriptstyle 1
      $}} } 
\newcommand{\ttwo}[1]{ \left[{#1}\right]_{\raisebox{-.05ex}{$\scriptstyle 2
      $}} } 
\newcommand{\tp}[1]{ \left[{#1}\right]_{\raisebox{-.05ex}{$\scriptstyle 2P
      $}} } 
\newcommand{\tpj}[1]{ \left[{#1}\right]_{\raisebox{-.05ex}{$\scriptstyle 2p_j
      $}} } 
\newcommand{\fourp}[1]{ \left[{#1}\right]_{\raisebox{-.05ex}{$\scriptstyle 4P
      $}} } 
\newcommand{\pj}[1]{ \left[{#1}\right]_{\raisebox{-.05ex}{$\scriptstyle p_j
      $}} } 
\newcommand{\po}[1]{ \left[{#1}\right]_{\raisebox{-.05ex}{$\scriptstyle 2p_1
      $}} } 
\newcommand{\brN}[1]{ \left[{#1}\right]_{\raisebox{-.05ex}{$\scriptstyle N$}} } 

\newcommand{\Mss}{ \M_{\raisebox{-.35ex}{$\scriptstyle >s$}} }

\newcommand{\begla}{\begin{equation}}
\newcommand{\beglab}[1]{\begin{equation}	\label{#1}}
\newcommand{\edla}{\end{equation}}

\newcommand{\beglano}{\begin{equation*}}
  \newcommand{\edlano}{\end{equation*}}

\newcommand{\rhs}{right-hand side}
\newcommand{\lhs}{left-hand side}

\newcommand{\bJ}{\mathbf{J}}

\newcommand{\sN}{\cN_{\hspace{-.4ex}*}}

\newcommand{\indfS}{\mathds1_{\fS}}
\newcommand{\induh}{\mathds1_{\fS^\uh}}
\newcommand{\indul}{\mathds1_{\fS^\ul}}
\newcommand{\indulp}{\mathds1_{\fS^{\ul'}}}
\newcommand{\indusl}{\mathds1_{\fS^{\sig_1(\ul)}}}

\DeclareMathOperator{\med}{med}
\newcommand{\median}{_{\med}}

\newcommand{\matr}[2]{\begin{pmatrix}#1\\#2\end{pmatrix}}

\newcommand{\compl}{ \raisebox{.6ex}{$\scriptstyle\complement$}\!}
\newcommand{\ccompl}{ \raisebox{.4ex}{{\scalebox{.65}{$\scriptstyle\complement$}}}\!}

\newcommand{\CEv}{ \raisebox{.6ex}{$\scriptstyle\complement$}\hspace{-.25ex} Ev}
\newcommand{\cEv}{
  \raisebox{.4ex}{{\scalebox{.65}{$\scriptstyle\complement$}}}\hspace{-.15ex} Ev}

\makeatletter
\DeclareRobustCommand\dalemb{\mathpalette\inner@dalemb{}}
\def\inner@dalemb#1{%
  \add@dalemb#1{03}%
  \add@dalemb#1{06}%
  \square
}
\def\add@dalemb#1#2{%
  \sbox0{\scalebox{1.#2}{$#1\square$}}%
  \rlap{\lower0.#2\ht0\box0}%
}
\makeatother

\newcommand{\pp}[1]{^{(#1)}}

\newcommand{\defeq}{:=}

\newcommand{\ie}{i.e.\ }

\newcommand{\fakeheight}[2]{\mbox{\smash{#1}\vphantom{#2}}}

\newcommand{\modM}{\;\textrm{mod}\,2P}
\newcommand{\modcM}{\;\textrm{mod}\,2cP}

\usepackage[style=alphabetic,backend=bibtex8,maxbibnames=10]{biblatex}
\begin{document}

\title{\textbf{A proof of Witten's asymptotic expansion conjecture for WRT invariants of Seifert fibered homology spheres}}

\author{J\o rgen Ellegaard Andersen, Li Han, Yong Li, 
\\ William Elb\ae k Misteg\aa rd, David Sauzin, Shanzhong Sun}

\date{}
\maketitle


\begin{abstract}
  %
  Let~$X$ be a general Seifert fibered integral homology $3$-sphere
  with $r\ge3$ exceptional fibers.
  For every root of unity $\ze\not=1$, we show that the $\SU(2)$ WRT
  invariant of~$X$ evaluated at~$\ze$ is (up to an elementary 
  factor) the non-tangential limit at~$\ze$ of the GPPV invariant
  of~$X$, thereby generalizing a result from \cite{AM22}.
  Based on this result, we apply the quantum modularity results
  developed in \cite{han2022resurgence} to the GPPV invariant of~$X$
  to prove Witten's asymptotic expansion conjecture \cite{Witten98}
  for the WRT invariant of~$X$.
  We also prove that the GPPV invariant of~$X$ induces a higher depth
  strong quantum modular form.
  Moreover, when suitably normalized, the GPPV invariant
    provides an ``analytic incarnation'' of the Habiro invariant.
  %
  %
\end{abstract}\blfootnote{\datestamp}



\section{Introduction}

\subsubsection*{\emph{Witten's asymptotic expansion conjecture}}

Let~$Y$ be a closed oriented $3$-manifold.
%
%
For $k\in\Z_{\ge2}$, let $\WRTk(Y) \in \C$ 
%
%
denote the level-$(k-2)$ Witten-Reshetikhin-Turaev invariant of~$Y$
constructed by Reshetikhin and Turaev in \cite{RT91,RT90} and
motivated by Witten's study \cite{Witten98} of quantum Chern-Simons
field theory with gauge group $\SU(2)$ and the Jones polynomial
\cite{Jones87,Jones85}. We work with the normalization given by $\WRTk(S^3)=1$.

Classical Chern-Simons theory \cite{ChernSimons74} is a gauge theory
with a Lagrangian formulation \cite{Freed}, which we now
present. Recall that every principal $\SU(2)$-bundle on~$Y$ is
trivializable. The action of the gauge equivalence class of an
$\SU(2)$-connection $A \in \Omega^1(Y,\mathfrak{su}(2))$ 
on the trivial principal $\SU(2)$-bundle is given by
\begin{equation}   \label{eqdefaction}
    \cS([A])= \frac{1}{8\pi^2} \int_Y \Tr(A \wedge dA+\frac{2}{3}A \wedge A \wedge A) \in \mathbb{R} / \Z.
\end{equation}
The space of solutions to the Euler-Lagrange equation $\delta\!\cS=0$
is equal to the moduli space $\M(Y)$ of flat $\SU(2)$-connections, and
we write $\CS(Y):=\cS(\M(Y))$.
The moduli space $\M(Y)$ is compact and the set $\CS(Y)$ is
finite. Viewing $\WRTk(Y)$ as the mathematical formalization of the
partition function of quantum Chern-Simons theory \cite{Witten98}
 motivates the following

\begin{conjecture}[The asymptotic expansion conjecture
  \cite{Andersen02,Andersen13,Witten98}]
  \label{ConjAEC}
  Let~$Y$ be a closed oriented $3$-manifold. For each
  $\SU(2)$ Chern-Simons action
  $S \in \CS(Y)$, there exists a Puiseux series
  %
  %
$W_{S}(\tau) \in \bigcup_{n=1}^{\infty} \C((\tau^{\frac{1}{n}}))$ such
that the WRT invariant of~$Y$ has the following Poincar\'e asymptotic
expansion 
  \begin{equation} \label{eq:AEC}
  \WRTk(Y) \sim  \sum_{S\in \CS(Y)} e^{2\pi i k S}\, W_S(k\ii)
  \quad \text{as $k\to\infty$.}
  \end{equation}
  %
  %
  %
  %
  %
  %
  %
\end{conjecture}

This conjecture is one of the central open problems in quantum
topology. It was discussed in \cite{Witten98}, from the point of view
of path integrals and perturbation theory. The above formulation is
independent of path integral techniques, and if true,
%
%
the collection of Puiseux series $(W_S)_{S \in \CS(Y)}$
will be uniquely determined by the asymptotic behaviour of $\WRTk(Y)$, and will therefore be a
topological invariant of~$Y$. 
The asymptotic expansion conjecture is connected to the use of
resurgence \cite{Ecalle81,Ecalle81b,MS16}. In recent years, there has
been a fruitful interplay between quantum topology, complexification,
asymptotic theory and resurgence, resulting in a large body of works
including
\cite{AM22,AM24,AP19,OG11,FantiniWheeler,G08,GGM23,GGM21,GK21,GL08,GLM08,GZ23,GMP2016resurgence,LZ99,Marino14,MistegaardMurakami,Witten11}.
This article is a contribution to this interplay.

%
%
We highlight that, complementary to Conjecture~\ref{ConjAEC}, there
are also the so-called growth rate conjecture
\cite[Conjecture~1.2]{Andersen13}, which gives an explicit conjecture
for the order of the leading terms of the expansion~\eqref{eq:AEC},
and
Witten's semi-classical approximation conjecture
\cite{Witten98} (see also \cite[Conjecture~1.3]{AH} and references in this paper),
which gives an explicit formula for the
coefficient of the leading terms of the expansion~\eqref{eq:AEC}. Both
of these conjectures are formulated in terms of gauge-theoretic
invariants.

 The asymptotic expansion conjecture is connected to the theme of
 integrality in quantum topology. For a general integral homology
 sphere $Y$, the number-theoretic nature of the WRT invariants has
 been well studied. It is known that $\WRTk(Y)\in \Z[e^{2\pi i/k}]$
 for all~$k$ (\cite{MasbaumRoberts1997,Murakami94,Habiro08}). Let
 $\gR \subset \C^{\times}$ denote the group of roots of unity.
By \cite{RT91} and \cite{Habiro08},
%
%
there exists a topological invariant in the form of a map
$\WRT(Y,\cdot): \gR \rightarrow \C$, such that
\begin{equation}   \label{eqWRTkze}
  %
  %
  %
\WRT(Y,\sigma\cdot e^{2\pi i/k})=\sigma\cdot\WRTk(Y)
  \quad \text{for every} \ens k\in\Z_{\ge1}
  \;\text{and}\; \sigma \in \Gal(\Q(e^{2\pi i/k}): \Q) 
\end{equation}
%
%
with the convention $\WRTone(Y)=1$.
%
%
Number-theoretic considerations have allowed Ohtsuki to extract a
formal power series invariant (\cite{Ohtsuki})
\begla
  \la_Y(q) = 1 + \sum_{n\ge1} \la_{Y,n} (q-1)^n,
\edla
the coefficients of which are now known to be integers (\cite{R06}),
the first nontrivial one being $\la_{Y,1} = 6 \la$, where $\la\in\Z$
is the Casson invariant of~$Y$ (\cite{Murakami94}).
%

\subsubsection*{\emph{Statement for Seifert fibered homology spheres}}

All the results in this paper are relative to the case
 where~$Y$ is a Seifert fibered integral homology sphere.
 Let $r\in \Z$ with $r\geq 3$. For each $j\in \{1,\ldots,r\}$, let
 $p_j,q_j$ be pairwise coprime non-zero integers, such that 
 $p_1,\ldots,p_r$ are positive and pairwise coprime and
 %
\beglab{eqdefP}
 P\sum_{j=1}^r \frac{q_j}{p_j}=1, \qquad \text{where}\ens P := p_1\cdots p_r.
 \edla
 Without loss of generality we assume that
$p_2,\ldots,p_r$ are odd.
 %
 %
 Let
 \begin{equation}   \label{defX}
   \begin{split}
   X := & \; \text{the closed
     oriented Seifert fibered $3$-manifold} \\
     & \; \text{with Seifert invariants
       $\{0;(p_1/q_1),\ldots,(p_r/q_r)\}$,}
   \end{split}
\end{equation}
where we follow the convention for Seifert invariants introduced in \cite{JN83}.
The $3$-manifold~$X$ is an integral homology sphere.


In the case of the Seifert integral homology sphere~$X$,
a large~$k$ expansion of $\WRTk(X)$ of the form~\eqref{eq:AEC} where
one is allowed to sum over a finite set of rationals~$S$ was
proven in \cite{LR99}.
Note that $0\in\CS(X) \subset \Q/\Z$ in this case;
%
%
it follows from the arguments in that article that the trivial
connection contribution~$W_0$ is a normalization of the aforementioned
Ohtsuki series: 
%
%
\beglab{eqdefWzerolaX}
W_0(\tau) = \la_X(e^{2\pi i\tau})  \in \Q[[2\pi i\tau]],
\edla
and that, for each non-zero~$S$, 
$W_S(k\ii)$ is a formal Laurent series in~$k^{-1/2}$ (i.e.\ the sum of a polynomial in $k^{1/2}$ and a
  formal series with non-negative integer powers of~$k^{-1/2}$).


Conjecture~\ref{ConjAEC} was proven for~$X$ in the case of $r=3$ in
\cite{K05a,LZ99} and for $r=4$ in the works \cite{K05b,K06}.
Our main result is that Conjecture~\ref{ConjAEC} holds for any Seifert
fibered integral homology sphere. More precisely:
 
 \begin{theorem} \label{Thm:main}
   For the Seifert integral homology sphere~$X$, the formal
   series~$W_0(\tau)$ of~\eqref{eqdefWzerolaX} is
   resurgent\footnote{Throughout this paper,   \label{ftnDefResur}
     we say that a formal series
  $\wt\Theta(\tau)=\sum_{p\ge0} a_p\tau^p$ is resurgent
  in~$\tau$ if its formal Borel transform
  $\wh\Theta(\xi):=\sum_{p\ge0} a_{p+1}\xi^p/p!$
  is convergent for~$|\xi|$ small enough and has ``endless analytic
  continuation'' with respect to~$\xi$;
  see \cite{Ecalle81,Eca93} or \cite{MS16,Encycl25}
  and beware that we slightly depart from the standard terminology,
  for which the above $\wt\Theta(\tau)$ would rather be
  considered resurgent in~$1/\tau$.
  }
and Borel-summable\footnote{Given a formal series   \label{ftnBLaplSum}
  of the same form as in footnote~\ref{ftnDefResur}, its Borel sum in
  a direction~$\tht$ is
  $\bS^\tht\wt\Theta(\tau) := a_0+\cL^\tht\wh\Theta(\tau)$ for
  $\arg\tau\in(\tht-\frac\pi2,\tht+\frac\pi2)$,
  with the Laplace transform operator~$\cL^\tht$ defined
  by~\eqref{eqdefcLtht} and 
  under suitable
  conditions (in particular $\wh\Theta(\xi)$ is supposed to be
  convergent for~$|\xi|$ small enough with analytic continuation along
  the ray $\R_{>0}e^{i\tht}$),
  and this function has Poincar\'e asymptotic expansion
  $\bS^\tht\wt\Theta(\tau) \underset{\tau\to0}{\sim}\wt\Theta(\tau)$
  with $1$-Gevrey qualification.
  }
in the directions of $(-\frac{3\pi}2,\frac\pi
2)$,
and there is an exact formula
   %
   %
   %
   \begin{gather} \label{eq:exactAEC}
    \WRTk(X) = (\bS^0 W_0)(k\ii) +
    \sum_{S \in \CS(X) \setminus \{0\}} e^{2\pi ikS}\,
   k^{3/2} \cE(k\ii) H_S(k)  
\\[-1.75ex] \intertext{with a convergent series}
 \label{eqdefcE}
  \cE(\tau) := \frac{ (-1)^r \tau\, e^{-i\pi\tau\phi/2} }{ 4 i \sin(\pi\tau) }
   = \frac{ (-1)^r }{ 4 \pi i } + O(\tau) \in \C\{\tau\}
  \quad \text{and} \ens
  \phi \; \text{as in~\eqref{eqdefphi} below,} 
\end{gather}
   where the $\SU(2)$ Chern-Simons actions $S\in\CS(X)$ are
     described in~\eqref{eqdefSl}--\eqref{eqformulCSX}
   below and, for each non-zero $S\in \CS(X)$, $H_S(k)$ is
   a polynomial in~$k$ satisfying
   \begin{equation}   \label{eqWSkZskdS}
     \deg(H_S) \le \tfrac{d_S}2 \ens
     \text{with $d_S :=$ maximum of the dimensions of the components of $\cS\ii(S)$.}
     %
     %
   \end{equation}
   %
   %
\end{theorem}

More will be said on the resurgent structure of~$W_0(\tau)$ later.
In view of the properties of the Borel-Laplace
summation operator in the direction $\tht=0$,
formula~\eqref{eq:exactAEC} implies the asymptotic
expansion~\eqref{eq:AEC} with $1$-Gevrey qualification,
%
 with $W_S(\tau) := \tau^{-3/2} \cE(\tau) H_S(\tau\ii)$
for $S\neq0$.

%
%
%

Comparing with \cite{LR99}, our main contribution is to show that, for all $S \in \Q/\Z$ with
  non-zero~$W_S$, we have $S \in \CS(X)$.
%
%
 Moreover, our bound~\eqref{eqWSkZskdS} on the degree of the
 polynomial~$H_S(k)$
in terms of the dimension of
the preimage of~$S$ in $\M(X)$ by the action functional~$\cS$
of~\eqref{eqdefaction}
%
%
is in agreement with the growth rate conjecture
\cite[Conjecture~1.2]{Andersen13}.
%
Further, our identification of~$H_S(k)$ in
Section~\ref{secpfthmmain} below 
together with other results
of this paper provides a first step towards 
proving the semi-classical approximation conjecture
\cite[Conjecture~1.3]{AH} for~$X$, as will be explained in the next
paragraph.
We emphasize that Hikami has proven results in this direction in the
case of $r=3$, but this case is much easier than for large~$r$. This
is because for $r=3$, the gauge-theoretic invariants appearing in
\cite[Conjecture~1.3]{AH} are defined with reference to discrete
moduli spaces of flat connections in this case, but in general, the
relevant moduli spaces have components of dimension up to as high as
$2r-6$.
%
%

Our proof of Theorem~\ref{Thm:main} depends on our
Theorems~\ref{thm:QM} and~\ref{thmMflat}, both of which are of
independent interest.
Theorem~\ref{thm:QM} demonstrates that the WRT invariant of~$X$ at a
general root of unity is,
up to an elementary factor,
%
%
the limit of the GPPV invariant of~$X$ \cite{GPPV} (introduced below)
at that root of unity;
it also demonstrates the quantum modularity of the GPPV invariant.
%
Theorem~\ref{thmMflat} gives a new parametrization of $\M(X)$
in terms of moduli spaces of flat $\SU(2)$-connections on the orbifold
surface of~$X$, with prescribed holonomy at exceptional orbits. This
is used in Corollary~\ref{Cor:CS} to determine the set of classical
Chern-Simons invariants $\CS(X)$.
%
%
Theorem~\ref{thmMflat} is also a first step towards proving the
semi-classical approximation conjecture \cite[Conjecture~1.3]{AH}. The
latter expresses the leading term coefficient of the
expansion~\eqref{eq:AEC} as an integral of gauge-theoretic functions
over components of the moduli space of flat connections on~$X$. In
general, this integration is difficult. However, our
Theorem~\ref{thmMflat} allows us to pull back these integrals to
smooth and compact moduli spaces of flat connections on punctured
spheres with prescribed holonomomy around the punctures. The advantage
is that cohomology generators and intersection pairings for these
moduli spaces have been thoroughly studied in the literature
\cite{JeffreyKirwan,thaddeusConformalFieldTheory1992,Witten92,zagierCohomologyModuliSpaces1995,meinrenkenWittenFormulasIntersection2005}.
The results are referred to as Witten's formulas for intersection
pairings and these are understood in sufficient generality for our
purposes. A proof of the semi-classical approximation conjecture
\cite[Conjecture~1.3]{AH} for~$X$ using Theorem~\ref{thmMflat} and
Witten's formulas for intersection pairings is planned to appear in a
separate publication.
%

\subsubsection*{\emph{The GPPV invariant of Seifert fibered homology spheres}}

We now present Theorem~\ref{thm:QM}. Being a Seifert fibered
$3$-manifold,~$X$ is also a graph $3$-manifold \cite{Waldhausen67}
and, as detailed in \cite{GM,AM22}, it admits a negative definite
plumbing graph (this notion is recalled in Section~\ref{sec:defGPPV}).
%
%
Consider the GPPV invariant
%
%
$\hh{Z}_0(X;q) \in q^{-\De_X} \Z[[q]]$,
%
where~$\De_X$ is the rational number defined 
by~\eqref{eqdefDeX} below. 
The GPPV invariant was introduced in \cite{GPPV} for pairs consisting of a
$3$-manifold with a $\spin^c$-structure, by use of physics arguments,
and it was proven in \cite{GM} to be a topological invariant of a
graph $3$-manifold with negative definite plumbing graph
and equipped with a $\spin^c$ structure (since our~$X$ is an integral
homology sphere, there is only one such). {For more on GPPV
  invariants, see Section~\ref{sec:defGPPV}.}

In our case, the coefficients of the \emph{the normalized GPPV invariant of~$X$}
\begin{equation}   
 Z^*(q) := q^{\De_X}\hh{Z}_0(X;q) \in \Z[[q]] 
\end{equation}
can be obtained as follows.
Define $m_0 \in \Z$ and the sequence of integers
$\big(\tilde\chi(m)\big)_{m=m_0}^{\infty}$ by the Laurent expansion
\begin{equation} \label{eq:Taylor}
  G(z) := (z^P-z^{-P})^{-(r-2)}\, \prod_{j=1}^r (z^{P/p_j}-z^{-P/p_j})
  = (-1)^r\sum_{m=m_0}^{\infty} \tilde\chi(m) z^m,
\end{equation}
where we use the notation~\eqref{eqdefP} and~$z$ is a new indeterminate.
One readily checks that
\beglab{eqformulamzero}
m_0 = \Big(r-2 - \sum_{1\le j\le r} \frac1{p_j}\Big) \, P
\edla
and $(-1)^r\ti\chi(m_0) = 1$.
By \cite[Theorem~3]{AM22},\footnote{%
  In \cite{AM22}, the quantity $\De_X + \frac{m_0^2}{4P}$
  %
  %
  is denoted by~$\De$ and computed in \cite[(4.2)]{AM22}.}
we have
%
%
%
%
\begin{equation}   \label{eqhatZz}  
  %
  %
  Z^*(q) =  
  %
  %
  %
  %
  \sum^{\infty}_{m=m_0}  \tilde\chi(m) q^{\frac{m^2-m_0^2}{4P}} 
\end{equation}
(it is a fact that $4P$ divides $m^2-m_0^2$ for all~$m$ in the support
of~$\tilde\chi$---see Proposition~\ref{propsuppchijmz} below).
%
%
%
%
This series is convergent for~$q$ in the open unit disc~$\DD$,
or equivalently for
\[
  q=e^{2\pi i \tau} \quad \text{with} \ens \tau\in\HH, \qquad
\HH := \{\tau \in \C \mid \IM(\tau)>0\}.
\]
%
%
We can thus define the normalized GPPV invariant of~$X$ as the
holomorphic function~$Z^*$ obtained as sum of~\eqref{eqhatZz} for $|q|<1$ or, equivalently,
\begin{equation}   \label{partialthetaseriesPsist}
  %
  %
 \Psi^*(\tau)=Z^*(e^{2\pi i\tau}) =
%
%
  %
  \sum_{m\ge m_0} \tilde\chi(m) e^{\frac{i\pi (m^2-m_0^2)\tau}{2P}},
  %
%
  \quad
  \text{$1$-periodic function of}\;\, \tau\in\HH. 
\end{equation}
%
%

In \cite{han2022resurgence} quantum modularity properties are analyzed
for partial theta series with coefficients given by a periodic
sequence multiplied by a monomial (we recall the definition of quantum
modularity in Section~\ref{secqmf}). 
Below we will see that \emph{the modified GPPV invariant}
\beglab{partialthetaseries}
 \Psi(\tau) := e^{\frac{i\pi m_0^2\tau}{2P}}\Psi^*(\tau) 
= \sum_{m\ge m_0} \tilde\chi(m) e^{\frac{i\pi m^2\tau}{2P}}
\edla
is a linear combination of functions of this form
(beware that it is not $1$-periodic in~$\tau$).
In this article, we apply the techniques from \cite{han2022resurgence}
to prove the following generalization of \cite[Theorem~4]{AM22}:


\begin{theorem} \label{thm:QM}
There is a family of formal series indexed by~$\gR$,
  %
%
%
%
%
\begin{equation} \label{eq:HatAlpha}
 \wt Z^*_{\ze}(q) := \sum_{m\ge0} Z^*_{\ze,m} (q-\ze)^m \in
  \C[[q-\ze]]   \qquad (\ze\in\gR), 
  %
  %
  %
\end{equation}
%
%
\vspace{-3.5ex}

\noindent such that $\wt Z^*_{\ze}(q)$ is resurgent in $q-\ze$ for
each~$\ze$ and:
\medskip

\noindent \emph{\textbf{(i)}}
The normalized GPPV invariant~$Z^*$ of~$X$
%
%
enjoys the asymptotic expansion property
\begin{equation} \label{eq:resurgentexpansionZst}
 Z^*(q) \sim \wt Z^*_\ze(q)  \quad \text{as}
  %
  %
  %
  \; q\to\ze \; \text{non-tangentially from within~$\DD$}
  %
  %
\end{equation}
for each $\ze\in\gR$.
In particular the constant term $Z^*_{\ze,0}$ is the non-tangential limit of $Z^*$ at~$\ze$. 
\medskip

\noindent \emph{\textbf{(ii)}}
%
%
The GPPV invariant and the WRT invariant are related as follows:
\begin{equation} \label{eq:Zhatradiallimit}
  Z^*_{\ze,0} = 2(-1)^r (\ze-1) \ze^{n_*-6\la} \WRT(X,\ze)
  \quad \text{for all $\ze \in \gR$},
  %
  %
\end{equation}
where~$\la$ is the Casson invariant of~$X$ and $n_* \in\Z$ is
defined by
\beglab{eqdefnstar}
-\frac{(r-1)(r-2)}2 P + (r-2) \sum_{1\le i \le r} \frac P{p_i}
- \sum_{1\le i < j \le r} \frac P{p_i p_j} = 2 n_* + 1
\edla
(it is a fact that the \lhs\ is an odd integer).
%

\noindent \emph{\textbf{(iii)}}
The family of formal series $(\wt\Psi_\al)_{\al\in\Q}$ defined from
the family $(\wt Z^*_\ze)_{\ze\in\gR}$ by the formula
\beglab{eqdefwtPsifromwtZst}
%
%
\wt\Psi_\al(\tau) := e^{\frac{i\pi m_0^2\tau}{2P}} \, \wt Z^*_\ze( e^{2\pi i\tau} )
%
%
%
= \sum_{m=0}^{\infty} {\Psi}_{\al,m} (\tau-\al)^m \in \C[[\tau-\al]]
%
%
%
\quad \text{with}\ens \ze = e^{2\pi i\al}
\edla
  satisfies the following properties:
  $\wt\Psi_\al(\tau)$ is resurgent in $\tau-\al$ for each $\al\in\Q$ and
  the modified GPPV invariant enjoys the asymptotic expansion property
\begin{equation} \label{eq:resurgentexpansion}
 \Psi(\tau) \sim \wt\Psi_\al(\tau) \quad \text{as}
 \; \tau\to\al \; \text{non-tangentially from within~$\HH$;}
\end{equation}
the function
\begla
\al\in\Q \mapsto
\Psi_{\al,0} = e^{\frac{i\pi m_0^2\al}{2P}} Z^*_{e^{2\pi i\al}, 0}
\edla
   %
   %
is a depth $r-2$ quantum modular form 
with weight $r-\frac52$ on the congruence
subgroup\footnote{\label{ftncongsubgp} With the
  standard notation
  $\Ga_1(N):=\{\, \begin{psmallmatrix} a&b \\
    c&d \end{psmallmatrix}\in\SL(2,\Z) \mid a=d=1\!\mod N,
  \; c=0\!\mod N\,\}$ for $N\ge1$.}~$\Ga_1(4P)$, and
it is a component
%
%
of a vector-valued depth $r-2$ quantum
modular form with weight $r-\frac52$ on the full modular group
$\SL(2,\Z)$;
   %
   the map $\al\in\Q \mapsto \wt\Psi_\al$ is a strong higher depth
   quantum modular form (with the same qualifications).
\end{theorem}

Theorem~\ref{thm:QM} will enable us 
to prove Theorem~\ref{Thm:main} by studying the GPPV invariant
and making use of~\eqref{eq:Zhatradiallimit} with $\ze=e^{2\pi i/k}$.
%
%
We will find that the series $W_0(\tau) = \la_X(e^{2\pi i\tau})$
of~\eqref{eqdefWzerolaX} is directly obtained from the asymptotic
expansion of the modified GPPV invariant as $\tau\to0$ by
formula~\eqref{eqWzintermsofcEPsiz} below.
Furthermore, we remark
that the resurgence and quantum modularity structure of all the formal
series $\wt\Psi_\al$ and thus of~$W_0$ and all the formal series $\wt Z^*_\ze$ is
completely understood by the results obtained in
\cite{han2022resurgence}, as explained below.

%
For $\alpha=(2k)^{-1}$, an asymptotic expansion of the form
\eqref{eq:resurgentexpansion} was obtained in \cite{AM22} and
for $\ze=e^{2\pi i/k}$
%
%
the identity \eqref{eq:Zhatradiallimit} was
conjectured in \cite{GPPV,GM}.
A similar result was obtained by different methods in \cite{FIMT21}
(for $\ze=e^{2\pi i/k}$)
%
%
and for $r=3$ in the work \cite{LZ99}.
The radial limit conjecture of \cite{GPPV,GM} was solved for general
plumbed $3$-manifolds with negative definite plumbing graph in
\cite{Murakami24}.
Further, we remark that quantum modularity for the modified GPPV invariant
of~$X$ was previously proven in the works \cite{BMM20a} by a different
method. In this article, we give a new proof of quantum modularity,
which uses resurgence to illuminate how quantum modularity is
connected to the ``Stokes phenomenon'', as explained in full detail below.

\begin{remark} \label{remstrongerf}
  The asymptotic property~\eqref{eq:resurgentexpansionZst} holds with
  $1$-Gevrey qualification for each $\ze=e^{2\pi i\al}\in\gR$, and
  similarly for~\eqref{eq:resurgentexpansion}.
This is a consequence of the following stronger facts:
\begin{quote}
\emph{\noindent The function $\Psi(\al+T)$ is the median Borel sum\footnote{With   \label{ftnLaplMedSum}
    reference to footnote~\ref{ftnBLaplSum}, in the present situation,
    we cannot use $\tht=\frac\pi2$ for $\wt\Psi_\al(\al+T)\in\C[[T]]$ due to
    the presence of singularities along $e^{i\frac\pi2}\R_{>0}$, but
    there are two well defined lateral Borel sum
    $\bS^{\frac\pi2\pm\epsilon}\wt\Psi_\al$
    independent of~$\epsilon$ small enough; their arithmetic average
    happens to coincide with the so-called ``median'' Borel sum in the
    direction~$\frac\pi2$ in this case (see \cite[Sec.~1.4]{Eca93},
    \cite{Menous}, \cite[p.~253]{han2022resurgence}),
     which we denote by $\bS^{\frac\pi2}\median\wt\Psi_\al$.}
in the direction~$\frac\pi2$ of the resurgent series $\wt\Psi_\al(\al+T)$,
  and the function $Z^*(\ze+Q)$ is the median Borel sum in the direction $2\pi\al+\pi$ of
  the resurgent series $\wt Z^*_\ze(\ze+Q)$.  }
\end{quote}
We will also see, in Section~\ref{secWRTmedsum}, that the WRT
  invariant at~$k$ is itself the limit of the median sum of the
  resurgent-summable series~$W_0(\tau)$ as $\tau\to 1/k$ non-tangentially from within~$\HH$.
\end{remark}

\smallskip

\begin{remark}
Given $\al\in\Q$ and $\ze=e^{2\pi i\al}$,
\begla
Z^*_{\ze,m}\in \Q(\ze), \quad
(2\pi i)^m\Psi_{\al,m} \in \Q(e^{2\pi i\al})
\quad \text{for all}\ens m\ge0.
\edla
For $\ze=1$ the constant terms vanish, $Z^*_{1,0} = \Psi_{\ell,0} =0$
for all $\ell\in\Z$, and the rational numbers $Z^*_{1,m}$ and $(2\pi
i)^m\Psi_{\ell,m}$ are related to the coefficients of the Ohtsuki
series $\la_X(q)$:
\beglab{eqZstonerelOhtsuki}
\wt Z^*_1(q) = 2(-1)^r (q-1) q^{n_*-6\la} \la_X(q),
\quad
\wt\Psi_0(\tau) = 2(-1)^r e^{2\pi i \big( \frac{m_0^2}{4P} + n_*-6\la \big) \tau}
( e^{2\pi i \tau} -1 ) \, \la_X( e^{2\pi i \tau} ).
\edla
 In view of~\eqref{eqdefWzerolaX} and~\eqref{eqmphipmzsqP} below,
  %
%
the above formulas are equivalent to
\beglab{eqWzintermsofcEPsiz}
 W_0(\tau) = \cE(\tau) \wt\Psi_0(\tau) / \tau
\quad\text{with}\;\, \cE(\tau)\;\, \text{as in~\eqref{eqdefcE}}
\edla
 (the formal series $\wt\Psi_0(\tau)$ is divisible by~$\tau$ since
  $\Psi_{0,0} =0$). 
\end{remark}


\subsubsection*{\emph{Plan of the article}}
In Section \ref{sec:quantuminvariants} we recall the definitions of
WRT invariants and GPPV invariants, as well as the definition of
a quantum modular form.

In Section \ref{sec:QMresurgence} we recall key elements from
\cite{han2022resurgence}. (Section \ref{sec:quantuminvariants} and
Section \ref{sec:QMresurgence} contain no new results, except for Proposition~\ref{propUnivPol}).

In Section \ref{sec:proofs} we first analyze the GPPV invariant in detail and describe it in terms of so-called Hikami functions, which play a central role. We then proceed to prove Theorem~\ref{thm:QM}.

Section~\ref{secWAECSFHS} is devoted to the proof of
  Theorem~\ref{Thm:main}.
In Section \ref{sec:Mflat} we parametrize the set of components of
$\M(X)$ and determine the set $\CS(X)$. The components of
the moduli space $\M^{\Irr}(X)$ of irreducible flat $\SU(2)$-connections
are shown to be homeomorphic to moduli spaces of flat
$\SU(2)$-connections on the orbifold surface of~$X$ with punctures inserted
at the exceptional orbits.
%
%
This is the content of Theorem~\ref{thmMflat}, which builds on the
works \cite{AM22,FS90,JM05,KK91}. We remark that
Theorem~\ref{thmMflat} is of independent interest, and that
Section~\ref{sec:Mflat} can be read independently of the rest of the
article.
Section~\ref{secpfthmmain} contains the proof of Theorem~\ref{Thm:main}.
%


 Appendix~\ref{Appendix:Normalizations} discusses normalization
  issues about the WRT invariants.
Appendix~\ref{sectiontoolbox} collects the technical computations that
are necessary to study the so-called generalized Hikami functions and their
discrete Fourier transforms; this is a class of periodic sequences,
some of which appear as elementary components in a decomposition of
the sequence~$\ti\chi$ of~\eqref{eq:Taylor}.

\subsubsection*{\emph{Announcement about the Habiro invariant}}

We conclude this introduction by seizing the opportunity for
announcing a new result about the Habiro invariant \cite{Habiro08} of Seifert fibered
homology spheres, that is closely related to our work on the GPPV invariant:
%
%
\begin{quote} \noindent \emph{At each roof of unity, the asymptotic
  expansion of the normalized GPPV invariant coincides with the Taylor
  expansion of the Habiro invariant (itself suitably normalized),
  which implies integrality of the former and resurgence-summabililty of the latter.}
\end{quote}

The precise statement (including a presentation of the
  normalizations) is given in Theorem~\ref{thm:HabiroGPPV} below. We
  now give the details. For $n\geq 0,$ let
$(q)_n:=\prod_{j=1}^n(1-q^j)$ be the $n$th Pochhammer symbol and
consider the Habiro ring $\wh{\Z[q]}:=\varprojlim \Z[q]/((q)_n)$
introduced in \cite{Hab04,Habiro08}.
It is easily seen that, for each root of unity~$\ze$, there is a
natural ring homomorphism
\begla
T_\ze \col \wh{\Z[q]} \to \Z[\ze][[q-\ze]],
\edla
which is proved to be injective in \cite{Hab04}.
For $\bJ\in \wh{\Z[q]}$, the formal series $T_\ze\bJ$ can be viewed as
the Taylor expansion at~$\ze$ of~$\bJ$ and its constant term $\ev_\ze\bJ$ as
the evaluation at~$\ze$ of~$\bJ$.
Collecting the constant terms by defining
$(\ev\bJ)(\ze):=\ev_\ze\bJ$,
we get a ring homomorphism $\bJ \mapsto \ev\bJ$
from~$\wh{\Z[q]}$ to the ring of functions on~$\gR$, which happens to
be injective too, by \cite{Hab04}.
These injectivity properties are like ``arithmetic quasianalyticity''
results, leading us to view the elements of~$\wh{\Z[q]}$ as ``analytic
functions on the space of roots of unity''.

In \cite{Habiro08}, K.~Habiro defined for every integral homology
three-sphere~$Y$ a topological invariant $\bJ_Y \in \wh{\Z[q]}$, now
called the Habiro invariant, which unifies the WRT invariants of~$Y$
in the sense that 
\begin{equation} \label{eq:HabiroEvaluatestoWRT}
  \ev_\zeta (\bJ_Y)=\WRT(Y,\zeta)
  \quad \text{for each}\ens \ze\in\gR.
\end{equation}
This simultaneously provided a unification of the WRT invariants at
different roots and generalized the integrality results of
\cite{MasbaumRoberts1997,Murakami94}
available for~$\ze$ of odd prime order.
%
%
Further, the Habiro invariant also dominates the Ohtsuki series in the sense that
\begla
T_1 \bJ_Y=\lambda_Y(q).
\edla
However, Habiro posed the challenge of interpreting the invariant
$\bJ_Y$ from the point of view of quantum Chern-Simons theory and to
elucidate its analytic properties.

We propose a solution to this in the form of
Theorem~\ref{thm:HabiroGPPV}, according to which the Taylor expansion
of the Habiro invariant at each~$\ze$ is equal to the asymptotic expansion of the
GPPV invariant suitably normalized.
This provides a physical explanation, as the GPPV invariant is a
nonperturbative mathematical model of the partition function of
quantum Chern-Simons theory with complex gauge group $\SL(2,\C)$, and
it explains the analytic properties as arising from the fact that the
collection of Taylor series is the collection of resurgent expansions
of a quantum modular form.
 
\begin{theorem} \label{thm:HabiroGPPV}
Formula~\eqref{eq:Zhatradiallimit} can be upgraded to
  \beglab{eqconjwtZstHabUnif}
  \wt Z^*_\ze = T_\ze\big(2(-1)^r (q-1) q^{n_*-6\la} \bJ_X\big)
  \quad \text{for each} \ens \ze\in\gR.
  \edla
\end{theorem}
  
We thus may consider the holomorphic function
\begla
q \in \DD^* \mapsto \JX(q) :=
\frac{Z^*(q)}{2(-1)^r (q-1) q^{n_*-6\la}}
\edla
as the ``analytic incarnation'' of the  Habiro
invariant~$\bJ_X$ in the sense that not only~$\JX$ has limits at the
roots of unity that match the evaluation of~$\bJ_X$, but also the
various expansions
$T_\ze\bJ_X \in \Z[\ze][[q-\ze]] \subset \C[[q-\ze]]$ are
resurgent series admitting median summation,
each of them producing the same function, namely $\JX(q)$.

Note that, since $T_1\bJ_X$ is nothing but the Ohtsuki
series~$\la_X(q)$, formula~\eqref{eqZstonerelOhtsuki} already says
that~\eqref{eqconjwtZstHabUnif} holds true for $\ze=1$.

The proof of Theorem~\ref{thm:HabiroGPPV}
will appear in a separate publication.






\section{Definitions: quantum invariants and quantum modular forms} \label{sec:quantuminvariants}

We briefly recall the definitions of the relevant quantum
invariants, 
first the WRT-invariants $\WRTk(Y)$
and then the GPPV invariants $\hh{Z}_a(Y;q)$.   

\subsection{WRT invariants}

Let $k\in\Z_{\ge2}$, $\Lambda_k:=\{1,\ldots,k-1\}$, $\zeta_k:=\exp(2\pi i/k)$.
For each $m\in \Lambda_k$, define the quantum integer
$[m]_k :=\sin(\pi m/k)/\sin(\pi/k)$.
For an oriented framed link $L \subset S^3$ with a labelling $\lambda \in \Lambda_k^{\pi_0(L)}$, we denote by
$J_{\lambda}(L,\zeta_k) \in \Z[\zeta_k^{\pm 1/4}]$ the colored Jones
polynomial of $(L,\lambda)$ evaluated at $\zeta_k$.
Originally defined by Jones \cite{Jones85, Jones87} using von
  Neumann algebras, this invariant can be defined in an elementary
  fashion using the Kauffman bracket polynomial \cite{Kauffman87, KauffmanLins94}.
Our normalization is such that for all $n \in \Z$ and $m\in\Lambda_k$, we have that 
\begin{equation} J_m(U_n,\zeta_k)=\zeta_k^{\frac{n(m^2-1)}{4}}[m]_k,
\end{equation} where $U_n$ is the $n$-framed unknot.
%
%
For $\epsilon \in \{-1,1\}$,
we define $G_{k,\epsilon}:=\sum_{m\in\Lambda_k} [m]_k J_m(U_{\epsilon},\zeta_k)$,
which is nonzero, as can be seen from explicit formulas in
terms of Gauss sums, and
\begin{equation}   \label{eqGkzero}
  G_{k,0}:=\frac{i\sqrt{2k}}{\zeta_k^{1/2}-\zeta_k^{-1/2}}
  = G_{k,0}\ii \sum_{m\in\Lambda_k} [m]_k J_m(U_0,\zeta_k).
\end{equation}
By \cite{Lickorish,Wallace}, every closed oriented $3$-manifold~$Y$
can be obtained by Dehn surgery on a framed oriented link $L\subset
S^3$, which is unique up to Kirby equivalence \cite{Kirby};
we then use the notation $Y=S^3_L$.
The notion of Dehn surgery is explained in detail in
  Section~\ref{sec:RationalSurgery} below.
Let $n_{\pm}(L)$ denote the number of positive/negative eigenvalues of
the linking matrix of $L$, and let $n_0(L)=b_1(S^3_L)=\rank(H_1(S^3_L,\Z))$.

\begin{definition}[\cite{RT91,RT90}]
\label{ddef:WRT}
  The $\SU(2)$ level-$(k-2)$ WRT invariant of $S_L^3$ is by definition
  \begin{equation} \label{def:WRT}
    \WRTk(S_L^3) := G_{k,0}^{-n_0(L)}G_{k,+}^{-n_+(L)}G_{k,-}^{-n_{-}(L)}\sum_{\lambda
      \in \Lambda_k^{\pi_0(L)}} J_{\lambda}(L,\zeta_k) \prod_{j \in \pi_0(L)} [\lambda_j]_k.
    %
    %
  \end{equation}
\end{definition}
It was proven by Reshekthin and Turaev \cite{RT91,RT90} that the
complex number on the right hand side of \eqref{def:WRT} is an
invariant of the Kirby equivalence class of $L$, i.e.\ the set of all
links $L'$ which can be obtained from $L$ by a finite sequence of
Kirby moves, and therefore defines an invariant of the $3$-manifold
$S_L^3$ (in \cite{RT91}, they in fact worked with a slightly different
normalization as detailed in Appendix~\ref{Appendix:Normalizations}).
With the above normalization, we have
\[
\WRTk(S^3)=1, \qquad
\WRTk(S^1 \times S^2)=G_{k,0} = \frac{i\sqrt{2k}}{\exp(\pi i/k)-\exp(-\pi i/k)},
\]
since $S^3=S^3_{U_{\pm 1}}$ and $b_1(S^3)=0$, and $S^1\times S^2=S^3_{U_0}$ and $b_1(S^1 \times S^2)=1$.

\subsubsection{Integrality} If $S_L^3$ is an integral homology sphere (i.e.\  if
$n_0(L)=b_1(S_L^3)=0$), then the invariant $\WRTk(S_L^3)$ is equal to
the invariant denoted by $\tau_{\zeta_k}(S_L^3)$ in
\cite{Habiro08},
%
%
and for such a $3$-manifold, we have by
\cite{Habiro08} that $\WRTk(S_L^3) \in \Z[\zeta_k]$. For every
primitive $k$th root of unity $\ze$, there exists a unique Galois
transformation $\sigma \in \Gal(\Q(\zeta_k): \Q)$ such that
$\sigma\cdot \zeta_k=\ze$, and
%
%
\begin{equation} \label{eq:GaloisEquivariance}
    \WRT(S_L^3,\ze)=\sigma\cdot \WRTk(S_L^3) \in \Z[\ze_k]
  \end{equation}
  by~\eqref{eqWRTkze}.

\subsubsection{A formula for WRT invariants in terms of rational surgery presentations} \label{sec:RationalSurgery}

Let $L \subset S^3$ be a framed oriented link. Let
$\{L_j \}_{j\in \{1,\ldots,m\}}$ be the set of components of $L$. For
each $j \in \{1,..,m\}$ the framing determines an orientation
preserving diffeomorphism
$(\nu(L_j),L_j) \cong ( B^2\times S^1,\{0\} \times S^1)$, where
$\nu(L_j) $ is a tubular neighbourhood of $L_j$ and $B^2 \subset \R^2$
is the unit disc. For each $j\in \{1,\ldots,m\}$, let $a_j,b_j \in \Z$ be
co-prime integers and let $B_j \in \SL(2,\Z)$ be a matrix such that
the first column of $B_j$ is equal to the tranpose of $(a_j,b_j)$. Let
$B=(B_j)_{j \in \{1,\ldots,m\}}$. Recall that each $B_j$ acts by an
orientation-preserving diffeomorphism on $S^1 \times S^1$ through the
identification $S^1 \times S^1 \cong (\R / \Z)^2$. Set
$\nu(L)=\bigsqcup_{j=1}^m \nu(L_j)$. The $3$-manifold
$S^3_{L,B}$ obtained through surgery on $L$ with rational surgery data
$B \in \SL(2,\Z)^{\pi_0(L)}$ is given by the quotient space
\[
S^3_{L,B} := \left( S^3 \setminus \interior \nu(L)\right) \bigsqcup_{j=1}^m (B^2 \times S^1)_j / \sim,
\]
where the quotient is with respect to the equivalence relation
generated by the identifications
$B_j: \partial (B^2 \times S^1)_j \rightarrow \nu(L_j)$ for $j\in
\{1,\ldots,m\},$ through the usual identification
$\partial (B^2 \times S^1)=S^1 \times S^1=\partial \nu(L_j).$ The
class of $S^3_{L,B}$ as an oriented smooth manifold depends only on
the tuple $B$ through the tuple of rationals
$(a_j/b_j)_{j\in \{1,\ldots,m\}}$, and therefore the notation
$S^3_{L,(a_j/b_j)}$ is commonplace. Performing standard Dehn surgery
on a component $L_j$ corresponds to the assignment $a_j=0, b_j=1$.

In \cite{Jeffrey92} a formula is given for the WRT invariant of
$Y=S^3_{L,B}$ in terms of the colored Jones polynomial of $L$. To
state this formula, we need to recall a certain representation
$\rho_k\col \PSL(2,\Z) \rightarrow \GL(k-1,\C)$, which is known from the
study of affine Lie algebras \cite{Kac90}, and we need to recall the
Rademacher $\Phi$ function.
Recall that $\SL(2,\Z)$ can be generated by the two matrices $T:=\begin{psmallmatrix} 1&1 \\
  0&1 \end{psmallmatrix}$
and $S:=\begin{psmallmatrix} 0&-1 \\ 1&0 \end{psmallmatrix}$.
The representation $\rho_k$ is determined
by the following explicit formulas for the matrix entries, where $j,\ell$
range through $\Lambda_k=\{1,\ldots,k-1\}$,
\[
  \rho_k(S)_{j,\ell}=\sqrt{\frac{2}{k}}
  \sin \left(\frac{\pi j\ell}{k}\right), 
\qquad \rho_k(T)_{j,\ell}= e^{-\pi i /4}\zeta_{4k}^{j^2} \, \delta_{j,\ell}.
%
%
\]
For coprime integers $a,b$ 
we use the notation $s(a,b)$ for the Dedekind sum. For
$\ga = \begin{psmallmatrix} a & c \\ b & d
\end{psmallmatrix} \in \SL(2,\Z)$
the Rademacher function is given by
\[
\Phi(\ga) := \begin{cases} &\frac{a+d}{b}-12s(a,b) \qquad\text{if $b\not=0$}, 
\\ & \frac{c}{d} \qquad \text{otherwise.}
\end{cases}
\]
Finally, define
\[
\Phi(L,B):= \sum_{j=1}^m \Phi(B_j)-3 (n_+(L)-n_{-}(L)),
\]
where, as above, $n_\pm(L)$ denotes the number of positive/negative
eigenvalues of the linking matrix of $L$.
We then have 
\begin{equation} \label{eq:rationalsurgeryformula}
\WRTk(S^3_{L,B})= \exp\left({\frac{\pi i}{4}\left(\frac{k-2}{k}\right) \Phi(L,B)}\right) \sum_{\lambda \in \Lambda_k^{\pi_0(L)}} J_{\lambda}(L,\zeta_k) \prod_{j \in \pi_0(L)} \rho_k(B_j)_{\lambda_j,1}.
\end{equation}
This formula is generalized in \cite[Corollary $8.3$]{Hansen01} and
note that, to compare, one must take into consideration the difference
in normalization explained in Appendix~\ref{Appendix:Normalizations}.

\bblack

\subsection{GPPV invariants} \label{sec:defGPPV}
Let $(\Gamma,b)$ be a weighted tree, i.e.\ $\Gamma$ is a tree together
with a map~$b$ from its set of vertices~$V$ to~$\Z$. Let
$B=B(\Gamma,b)$ be the $V\times V$ symmetric matrix with entries given by
\[
B_{v,w} := \begin{cases}
  b(v) & \text{if} \ v=w,
  \\
  1 & \text{if}  \ v \text{ and }  w \ \text{are joined by an edge,}
  \\
  0 & \text{otherwise.} \end{cases}
\]
  We say~$B$ is weakly negative definite if~$B$ is invertible and
  $B^{-1}$ is negative definite on the subspace of $\Z^V$ spanned by
  vertices of degree at most~$3$. Further, we say that the graph
    $(\Gamma,b)$ is negative definite (resp. weakly negative definite) if
    the adjacency matrix $B$ is negative definite (resp. weakly
    negative definite). 
  Assume that $B$ is weakly negative
  definite. Let
  \begin{equation}
    Y \; := \; \text{the oriented closed $3$-manifold with surgery
      link $L(\Gamma,b)$}
  \end{equation}
  where $L=L(\Gamma,b)$ is constructed as follows:
  for each vertex $v$ the
  link $L$ has an unknotted component $U_v$ with framing $b_v,$ and
  the linking number of $U_w \cup U_v$ is equal to unity if $v$ and
  $w$ are joined by an edge, and otherwise $U_{w} \cup U_v$ is a
  split-link of two unknots. Notice that $B$ is the linking matrix of
  $L$.

  Assume that $b_1(Y)=0$.
  As above, let $n_+(B)$ denote the number of positive eigenvalues of~$B$. Let $\sigma(B)$ be the signature of $B$, and set 
  \begin{equation} \label{label:Delta}
    %
    %
    \Delta(B) :=\frac{3\sigma(B)-\sum_{v\in V} b(v)}{4}.
  \end{equation}
  Set $\delta=(\deg(v))_{v\in V} \in \Z^V$ and set
  $\vec{b}=(b(v))_{v \in V} \in \Z^V$ As explained in detail in
  \cite{GM} we have isomorphisms $ \spin^c(Y) \simeq H_1(Y,\Z)  \simeq
  (\Z^V+\vec{b})/ 2B\Z^V \simeq (\Z^V+\delta)/ 2B\Z^V$.
  Let $a \in  (\Z^V+\delta)/ 2B\Z^V$.
  Define the formal series
  \begin{equation} \Theta_a^{-B}(q,\vec{z}):=\sum_{\vec{l}\in 2B\Z^V+a}
    q^{-\frac{(\vec{l},B^{-1}\vec{l}))}{4}} \prod_{v \in V} z_v^{l_v}
    \in \Z[z_v, v \in V][[q]]
  \end{equation}
  where~$q$ and $(z_v)_{v\in V}$ are indeterminates.
  %

      \begin{definition}[\cite{GPPV}] The GPPV invariant of $(Y,a)$ is by definition
  	\begin{equation} \label{eq:GPPV}
          \hh{Z}_a(Y;q) := (-1)^{n_+(B)} q^{\Delta(B)}
          %
          %
        \,   v.p. \oint_{
          {\displaystyle\,\underset{\fakeheight{$\scriptstyle v\in
                V$}{$\scriptstyle \ti V$}}{\fakeheight{$\bigtimes$}{$=$}} \scriptstyle\{| z_v |=1\}} }
          \,\,
          \prod_{v\in V} \frac{d  z_v}{2\pi iz_v}\left(z_v-z_v^{-1}\right)^{2-\text{deg}(v)}
          \Theta_a^{-B}(q,\vec{z}), 
	\end{equation}
        where $v.p.$ denotes the principal value of the integral.
      \end{definition}
      
      The topological invariance of \eqref{eq:GPPV} was proven in
      \cite{GM}.

      As $X$ is an integral Seifert fibered $3$-manifold, $X$ is also
      a graph $3$-manifold \cite{Waldhausen67} and, as detailed in
      \cite{GM,AM22}, it admits a negative definite plumbing
      graph $(\Ga,b)$. We set
      \beglab{eqdefDeX}
       \De_X := \De\big( B(\Gamma,b) \big). 
      \edla
      Further, as $X$ is a Seifert fibered integral homology
      sphere, there is up to isomorphism only one $\spin^c$-structure,
      which we denote by
      $0$. 




  \subsection{Quantum modular forms with higher depth}  \label{secqmf}

  The study of modular forms boasts a rich historical
  background. Following a significant example by Kontsevich, Zagier
  laid down the groundwork for what are now termed quantum modular
  forms (cf.\ \cite{Z01}, \cite{Zagier10}). Additionally, Lawrence and
  Zagier delved into exploring the interplay between quantum modular
  forms and WRT invariants \cite{LZ99}. In this section, we will
  revisit the definition of quantum modular forms as delineated in
  \cite{BMM20b}. Our objective is to demonstrate that the GPPV
  invariant qualifies as a quantum modular form of higher depth, as
  stated in Theorem~\ref{thm:QM}(iii). 

  To fix our notations, we recall that the left action of $\GA:=\SL(2,\Z)$ on~$\HH$ 
\beglab{eqdefleftactionGA}
\ga = \matr{a&b}{c&d} \Imp \ga\tau = \frac{a\tau+b}{c\tau+d}
\edla
extends to $\tau \in \HH\cup\Q\cup\{\infty\}$, and 
\beglab{eqdefJga}
J_\ga(\tau) := c\tau+d
\ens\text{satisfies}\ens
J_{\ga_1 \ga_2} = (J_{\ga_1} \circ \ga_2) J_{\ga_2}
\ens\text{for all}\; \ga_1,\ga_2\in\GA.
\edla

\bblack

  
\begin{definition}[adapted from \cite{BMM20b,Zagier10}]\label{definitionquantummodularforms}
  Let $\sQ\subset \Q$, $w\in\frac12 \mathbb{Z}$ and let~$\Ga$ be a
  subgroup of $\SL(2,\Z)$ leaving~$\sQ\cup\{\infty\}$ invariant.
  Given a function $\eps\col\Ga\to \C^*$, we say that a function
  $\ph\col\sQ \to \mathbb{C}$ is a quantum modular form on~$\Ga$ with
  weight~$w$, quantum set~$\sQ$ and multiplier~$\eps$
if, for any $\ga = \begin{psmallmatrix}
      a&b  \\    c&d   \end{psmallmatrix}\in\Ga$, 
    the ``modularity defect''
    \begin{equation}   \label{eqdefmodpbstr}
      \al \in \sQ \setminus\{-d/c\} \mapsto
      \ph(\al) - \eps(\gamma) J_\ga(\al)^{-w} \ph(\ga\al)
      \end{equation}
      belongs to $\cO(R_\ga)$ for some open subset~$R_\ga$ of~$\R$
      (i.e.\ extends to a holomorphic 
function on~$R_\ga$).
      %
      %
      The vector space of such functions is denoted by
      $\cQ_w^1(\sQ,\Ga,\eps)$.
\end{definition}

If~$w$ is integer, then the term $J_\ga^{-w} \, (\ph\circ\ga)$ that
appears in~\eqref{eqdefmodpbstr} is unambiguously determined (and is
related to the ``weight~$w$ left action''
$(\ga,\phi)\mapsto J_\ga^{-w} \, (\phi\circ\ga)$
of $\SL(2,\Z)$ on the space
of functions on~$\HH$).
Since~$w$ may be non-integer, we must specify which branch of
$J_\ga(\al)^{-w}$ we use then; 
our convention will be determined by:
\beglab{eqchoosebranch}
c\al+d>0 \imp J_\ga(\al)^{1/2} \in \R_{>0},
\qquad
c\al+d<0 \imp J_\ga(\al)^{1/2} \in \;
\left|\begin{split}
   -i\R_{>0} & \ens\;\text{if $c>0$}\\[.8ex]
  i\R_{>0} & \ens\;\text{if $c\le0$}
\end{split}\right.
\bblack
\edla
(see Appendix~\ref{secVVsqmf} for a better point of view, relying on
the use of the \emph{metaplectic double cover of $\SL(2,\Z)$}).

\begin{remark}   \label{remgamga}
  One can check that the modularity defects~\eqref{eqdefmodpbstr}
  associated with~$\ga$ and~$-\ga$ coincide when
  $\eps(-\ga)=i^{2w}\eps(\ga)$, because our convention implies that
  $J_{-\ga}^{1/2} = i J_\ga^{1/2}$ if $c>0$, or if $c=0$ and $d>0$.
\end{remark}
  
\begin{remark}   \label{remqmfupper}
  In subsequent discussions, when this does not cause any ambiguity, we will
  sometimes speak of a function~$\phi$ defined on the upper
  half-plane as a quantum modular form. This means that~$\phi$ has limits
  at the points of the quantum set that provide the function~``$\ph$''
  of Definition~\ref{definitionquantummodularforms}.
\end{remark}

Quantum modular forms with depth~$N$ are a generalization of quantum
modular forms (which are declared to have depth~$1$):

\begin{definition}[adapted from \cite{BMM20b}]\label{definitionhigherdepthscaler}
  Given $\sQ$, $w$, $\Ga$ and~$\eps$ as above, 
  the space $\cQ_w^N(\sQ,\Ga,\eps)$ of quantum modular forms with depth~$N$ is
  inductively defined as follows:
  $\cQ_w^0(\sQ,\Ga,\eps)=\C$, $\cQ_w^1(\sQ,\Ga,\eps)$ is as in
  Definition~\ref{definitionquantummodularforms} and, for $N\ge2$,
  $\cQ_w^N(\sQ,\Ga,\eps)$ is the space of all functions $\ph \col\sQ\to \C$ 
  such that, for any $\ga\in\Ga$,
  %
    %
    the modularity defect
    \begin{equation}
      \ph(\al) - \eps(\gamma)J_\ga(\al)^{-w} \ph(\gamma\al)
      \ens \text{belongs to} \ens
      \bigoplus_{j=1}^J \cO(R_\ga) \otimes \cQ_{w_j}^{N_j}(\sQ,\Ga,\eps_j)
    \end{equation}
    where~$R_\ga$ is an open subset of~$\R$ and $J\in\Z_{\ge1}$,
    %
    %
    for some weights $w_1,\ldots,w_J \in \frac12 \mathbb{Z}$ 
    and multipliers $\eps_1,\ldots\eps_J$, 
    and with $0\le N_j <N$ for each~$j$.    
\end{definition}


Vector-valued quantum modular forms are defined as follows:

\begin{definition}\label{definitionhigherdepthvv}
  Given $\sQ$, $w$, $\Ga$ as above, $M\in\Z_{\ge1}$ and
  $\eps=[\eps_{m,\ell}] \col \Ga \to \GL(M,\C)$,
  we define $\ocQ_w^N(\sQ,\Ga,\eps)$ by induction on~$N$:
  $\ocQ_w^0(\sQ,\Ga,\eps) := \C^M$ and, for $N\ge1$,
  $\ocQ_w^N(\sQ,\Ga,\eps) :=$ the space of tuples $(\ph_1,\cdots,\ph_M)$
  of functions $\ph_\ell \col \sQ \to \mathbb{C}$ such that, for any
  $\ga\in\Ga$,
    %
       %
       \begin{equation}   \label{eqvvmoddef}
       \Big( \ph_\ell(\al) - J_\ga(\al)^{-w} \sum_{m=1}^M \eps_{m,\ell}(\ga) \ph_m(\ga\al) \Big)_{1\leq \ell \leq M} 
       \ens \text{belongs to} \ens
       \bigoplus_{j=1}^J \cO(R_\ga) \otimes
       \ocQ_{w_j}^{N_j}(\sQ,\Gamma,\eps^{(j)}),
     \end{equation}
    where~$R_\ga$ is an open subset of~$\R$ and $J\in\Z_{\ge1}$,
    %
    %
    for some weights $w_1,\ldots,w_J \in \frac12 \mathbb{Z}$ 
    and matrix-valued multipliers $\eps^{(1)},\ldots,\eps^{(J)}$, 
    and with $0\le N_j <N$ for each~$j$.    
     %
     %
     %
     %
     %
  \end{definition}

Finally, the ``strong'' version of quantum modular forms is obtained by following the
lines of \cite{Zagier10} and replacing functions $\ph\col\sQ\to\C$ with
maps
\[
\al\in\sQ \mapsto  \wt \ph_\al = \sum_{m\ge0} \ph_{\al,m} T^m \in \C[[T]].
\]
Heuristically, $\wt \ph_\al(T)$ stands for
``$\ph(\al+T)$'',
where ``$\ph$'' should be the strong quantum modular form, except that
the formal series $\wt \ph_\al$ maybe very well be divergent for
all~$\al$.
This is formalized in Definitions~\ref{defstrongqmf} and~\ref{defstrongqmfvar}:

\begin{definition}   \label{defstrongqmf}
Given $\sQ$, $w$, $\Ga$ and~$\eps$ as in Definition~\ref{definitionquantummodularforms},
we say that a family of power series
$(\wt \ph_\al)_{\al\in\sQ}$ is a strong quantum modular form
on~$\Ga$ with weight~$w$, quantum set~$\sQ$ and multiplier~$\eps$ if:
\smallskip

\noindent \textbf{(i)}
the constant terms give rise to a quantum modular form
$\al\in\sQ \mapsto \ph_{\al,0}$, belonging to
$\cQ_w^1(\sQ,\Ga,\eps)$, thus with modularity defects
\begin{equation} \label{eqdefmodpbstrstrong}
  h_\ga(\al) := \ph_{\al,0} -  \eps(\ga) (c\al+d)^{-w} \ph_{\ga\al,0}
\end{equation}
extending to holomorphic 
functions $h_\ga\in \cO(R_\ga)$
for all $\ga=
    \begin{psmallmatrix}
         a&b  \\   c&d 
       \end{psmallmatrix}\in \Ga$,
\smallskip

\noindent \textbf{(ii)}
for each $\al\in R_\ga\cap\sQ \setminus\{-d/c\}$ , the formal series
\[
  \wt h_{\ga,\al}(T) :=
  \wt \ph_\al(T)-  \eps(\ga) \big( c(\al+T)+d \big)^{-w}
  \wt \ph_{\ga\al}\big(\ga(\al+T)-\ga\al\big) \in \C[[T]]
\]
coincides with the Taylor series of $h_\ga(\al+T)$ around $T=0$.
\end{definition}
  
A condition equivalent to (i)--(ii) is that for each $\ga=
    \begin{psmallmatrix}
         a&b  \\   c&d 
       \end{psmallmatrix}\in\Ga$
       there exists an open subset~$R_\ga$ of~$\R$ such that,
       for each $\al\in R_\ga\cap\sQ\setminus\{-d/c\}$, the formal series
\begin{equation}
            \wt \ph_\al(\tau-\al) - \eps(\ga) J_\ga(\tau)^{-w}
            \wt \ph_{\ga\al}(\ga\tau-\ga\al) \in \C[[\tau-\al]]
\end{equation}
is convergent and is the Taylor series at~$\al$ of a holomorphic
function~$h_\ga$ that does not depend on~$\al$.

\begin{definition}   \label{defstrongqmfvar}
  Strong quantum modular forms with higher depth, possibly
  vector-valued, are defined from Definition~\ref{defstrongqmf} by
  mimicking the passage from Definition~\ref{definitionquantummodularforms}
  to Definitions~\ref{definitionhigherdepthscaler}--\ref{definitionhigherdepthvv}.
\end{definition}

\begin{remark}     \label{remsqmfupper}
  Similarly to Remark~\ref{remqmfupper}, we will sometimes speak of a
  function $\phi \col \HH\to\C$ as a strong quantum modular form.  This
  means that it has asymptotic expansions $\wt \ph_\al(\tau-\al)$ at all
  points~$\al$ of a quantum set $\sQ\subset\Q$ that satisfy
  Definition~\ref{defstrongqmf} or~\ref{defstrongqmfvar}.
\end{remark}

\begin{remark}
  Zagier's seminal paper also mentions an extra property (``leaking''
  into the lower half-plane though~$\sQ$) that is
  sometimes encountered in the setting of Remark~\ref{remsqmfupper}:
  it may be the case that the formal series $\wt \ph_\al$ making up the
  strong quantum modular form occur as asymptotic expansions of one
  function~$\phi$ in~$\HH$ and also as asymptotic expansions of one
  function~$\phi^-$ in $\HH^-:=\{\IM\tau<0\}$.
  This is the case for the partial theta series considered in
  \cite{han2022resurgence} and next section, as explained in
  \cite{transseriescompletions},
  with~$\phi$ holomorphic in~$\HH$ and~$\phi^-$ real analytic (not
  holomorphic!) in~$\HH^-$; this will imply a similar property for the
  modified 
  GPPV invariant~$\Psi(\tau)$.
\end{remark}


  \section{Reminders on partial theta series} \label{sec:QMresurgence}




Given a positive integer~$M$, an $M$-periodic function $f:\Z
\rightarrow \C$ and a non-negative integer~$j$,
%
%
the corresponding \emph{partial theta series} is the holomorphic function
\begin{equation}   \label{eqdefThetafjM}
  \Theta(\tau; j,f,M):= \sum_{n=1}^{\infty} n^j  f(n)
    e^{i\pi n^2\tau/M}
    \qquad \text{for} \ens \tau\in\HH.
\end{equation}
These functions are studied in \cite{han2022resurgence} from the
viewpoint of Borel-Laplace summation and resurgence, with a view to
describing their asymptotic behaviour as~$\tau$ tends to a rational
number~$\al$ and their modularity or quantum modularity properties.

It turns out that the modified
GPPV invariant~$\Psi(\tau)$ of~\eqref{partialthetaseries} can be
recast as a sum of partial theta series
(up to a trigonometric polynomial of~$\tau$ in the case of the
Poincar\'e homology sphere)---see Proposition~\ref{propPsiassumpts}
below.
We thus recall now the key elements of the analysis from \cite{han2022resurgence}.

\bblack

\subsection{Partial theta series and Laplace transforms}   \label{secptsLt}

In this section we review some results of \cite{han2022resurgence} about
  the partial theta series of the form $\Theta(\cdot\,; j,f,M)$ under
  the assumption that~$M$ is even (it will be $2P$ in our
application),
with emphasis on the case
  \begin{equation}   \label{eqparity}
    (\text{$j$ is even and $f$ is an odd function})
    \ens \text{or} \ens
    (\text{$j$ is odd and $f$ is an even function}).
  \end{equation}

Let $\alpha\in \Q$. The first result from \cite{han2022resurgence}
that we present is formula~\eqref{eq:Taunearalpha} below; it is useful
to understand the asymptotics of the function $\Theta(\tau; j,f,M)$ for $\tau$
near~$\alpha$ and will also be used for studying its quantum
modularity properties.
Consider the function
\[
  f_{\alpha/M} \col m\in \Z \rightarrow f_{\alpha/M}(m):=f(m)\exp(\pi
  im^2 \alpha/M).
\]
%
%
%
Clearly $f_{\alpha/M}(m)$ is periodic, and we let $M_{\alpha} \in
\Z_{\ge1}$ be a period.
For concreteness, one can take the least common multiple of $M$ and
the denominator $\den(\alpha/M)$ of $\alpha/M$,
%
but we stress that all of the formulas below are valid for any choice
of period, e.g.\ $M\den(\al)$.
Consider the generating function $F_{j,f_{\alpha/M}}$ of the sequence
$m \mapsto m^jf_{\alpha/M}(m)$ defined as follows:
\begin{equation}   \label{eq:generatingfunction}
  F_{j,f_{\alpha/M}}(t) := \sum_{m=1}^{\infty}
  m^jf_{\alpha/M}(m)\exp(-mt)
  =
  \Big(\!-\frac{d}{dt}\Big)^{\!j}\bigg(\frac{\sum_{\ell=1}^{M_{\alpha}}
    f_{\alpha/M}(\ell)\exp(-\ell t)}{1-\exp(-M_{\alpha}t)}\bigg).
\end{equation}
By the rightmost equality in~\eqref{eq:generatingfunction}, we see
that $F_{j,f_{\alpha/M}}$ has a meromorphic continuation to~$\C$ with
potential poles at $2\pi im/M_{\alpha}$, $m \in
\Z$, 
%
%
and its principal part at the origin is 
%
%
$j!\langle f_{\alpha/M} \rangle t^{-j-1}$,
where $\langle f_{\alpha/M} \rangle =\frac1{M_\al}\sum_{m=1}^{M_{\alpha}} 
f_{\alpha/M}(m)$ is the mean value of $f_{\alpha/M}$.
%
%
We can thus implicitly define holomorphic germs $\hh{\phi}^{\pm}_{j,f,\alpha,M}(t) \in \C\{t\}$ by
\begin{equation}
     F_{j,f_{\alpha/M}}(t)= \frac{j!\langle f_{\alpha/M} \rangle}{t^{j+1}}+\pi^{1/2}\hh{\phi}^{+}_{j,f,\alpha,M}(t^2/C_{M_\al}^2)+\pi^{1/2}\frac{t}{C_{M_\al}}\hh{\phi}^{-}_{j,f,\alpha,M}(t^2/C_{M_\al}^2),
   \end{equation}
   where $C_{M_\al}:=\sqrt{4\pi/M_\al} \, e^{i\pi /4}$.
The germs $\hh{\phi}^{\pm}_{j,f,\alpha,M}(\xi)$ extend to meromorphic
functions on~$\C$ with potential poles at
$\xi_m=i\pi m^2/M_{\alpha}$, $m \in \Z_{\ge1}$. Let
$\tau \in \mathbb{H}.$ For sufficiently small $\epsilon>0$, the
following Laplace tranforms are well-defined holomorphic functions of~$\tau$:
\begin{equation}   \label{eqThetapmLaplTrsf}
  \Theta_{j,f,\alpha,M}^{\pm}(\tau):= \frac{1}{2\tau^{1/2}}\Big( \cL^{\pi/2-\epsilon}\pm\cL^{\pi/2+\epsilon}\Big)\left[\frac{\hh{\phi}^{\pm}_{j,f,\alpha,M}(\xi)}{\xi^{1/4\pm 1/4}}\right](\tau),
\end{equation}
with the notation
\beglab{eqdefcLtht}
\cL^\tht\wh\ph(\tau) := \int_0^{e^{i\tht}\infty}
e^{-\xi/\tau} \wh\ph(\xi) \, d\xi
\quad \text{for} \ens \arg\tau\in(\tht-\pi/2,\tht+\pi/2).
\edla
By \cite[Remark~2.1]{han2022resurgence} we have that
$\Theta(\al+\tau; j,f,M)=\Theta(M_{\alpha}\tau/M;
j,f_{\alpha/M},M_{\alpha})$, and the desired fomula
follows from \cite[Theorem~1 \& eqns~(3.4)--(3.6)]{han2022resurgence}:
\begin{equation} \label{eq:Taunearalpha}
  \Theta(\al+\tau; j,f,M)=
  \frac{1}{2}\Gamma\Big(\frac{j+1}{2}\Big)\langle f_{\alpha/M}
  \rangle \Big( \frac{\pi \tau}{M i}\Big)^{-\frac{j+1}{2}}+
  \Theta_{j,f,\alpha,M}^{+}\Big(\frac{M_{\alpha}\tau}{M}\Big)+
  \Theta_{j,f,\alpha,M}^{-}\Big(\frac{M_{\alpha}\tau}{M}\Big).
\end{equation}

\begin{remark}   \label{remparity}
  If~$j$ and~$f$ have opposite parities, i.e. in the
    case~\eqref{eqparity},
    then the function~$\hh{\phi}^-_{j,f,\alpha,M}$ happens to be zero
  and the third term in~\eqref{eq:Taunearalpha},
  $\Theta_{j,f,\alpha,M}^{-}$, is thus absent.
  This is what will happen with the function $f=\chi_j$ of Proposition~\ref{propPsiassumpts}.
  If moreover~$f$ is an odd function (thus assuming~$j$ even),
  the first term is trivially absent because $f_{\alpha/M}$ is
  odd too, whence $\langle f_{\alpha/M}\rangle=0$.
  Less trivially, as shown below, in the case $f=\chi_j$ the first term will always be
  absent, even when~$j$ is odd and~$\chi_j$ is even.
  \end{remark}

\subsection{A resurgent asymptotic expansion}  \label{sec:resurgence}
We now present the asymptotic expansion of
$\Theta(\al+\tau; j,f,M)$ for $\tau$ near $0$. First we observe
that the third term in~\eqref{eq:Taunearalpha} is always exponentially small:
$\Theta_{j,f,\alpha,M}^{-}(\tau)=\cO(e^{-c\IM(-1/{\tau})})$
for sufficiently small $c>0$ (this
follows easily from the fact that $\Theta_{j,f,\alpha,M}^{-}$ is the
difference of two Laplace transforms of the same function). As for the
second term, consider
the $L$-function
$ L(s,f_{\alpha/M})=\sum_{m=1}^{\infty} f_{\alpha/M}(m)m^{-s}$; it
has a meromorphic continuation to~$\C$, and for all positive integers~$n$, we have that
\begin{equation} \label{eq:Lseriesevaluation}
    L(-n,f_{\alpha/M})= -\frac{M_{\alpha}^n}{n+1}
    \sum_{m=1}^{M_{\alpha}}
    B_{n+1}\!\Big(\frac{m}{M_{\alpha}}\Big) f_{\al/M}(m), 
\end{equation}
where $B_{n+1}(x)=\sum_{k=0}^{n+1} \binom{n+1}{k}B_{n+1-k} \, x^k$ is the
$(n+1)^{\text{th}}$ Bernoulli polynomial, and $(B_\ell)_{\ell\ge0}$ is the sequence of Bernoulli numbers.
Define the formal series
\begin{equation} \label{eq:ResurgentSeries}
  \wt\Theta_{j,f,\alpha,M}(\tau):=
  \sum_{p=0}^{\infty}\frac{1}{p!}L(-2p-j,f_{\alpha/M}) \Big(\frac{\pi i}{M}\Big)^p\tau^p.
  %
\end{equation}
By \cite[Theorem~2 \& Remark~3.2]{han2022resurgence} the series
\eqref{eq:ResurgentSeries} is resurgent,
with a Borel transform all of whose singular points are of the form
$i\pi m^2/M$, $m\in\Z_{\ge1}$;
it is Borel summable in all direction except~$\pi/2$
%
%
and its median Borel sum in the direction $\pi/2$ is
%
%
\begin{equation}   \label{eqbSmed}
\frac12\big( \bS^{\frac\pi2-\epsilon} + \bS^{\frac\pi2+\epsilon} \big)
\wt\Theta_{j,f,\alpha,M}(\tau) =
\Theta_{j,f,\alpha,M}^{+}\Big(\frac{M_\al\tau}{M}\Big)
\quad \text{for} \ens \tau\in \HH.
\end{equation}
%
%
\bblack
%
%
It follows that the first term in~\eqref{eq:Taunearalpha} is the
dominant one if $\langle f_{\alpha/M} \rangle\neq0$, and in fact
\begin{equation}   \label{eqexistlimal}
  \lim_{\tau\to\al}\Theta(\tau; j,f,M) \;\text{exists}
  \ens \Longleftrightarrow \ens 
  \langle f_{\alpha/M} \rangle=0
\end{equation}
(with reference to a non-tangential limit, i.e.\ with
$\arg(\tau-\al)\in I$ for an arbitrary compact interval $I\subset (0,\pi)$).
Let us define
\begin{equation}   \label{eqdefsQfM}
  \sQ_{f,M} := \{\, \al \in \Q \mid \langle f_{\alpha/M} \rangle=0
  \,\}.
\end{equation}
We remark that $\lim\limits_{\tau\to\al}\Theta(\tau; j,f,M)$ exists for
all $\al \in \sQ_{f,M} \cup\{\infty\}$.

\begin{remark}    \label{remtrivsQ}
  Trivially,
    if~$f$ is an odd function, then $\sQ_{f,M} = \Q$.
  \end{remark}

When $\al \in \sQ_{f,M}$, we thus have the
following Poincar\'e asymptotic expansion
\begin{equation}  \label{eq:generalexpansion}
  \Theta(\al+\tau; j,f,M)\underset{\tau\to0}{\sim}
  \wt\Theta_{j,f,\alpha,M}(\tau) 
  %
  %
\end{equation}
and, in particular,
\begin{equation}\label{equationlimitvalueTheta}
  \lim_{\tau\to\al}\Theta(\tau; j,f,M) =
L(-j,f_{\al/M}) = - \frac{M_{\alpha}^j}{j+1}\sum_{m=1}^{M_{\alpha}}
B_{j+1}\!\Big(\frac{m}{M_{\alpha}}\Big)f(m)\exp(\pi im^2 \alpha/M). 
\end{equation}

\subsection{Quantum modularity of partial theta series} \label{sec:QM}

We quote here \cite[Theorem~7]{han2022resurgence}: 

\begin{theorem} \label{thmqmpts}
  Suppose that $f(0)=0$ and there exists $n_0\in\Z$ such that, for all $n\in\Z$,
  \begin{equation}   \label{eqsupportcondition}
    f(n)\neq0 \iimp n^2 = n_0^2 \!\mod 2M.
  \end{equation}
  Then $\sQ_{f,M}$ is a dense subset of~$\Q$ such that
  $\sQ_{f,M}\cup\{\infty\}$ is invariant under the action of~$\Ga_1(2M)$.

  Suppose moreover that $j=0$ or~$1$ and~\eqref{eqparity} holds. Then
  $\Theta(\cdot\,;0,f,M)$ (resp.\ $\Theta(\cdot\,;1,f,M)$) is a strong
  quantum modular form on $\Gamma_1(2M)$ with quantum set~$\sQ_{f,M}$
  and weight~$\frac12$ (resp.~$\frac32$).
\end{theorem}


Notice we are following the convention of
Remark~\ref{remsqmfupper}: we call the partial theta series
$\tau\mapsto\Theta(\tau;j,f,M)$, $j=0,1$, a strong quantum modular
form instead of referring to the family of formal series defined by
its asymptotic expansions, $\wt\ph_\al(\tau) :=
\wt\Theta_{j,f,\alpha,M}(\tau)$ for $\al\in\sQ_{f,M}$.

We can be more specific. In the situation described by Theorem~\ref{thmqmpts},
let $\gamma= \begin{psmallmatrix} a&b \\ c&d \end{psmallmatrix}\in\Gamma_1(2M)$
and take any~$n_0$ in the support of~$f$.

\noindent
-- If $c=0$, then~$\ga$ acts on~$\HH$ like an integer translation and
\[ \Theta(\tau;j,f,M) = e^{-i\pi n_0^2/M}\Theta(\tau+1;j,f,M). \]
In fact,~\eqref{eqsupportcondition} ensures that $e^{-i\pi
  n_0^2\tau/M}\Theta(\tau;j,f,M)$ is a holomorphic function of $q=e^{2\pi i \tau}$.

\noindent
-- If $c\neq0$, then we can assume $c>0$ without loss of generality and,
combining formulas~(7.5),~(7.7) and~(7.9)
from\cite{han2022resurgence}, we get a set of two identities for each
parity case:
\begin{align}
\label{equationmodularobstructionfodd}
  \text{$f$ odd} &\imp
  \Theta(\tau;0,f,M) - \eps(\ga)
  J_\ga(\tau)^{-\frac12} \Theta(\ga \tau;0,f,M) =
  \bS^{\frac\pi2 \mp \epsilon}\wt\Theta_{0,f,-\frac{d}{c},M}(\tau+\tfrac dc)
  \\[1ex]
  \label{equationmodularobstructionfeven}
  \text{$f$ even} &\imp
  \Theta(\tau;1,f,M) - \eps(\ga)
  J_\ga(\tau)^{-\frac32} \Theta(\ga \tau;1,f,M) =
  \bS^{\frac\pi2 \mp\epsilon}\wt\Theta_{1,f,-\frac{d}{c},M}(\tau+\tfrac dc)
%
%
\end{align}
where~$\eps$ involves the Jacobi symbol: $\eps(\ga) := \big(\frac{2Mc}{|d|}\big) e^{-{i\pi n_0^2 b}/{M}}$,
and the branch of the square root of $J_\ga(\tau)=c\tau+d$ to be used in the
left-hand side of~\eqref{equationmodularobstructionfodd} (through its
inverse) or~\eqref{equationmodularobstructionfeven} (through the cube
of its inverse) depends on the choice of sign~`$\mp$' in the
right-hand side, namely
\begin{itemize}
  \item
        choosing~`$-$': the lateral summation
        $\bS^{\frac\pi2-\epsilon}$ gives right-hand sides that
        extend holomorphically to the cut plane $\{\arg(\tau+\frac
        dc)\in(-\pi,\pi)\} = \C\setminus(-\infty,-\frac dc]$
        and
        \eqref{equationmodularobstructionfodd}--\eqref{equationmodularobstructionfeven}
        hold true there provided the left-hand sides involve the
        principal branch of $J_\ga(\tau)^{\frac12}$ (with positive real part);
%
    %
\item
  choosing~`$+$': the lateral summation $\bS^{\frac\pi2+\epsilon}$
  gives right-hand sides that extend holomorphically to the cut
  plane
  $\{\arg(\tau+\frac dc)\in(0,2\pi)\} = \C\setminus[-\frac
  dc,+\infty)$
  and
  \eqref{equationmodularobstructionfodd}--\eqref{equationmodularobstructionfeven}
  hold true there provided the left-hand sides involve the opposite of
  the analytic continuation of the principal branch of
  $J_\ga(\tau)^{\frac12}$ (i.e.\ we use the branch of
  $J_\ga(\tau)^{\frac12}$ that has negative imaginary part).
\end{itemize}


In the notation of Section~\ref{secqmf}, we thus get a quantum modular
form
$\al\in\sQ_{f,M} \mapsto \ph(\al) := \lim\limits_{\tau\to\al}
\Theta(\tau;j,f,M)$ whose modularity defect extends analytically to
$\R\setminus\{-\frac dc\}$.

%
We emphasize that the extension property directly stems from the domains of analyticity
of the lateral sums of $\wt\Theta_{j,f,-\frac{d}{c},M}$:
since the only singularities of the Borel transform are on
$e^{i\frac\pi2}\R_{>0}$,
we can freely vary $\tht=\frac\pi2-\epsilon$ in $(-\frac\pi2,\frac\pi2)$ and
the corresponding Borel sums $\bS^\tht\wt\Theta_{j,f,-\frac{d}{c},M}(\tau)$
mutually extend (including the standard Borel sum $\bS^0$ associated
with the usual Laplace transform~$\cL^0$), resulting in the large
domain of analyticity $\arg\tau\in(-\pi,\pi)$ indicated above (actually, we can even
decrease~$\tht$ below~$-\frac\pi2$ provided we stop before
$-\frac{3\pi}{2}$, and get a domain as large as
$\arg\tau\in(-2\pi,\pi)$).
Similarly,
$\bS^{\frac\pi2+\epsilon}\wt\Theta_{j,f,-\frac{d}{c},M}(\tau)$ extends
as far as $\arg\tau\in(0,3\pi)$.
It is only when we consider both lateral sums simultaneously, as
in~\eqref{eqbSmed}, that we must restrict to $\arg\tau\in(0,\pi)$, i.e.\ to~$\HH$.

%
%
\bblack

There are also formulas for the action on $\Theta(\cdot\,;j,f,M)$ of
an arbitrary element of the full modular group~$\SL(2,\Z)$---see
\cite[Sec.~7]{han2022resurgence} and Appendix~\ref{secVVsqmf} below---which imply that
$\Theta(\cdot\,;j,f,M)$ is the first component of a strong quantum
modular form on~$\SL(2,\Z)$.
The idea of the proof of all these formulas is to analyze the action
of the generators $T=\begin{psmallmatrix} 1&1 \\
  0&1 \end{psmallmatrix}$
and $S=\begin{psmallmatrix} 0&-1 \\ 1&0 \end{psmallmatrix}$.
The key point is that we can write the median Borel sum
$\Theta_{j,f,0,M}^+(\tau)$ of~\eqref{eqThetapmLaplTrsf} as the sum of
a lateral Borel sum plus half the difference of two lateral Borel
sums, and compute the latter difference as a sum of the contributions
of the singularities of the Borel transform; we are then naturally led to
consider $\Theta(S\tau;j,\wh f,M)$, where $\wh f=U_M f$ is the Discrete
Fourier Transform\footnote{This is the parity-preserving
  operator defined by \label{ftnUM}
\[
U_M \col f \mapsto \wh f, \quad
\wh f(n) := \frac{1}{\sqrt M}\, \sum_{\ell\!\!\mod M}f(\ell) e^{-2 i\pi \ell n/M}
\ens\text{for all}\; n\in\Z.
\]
  } (DFT) of~$f$.

The computation is given in \cite{han2022resurgence} in terms of
\'Ecalle's alien derivations~$\De_\om$, which are fundamental tools in
Resurgence Theory.
When acting on a resurgent formal series, $\De_\om$ measures the
singular behaviour of its Borel transform at~$\om$ (it thus
annihilates any series whose Borel transform has all its branches
regular ar~$\om$) and satisfies the product rule---see
\cite{Ecalle81}, \cite{Eca93}, \cite{MS16}, \cite[Sec.~6]{han2022resurgence}.
Here we must use $\om=\xi_m=i\pi m^2/M$, $m \in
\Z_{\ge1}$.
%
%
Here is a sample of ``alien calculus'' that will help us to grasp
quantum modularity for $j\ge2$:

\begin{proposition}   \label{propUnivPol}
  The formulas
  \beglab{equationlemmarecursionpolynomial}
        P_0(x) =1, \quad P_1(x) =-x, \quad
        P_j(x) = \big( 2x^2 - (j-1)\big) P_{j-2}(x) - xP_{j-2}'(x)
        \quad \text{for} \ens j\ge2,
        \edla
  inductively define a sequence of integer polynomials $(P_j)_{j\ge0}$
  of the form
%
  %
  \beglab{eqformofPj}
        P_{j}(x) = \sum_{0\leq \nu \leq j \atop \nu \equiv j\
          [2]} P_{j,\nu} \, x^{\nu} \in \mathbb{Z}[x]
        \quad \text{with}\ens P_{j,j}=(-1)^j \, 2^{[\frac{j}{2}]},
        \edla
      where $[\frac{j}{2}]$ is the greatest integer $\le \frac{j}{2}$,
      %
      and for any $M \in \mathbb{Z}_{\geq 1}$,
      $j\in \mathbb{Z}_{\geq 0}$ and $f$ $M$-periodic function
      satisfying~\eqref{eqparity},
    \begin{equation}\label{equationalienhigherdepth}
      \Delta_{\frac{\pi i n^2 }{M}}\wt\Theta_{j,f,0,M}(\tau)
      = A_{j,M} \, \wh{f}(n)\,\tau^{-\frac{j+1}{2}} P_j\Big(n\Big(\frac{\pi i}{M}\Big)^{\!1/2}\tau^{-1/2}\Big)
      \quad \text{with} \ens
      A_{j,M} := 2^{-[\frac{j}{2}]+1} i^{\frac12} \, \Big( \frac{M}{\pi i}\Big)^{\!j/2},
    \end{equation}
    where $\wh f := U_M f$ is the DFT of~$f$.
As a consequence, if moreover $\langle f\rangle =0$, we obtain a set
of two identities (one for each choice of sign):
%
%
    \begin{equation}\label{equationShigherdepth}
      \Theta(\tau;j,f,M) = \bS^{\frac\pi2\mp \epsilon}\wt\Theta_{j,f,0,M}(\tau)
      \mp 2^{-[\frac{j}{2}]} \, i^{\frac12}
          \sum\limits_{0 \leq \nu \leq j \atop \nu \equiv j \ [2]}\Big( \frac{M}{\pi i}\Big)^{\!\frac {j-\nu}2} 
          P_{j,\nu} \, \tau^{-\frac{j+\nu+1}{2}}
          \Theta(-\tau^{-1};\nu,\wh f,M)
        \end{equation}
for $\tau \in \HH$ (here we use the principal branch of the square root:
$\arg(\tau^{1/2})\in(0,\pi/4)$ in~\eqref{equationShigherdepth}
and $i^{\frac12}=e^{i\pi/4}$ in~\eqref{equationalienhigherdepth}--\eqref{equationShigherdepth}).
\end{proposition}

Note that~\eqref{equationalienhigherdepth} says that
$\Delta_{\xi_n}\wt\Theta_{j,f,0,M}(\tau)$ (with $\xi_n=\pi i n^2/M$)
is a sum of terms proportional to $\tau^{-k-\frac12}$, with~$k$ integer between
$[\frac{j+1}2]$ and~$j$; the latter monomial represents a singularity of the form
$\frac1{2\Ga(k+\frac12)} (\xi-\xi_n)^{-k-\frac32} + $ 
regular germ at~$\xi_n$ for the Borel transform of $\wt\Theta_{j,f,0,M}$.

\bblack

\begin{proof}
  The cases $j=0$ and~$1$ of~\eqref{equationalienhigherdepth} are
  in \cite[eqns~(6.1)--(6.2)]{han2022resurgence}.
For $j\ge2$, making use of the relations 
\begin{equation}   \label{eqderivThetaj}
  \Theta(\tau;j,f,M) = \frac{M}{\pi i}\frac{d}{d\tau} \Theta(\tau;j-2,f,M) 
  \quad \text{and} \quad
  \Delta_{\frac{\pi in^2}{M}}\frac{d}{d\tau} =
  \Big( \frac{\pi in^2}{M} \tau^{-2} + \frac{d}{d\tau} \Big)
  \Delta_{\frac{\pi in^2}{M}},
  \end{equation}
and setting $x:=n\left(\frac{\pi i}{M}\right)^{1/2}\tau^{-1/2}$,
we compute by induction that 
\begin{align*}
  \Delta_{\frac{\pi in^2}{M}} \widetilde{\Theta}_{j,f,0,M} 
  &= \frac{M}{\pi i}\Delta_{\frac{\pi in^2}{M}} \frac d{d\tau} \widetilde{\Theta}_{j-2,f,0,M} 
  \\[1ex]
  &= \frac{M}{\pi i} A_{j-2,M} \wh f(n)
    \Big( x^2 \tau\ii + \frac{d}{d\tau} \Big)
    \big[ \tau^{-\frac{j-1}2} P_{j-2}(x)\big]
  \\[1ex]
  &= 2 A_{j,M} \wh f(n) \big( x^2 \tau^{-\frac{j+1}2} P_{j-2}(x)
    - \tfrac{j-1}2 \tau^{-\frac{j+1}2} P_{j-2}(x)
    -\tfrac12 \tau^{-\frac{j-1}2} \tau\ii x P_{j-2}'(x) \big)    
  \\[1ex]
  &= A_{j,M} \wh f(n) \tau^{-\frac{j+1}2} 
     \big(2x^2 P_{j-2} (x) - (j-1) P_{j-2}(x) - xP_{j-2}^\prime(x)\big),
\end{align*}
which yields~\eqref{equationalienhigherdepth}.

It is obvious that each coefficient $P_{j,\nu}$ is integer, that~$P_j$
has a parity matching that of the index~$j$, and that the degree of
the polynomial~$P_j$ is precisely~$j$ with leading coefficient as in~\eqref{eqformofPj}.


Formula~\eqref{equationShigherdepth} is obtained
from~\eqref{eq:Taunearalpha} and~\eqref{eqbSmed} by writing
$\Theta(\tau;j,f,M) - \bS^{\frac\pi2\mp
  \epsilon}\wt\Theta_{j,f,0,M}(\tau)$
as half the difference of two Laplace transforms, the integrand having
singularities at the points $\xi_n=i\pi n^2/M$. 
For instance, 
$\Theta(\tau;j,f,M) - \bS^{\frac\pi2+\epsilon}\wt\Theta_{j,f,0,M}(\tau)$
is half of $(\bS^{\frac\pi2-\epsilon}- \bS^{\frac\pi2+\epsilon})\wt\Theta_{j,f,0,M}(\tau)
= \sum\limits_{n\ge1} e^{-\xi_n/\tau}
\bS^{\frac\pi2+\epsilon}\De_{\xi_n}\wt\Theta_{j,f,0,M}(\tau) =
A_{j,M} \tau^{-\frac{j+1}{2}} \sum\limits_{n\ge1}
e^{-\xi_n/\tau} \wh{f}(n)\, P_j\big(n\big(\frac{\pi i}{M}\big)^{1/2}\tau^{-1/2}\big)$.
\end{proof}

We will make use of the universal polynomials~$P_j$ in our proof
  of Theorem~\ref{Thm:main} in Section~\ref{secWAECSFHS}.
  We now return to the context of Theorem~\ref{thmqmpts} and deduce
  from it the higher
depth version for general values of~$j$:

\begin{corollary}\label{corollarypartialthetaserieshigherdepth}
  Suppose that $f(0)=0$ and there exists $n_0\in\Z$ such
  that~\eqref{eqsupportcondition} holds.
  Suppose moreover that $j\in\Z_{\ge0}$ and~\eqref{eqparity} holds.
  Then $\Theta(\cdot\,;j,f,M)$ is a depth $[j/2]+1$ strong quantum modular form
  on~$\Ga_1(2M)$
  with quantum set~$\sQ_{f,M}$ and weight $j+\frac{1}{2}$.
\end{corollary}
\begin{proof}
  The function $\Theta(\tau;j,f,M)$ can be obtained from
  $\Theta(\tau;0,f,M)$ or $\Theta(\tau;1,f,M)$ by applying the first
  part of~\eqref{eqderivThetaj} $[j/2]$ times. The modularity defect
  is then represented (up to a constant factor in~$\C^*$) by the
  $[j/2]$-th derivative of
  formulas~\eqref{equationmodularobstructionfodd}
  or~\eqref{equationmodularobstructionfeven}. The desired result
  follows from the fact that
  $\frac{d}{d\tau} \big[\Theta(\cdot\,;0,f,M)\circ\ga\big] =
    J_\ga^{-2} \, \big[ \frac{d}{d\tau} \Theta(\cdot\,;0,f,M) \big]\circ\ga$.
 \end{proof}

See Remark~\ref{remcongrsubgp} for a slightly different viewpoint on
 Corollary~\ref{corollarypartialthetaserieshigherdepth}.


 \section{Quantum modularity of the GPPV invariant of a Seifert
   fibered integral homology sphere} \label{sec:proofs}

Let~$X$ be as in~\eqref{defX}, with $p_1,\ldots,p_r$ positive pairwise coprime integers
($r\ge3$), among which only~$p_1$ may be even.
This section aims at proving Theorem~\ref{thm:QM}.

We begin by recalling formulas from~\cite{AM22} for the WRT invariants
$\WRTk(X)$ and the modified 
GPPV invariant~$\Psi(\tau)$.
%
%
Recall that a rational function~$G$ was defined in
Equation~\eqref{eq:Taylor} and $P = p_1\cdots p_r$.


\subsection{The WRT invariant of~$X$} \label{sec:defWRT}

For the presentation of $\WRTk(X)$ with $k\in\Z_{\ge2}$ we follow
\cite{LR99}, and use almost identical notation. Some formulas simplify
as $P>0$ and as we assume that $H:=\lvert H_1(X,\Z)\rvert=1$.
Let~$y$ be a complex variable, and define
\begin{equation}   \label{eqdefFy}
  F(y):=(e^{y/2}-e^{-y/2})^{2-r}\prod_{j=1}^r(e^{y/(2p_j)}-e^{-y/(2p_j)}) = G(e^{\frac y{2P}}),
  \qquad
  g(y):=iy^2/(8\pi P).
\end{equation}
Let $C':=\mathbb{R} e^{\pi i/4} \subset \C$ be the oriented contour
with orientation induced by $t\in\R \mapsto te^{\pi i/4}\in C'$.
Recall the surgery procedure described in Section
\ref{sec:RationalSurgery}. Consider the framed oriented link $L$ given
by an unknot $U$ which is linked once with $r$ split disjoint unknots
$U_1,\ldots,U_r$. Consider the surgery data $B$ given by assigning the
tuple of rationals $(p_j/q_j)_{j\in \{1,\ldots,r\}}$ to tuple of unknots
$(U_j)_{j \in \{1,\ldots,r\}}$, and assigning $0/1$ to $U$. This is a
surgery link for $X$, i.e.\ we have an orientation preserving
diffeomorphism $S^3_{L,B} \cong X$.  By applying the surgery formula
\eqref{eq:rationalsurgeryformula} to $(L,B)$ Lawrence and Rozansky
prove the following identity \cite[eqn~(4.8)]{LR99}:\footnote{%
      %
      See Appendix~\ref{Appendix:Normalizations} for explanations on the normalization issues for WRT invariants.
      %
    %
}
\begin{equation} \label{eq:LRformula}
  %
  %
  \frac{4\, e^{\frac{i\pi}{4}+\frac{i\pi\phi}{2k}} \sqrt P}{G_{k,0}}
  \WRTk(X)=
  %
  %
  \frac{1}{2\pi i }\int_{C'} F(y)e^{kg(y)} dy -
  \sum_{m=1}^{2P-1} \Res\left(\frac{F(y)e^{kg(y)}}{1-e^{-ky}}, y=2\pi im\right) 
\end{equation}
with
%
%
$G_{k,0}$ as in~\eqref{eqGkzero}
%
and $\phi\in\Q$ defined in terms of Dedekind sums or,
  equivalently, retrieved from the Casson invariant~$\la$ by the
  formula \cite[eqn~(4.1)]{LR99}:
    \beglab{eqdefphi}
\phi = -24 \la -P \Big(r-2 - \sum_{1\le j\le r} \frac1{p_j^2} \Big).
    \edla
According to the proof of
    %
    %
    \cite[Theorem~2]{AM22}, we can rewrite the first term of the \rhs\
of~\eqref{eq:LRformula} as
\begin{equation} \label{eq:Borel}
  \frac{1}{2\pi i }\int_{C'} F(y)e^{kg(y)} dy= \int_0^{+\infty} e^{-k\xi}\cB_0(\xi) d \xi
  \quad \text{with} \ens
  \cB_0(\xi) := 
2 (2\pi i\xi/P)^{-1/2}
  %
%
    G\big( e^{(2\pi i\xi/P)^{1/2}} \big).
    %
  \end{equation}
  %
  %
In fact, according to \cite[eqns~(1.7)\&(2.6)]{AM22}, this function~$\cB_0$ is the
Borel transform of a suitable normalization of the Ohtsuki series
(this fact will also follow from Section~\ref{secpfthmmain}).


\subsection{The modified GPPV invariant of~$X$}   \label{secGPPVX}


Theorem~3 from \cite{AM22} essentially says that the modified GPPV
invariant of~\eqref{partialthetaseries} can be written as
%
%
\begin{equation}   \label{eqPsiasmedBzero}
\Psi(\tau)
  = \frac{(-1)^r e^{\pi i/4}}{\sqrt{2P}} \tau^{-1/2} \times \frac12
  \left(\cL^{\frac\pi2-\epsilon}+\cL^{\frac\pi2+\epsilon}\right)[\cB_0](\tau),
\end{equation}
i.e.\ that $\tau^{1/2} \Psi(\tau)$ is a multiple of the median Borel
sum of the Ohtsuki series since, when the Borel transform is
meromorphic, median Borel summation amounts to taking the arithmetic
average of lateral Borel sums---this is a slight extension of
  Footnote~\ref{ftnLaplMedSum} inasmusch as we are now dealing with
  half-integer powers: we shall see that
  $\cB_0(\xi) = \sum_{p\ge1} \cB_{0,p} \xi^{p-\frac12}/\Ga(p+\frac12)
  \in \xi^{1/2}\C\{\xi\}$ for some sequence of complex coefficients $(\cB_{0,p})_{p\ge1}$, thus
\beglab{eqasymptcBzero}
\tfrac12\left(\cL^{\frac\pi2-\epsilon}+\cL^{\frac\pi2+\epsilon}\right)[\cB_0](\tau)
\sim \sum_{p\ge1} \cB_{0,p} \tau^{p+\frac12}
 \quad \text{as}
 \; \tau\to0 \; \text{non-tangentially from within~$\HH$.} 
\edla

Formula~\eqref{eqPsiasmedBzero} can be recovered from \cite[Theorem~1]{han2022resurgence} and the beginning of
Section~3 from \cite{han2022resurgence} as follows:
since $\Psi(\tau) = \sum_{m\ge m_0} \tilde\chi(m) e^{\frac{i\pi
    m^2\tau}{2P}}$,
we consider $F_{\ti\chi}(t) := \sum_{m\ge m_0}\ti\chi(m) e^{-mt}$,
which is convergent for $\RE t>0$ only and coincides with 
\begla
(-1)^r G(e^{-t}) =
(e^{-Pt}-e^{Pt})^{-(r-2)} \prod_{j=1}^r (e^{Pt/p_j}-e^{-Pt/p_j}).
%
%
%
\edla
The function~$F_{\ti\chi}$ thus has an analytic continuation
to~$\C\setminus\frac{i\pi}P\Z$, that clearly is even and meromorphic in~$\C$.
It is easily seen that $t=0$ is not a pole but a zero of that function.
%
%
We follow \cite[eqns~(3.1)--(3.2)]{han2022resurgence} and
define implicitly $\wh\phi^+_{\ti\chi}(\xi)$ by
\begin{equation}\label{equationGPPVgerm}
  (-1)^r G(e^{-t})  
  = \pi^{1/2}\hh{\phi}^+_{\ti\chi}\Big(\frac{t^2}{C^2}\Big)
\quad \text{with} \ens
C := (2\pi i/P)^{1/2}.
\end{equation}
%
%
In our case, $\wh\phi^+_{\ti\chi}(\xi) \in\xi\C\{\xi\}$
and \cite[eqns~(3.4)--(3.5)]{han2022resurgence} yields
\begin{equation}   \label{eqPsimLaplphip}
  \Psi(\tau)
  = \tau^{-1/2} \times \tfrac12 (\cL^{\pi/2-\epsilon}+\cL^{\pi/2+\epsilon})
\big[\xi^{-1/2}\hh{\phi}_{\chi}^+(\xi)\big], 
\end{equation}
%
%
which is equivalent to~\eqref{eqPsiasmedBzero}
because~\eqref{equationGPPVgerm} gives
$\hh{\phi}^+_{\chi}(\xi)= (-1)^r{\pi}^{-1/2}G\big( e^{C \xi^{1/2}} \big)$,
which coincides with
$\frac{(-1)^re^{\pi i/4}}{\sqrt{2P}} \xi^{1/2} \mathcal{B}_0(\xi)$.
The Puiseux expansion that we have indicated for~$\cB_0(\xi)$ just
before~\eqref{eqasymptcBzero} stems from the Taylor expansion of
$\wh\phi^+_{\ti\chi}(\xi) \in\xi\C\{\xi\}$
(in particular the absence of a coefficient $\cB_{0,0}$ is equivalent
to the vanishing of~$\wh\phi^+_{\ti\chi}(\xi)$ at $\xi=0$).

Recall that the set~$\Omb$ of poles of
$\mathcal{B}_0$ is contained in $i \R_{>0}$. By applying Cauchy's
residue theorem, it was proven in \cite[Lemma~14]{AM22} that
%
%
\begin{equation}
  \Psi(\tau)
  =\frac{ (-1)^r e^{\pi i/4}}{\sqrt{2P}} \, \tau^{-1/2} \left( \int_0^{+\infty}
    e^{-\zeta/\tau}\mathcal{B}_0(\zeta) d \zeta- \pi i\sum_{\xi \in
      \Ombb}\Res(e^{-\zeta/\tau}\mathcal{B}_0(\zeta), \zeta=\xi)
  \right).
\end{equation}
The function
$\mathcal{R}(\tau):=\pi i\sum_{\xi \in
  \Ombb}\Res(e^{-\zeta/\tau}\mathcal{B}_0(\zeta), \zeta=\xi)$ can be
rewritten as a polynomial in $\tau^{-1}$ with coefficients in partial
theta series evaluated at $S\tau=-\tau^{-1}$ (the coefficients of these
partial theta series being given by periodic sequences). By applying
\cite[Proposition~15]{AM22},
which is a version of the asymptotic
expansion presented in Section~\ref{sec:resurgence},
it was shown in the proof of \cite[Theorem~4]{AM22} that, provided
that the relevant mean values vanish, we have that
\begin{equation}
  \lim_{\tau \rightarrow 1/k} R(\tau) =
  \sum_{m=1}^{2P-1} \Res\left(\frac{F(y)e^{kg(y)}}{1-e^{-ky}}, y=2\pi im\right),
\end{equation}
which is nothing but the second term of the \rhs\ of~\eqref{eq:LRformula}.
The relevant vanishing follows from Proposition~\ref{prop:vanishingmeanvalue}
below.
As
\begin{equation} 
  \frac{e^{\pi i/4} k^{1/2}}{\sqrt{2P}} =2
  \big( \exp(\pi i/k)-\exp(-\pi i/k)\big) 
  B G_{k,0},
\end{equation}
this implies by equation~\eqref{eq:LRformula} that the following limit
holds, which is a generalization of \cite[Theorem~4]{AM22}:
\begin{corollary} \label{cor:1/klimit}
  For every $k \in \Z_{\ge2}$, it follows that
  %
  %
  \begin{equation}  \label{eq1/klimit}
\lim_{\tau \rightarrow 1/k} \Psi(\tau) =2 (-1)^r
\big( \exp(\pi i/k)-\exp(-\pi i/k)\big) 
e^{\frac{i\pi\phi}{2k}} \WRTk(X).
 \end{equation}
\end{corollary}



\begin{remark}   \label{remcasekis1}
Formula~\eqref{eq1/klimit} also holds true when $k=1$, with~$0$ in the \rhs\ (with the usual
convention $\WRT_1(X)=1$). We even have
\beglab{eqlimellPsizero}
\lim_{\tau \rightarrow \ell} \Psi(\tau) =0
\quad \text{for every}\ens \ell\in\Z.
\edla
Indeed, in view of the absence of a coefficient $\cB_{0,0}$
in~\eqref{eqasymptcBzero},
the case $\ell=0$ of~\eqref{eqlimellPsizero} follows from~\eqref{eqPsiasmedBzero},
and since~\eqref{partialthetaseries} says that~$\Psi(\tau)$ is the
product of $e^{\frac{i\pi m_0^2\tau}{2P}}$ and a $1$-periodic function
of~$\tau$, we get the same result for any $\ell\in\Z$.
\end{remark}


\subsection{The normalized GPPV invariant of~$X$ and Hikami functions}
\label{secGPPVHik}





As announced in the introduction, we will now recast the modified GPPV
invariant of~\eqref{partialthetaseries} as a sum of partial theta
series of the form studied in \cite{han2022resurgence}, which will
allow us to make use of the results of Section~\ref{sec:QMresurgence}.

As a preliminary remark, we note that, for our Seifert fibered integral homology
sphere~$X$,
%
%
the integer~$m_0$ of~\eqref{eqformulamzero} is positive
%
%
%
%
in all cases but one: in the case of the Poincar\'e homology sphere, i.e.\ when $r=3$ and
$(p_1,p_2,p_3)=(2,3,5)$,
we have $m_0=-1$, whereas in all other cases $m_0\ge1$
%
%
(we leave it to the reader to check this elementary fact).
%

\begin{proposition}   \label{propPsiassumpts}
  There is an odd function $\chi \col \Z \to \Z$ of the form
  \begin{equation}   \label{eqtichisummjchij}
\chi(m) = \sum_{j=0}^{r-3} m^j \chi_j(m)
    \quad\text{for~$m\in\Z$}
  \end{equation}
with $2P$-periodic functions of alternating parities
  \begin{equation}
    \chi_0, \ldots, \chi_{r-3} \col \Z \to \Q,
    \quad
    \text{$j$ even} \iimp \text{$\chi_j$ odd function,}
    \quad
    \text{$j$ odd} \iimp \text{$\chi_j$ even function,}
  \end{equation}
  such that the integer coefficients~$\ti\chi(m)$
  of~\eqref{eq:Taylor}
  satisfy
  %
  %
  \begin{align}   \label{eqtichisummjchijti}
    \sum_{m\ge m_0} \ti\chi(m) z^m = & \sum_{m\ge1} \chi(m) z^m
    \hspace{-4em} &\text{if} \ens (p_1,\ldots,p_r)&\neq(2,3,5) \\[1ex]
    \label{eqtichisummjchijEXCEPT}
    \sum_{m\ge m_0} \ti\chi(m) z^m = -z\ii-z &+\sum_{m\ge1} \chi(m) z^m
    \hspace{-4em}  &\text{if} \ens (p_1,\ldots,p_r)&=(2,3,5).
    %
    %
  \end{align}
  %
  %
  %
  Consequently, with reference to notation~\eqref{eqdefThetafjM},
  the modified GPPV invariant~\eqref{partialthetaseries} of~$X$ can be written
  \begin{equation}   \label{eqPsicQsumThchij}
    \Psi(\tau) = \cQ(\tau) + \sum_{j=0}^{r-3} \Theta(\tau;j,\chi_j,2P)
    %
    %
    %
    %
  \end{equation}
  with $\cQ=0$ except when
    $(p_1,\ldots,p_r)=(2,3,5)$, in which case $\cQ(\tau) = -e^{-i\pi\tau/60}-e^{i\pi\tau/60}$.
  \end{proposition}

  
Subsection~\ref{secGPPVHik} is devoted to proving this proposition and deriving formulas for
the~$\chi_j$'s, which will be given in Remark~\ref{remdefcT}.
To that end, we introduce the notation
\begla
E := \{+1,-1\}^r,
\qquad
-\ue := (-\eps_1,\ldots,-\eps_r)
\quad \text{for} \ens
\ue := (\eps_1,\ldots,\eps_r) \in E,
\edla
and define the odd function
\beglab{eqdefwhpj}
\sN \col \ue \in E \mapsto \sN(\ue) :=
\sum_{j=1}^r \eps_j \whp_j \in \Z,
\qquad \text{with the notation} \ens
\whp_j := \frac P{p_j}.
\edla
One can check that~$\sN$ is injective. In fact, even the
composition
\beglab{eqtpcNoddinj}
\tp{\sN} := \tp{\ \cdot\ } \circ \sN \col E  \to \Z/2P\Z
\quad \text{is odd and injective,}
\edla
where we use the notation 
\begin{equation}   \label{eqdefnotabrN}
  \brN{\ \cdot\ } \col n\in \Z \mapsto \brN n :=n+N\Z \in \Z/N\Z
  \quad\text{for any $N\in\Z_{\ge1}$}
\end{equation}
for the canonical projection.
Indeed: suppose $\tp{\sN}(\ue')=\tp{\sN}(\ue)$ with
$\ue,\ue'\in E$;
notice that for each $i\in\{1,\ldots,r\}$,
$\frac12(\eps'_i-\eps_i)\in\{-1,0,1\}$ and
$\sum \frac12(\eps'_j-\eps_j)\whp_j$ is a multiple of~$P$, thus
of~$p_i$ too,
but $p_i \mid \whp_j$ for $j\neq i$;
thus $\frac12(\eps'_i-\eps_i)\whp_i$ is a multiple of~$p_i$ and,
since the $p_j$'s are pairwise coprime,
this implies $\eps'_i=\eps_i$.

Therefore, $\tp{\sN}$ induces a bijection from~$E$ to its range;
we will denote its inverse by
\begla
\tp{\sN}\ii \col \tp{\sN}(E) \to E.
\edla
Our first step towards the proof of Proposition~\ref{propPsiassumpts} is

\begin{lemma}   \label{lemchitichi}
  Let us define a subset of~$\Z$,
\begla
\fS := rP+\sN(E)+2P\Z =
\big\{\, m\in\Z \mid \tp{m-rP} \in \tp{\sN}(E) \,\big\},
\edla
and a function $\ue \col \fS \to E$,
\beglab{eqdefueonfS}
m\in\fS
\Imp \ue(m) := \tp{\sN}\ii\big(\tp{m-rP}\big).
\edla
Then the coefficients~$\ti\chi(m)$ defined for $m\ge m_0$
by~\eqref{eq:Taylor} satisfy \eqref{eqtichisummjchijti}--\eqref{eqtichisummjchijEXCEPT}
%
%
with a function $\chi \col \Z \to \Z$ defined as follows:
\begin{itemize}
\item
  if $m\in\fS$, then 
\beglab{eqdefchimell}
\chi(m) := \frac{(\ell+1)\cdots(\ell+r-3)}{(r-3)!} \eps_1 \cdots\eps_r
\quad \text{with}\ens
    \ell := \frac{m - (r-2)P - \sN(\ue)}{2P}
    \edla
    %
    %
    and $\ue :=\ue(m)$ (note that $\ell\in\Z$),
\item
  if $m\notin\fS$, then $\chi(m) := 0$.
\end{itemize}
\end{lemma}

For example,~\eqref{eqformulamzero} amounts to
\beglab{eqmzeroinfS}
m_0 = \sN(-\uo)+(r-2)P
\quad\text{with}\ens
\uo:=(1,\ldots,1)\in E,
\edla
hence $m_0\in\fS$ with $\ue(m_0)=-\uo$, and~\eqref{eqdefchimell} yields $\ell=0$ and
$\chi(m_0) = (-1)^r$.

\begin{proof}[Proof of Lemma~\ref{lemchitichi}]
  We can rewrite~\eqref{eq:Taylor} as
  \beglab{eqoriggenser}
  \sum_{m=m_0}^{\infty} \tilde\chi(m) z^m =
  (z^{-P}-z^P)^{-(r-2)} \prod_{j=1}^r (z^{\whp_j}-z^{-\whp_j}).
  \edla
  The second factor in the \rhs\ is
  \beglab{eqcomputgenpol}
    \prod_{j=1}^r \,\sum_{\eps\in\{\pm1\}} \eps z^{\eps\whp_j} =
  \sum_{\ue\in E} \,\prod_{j=1}^r \eps_j z^{\eps_j\whp_j} =
  \sum_{\ue\in E} \, \eps_1\cdots \eps_r \, z^{\sN(\ue)}.
\edla
The first one is
\begin{gather} \notag
  z^{(r-2)P} (1-z^{2P})^{-(r-2)} =
  \sum_{\ell\ge0} \frac{(\ell+1)\cdots(\ell+r-3)}{(r-3)!} z^{2\ell
    P+(r-2)P}, \\[-1.75ex] \intertext{therefore}
  \begin{split}
     \sum_{m=m_0}^{\infty} \tilde\chi(m) z^m &=
  \sum_{(\ue,\ell) \in E\times \Z_{\ge0}}
  \frac{(\ell+1)\cdots(\ell+r-3)}{(r-3)!} \eps_1\cdots \eps_r \,
  z^{\mu(\ue,\ell)} \\[1ex]
  & \qquad\qquad\qquad\qquad\text{with}\ens
  \label{eqdefmu}
  \mu(\ue,\ell):=\sN(\ue)+(r-2)P+2\ell P.
  \end{split}
\end{gather}

Notice that~$\mu$ induces a bijection $E\times\Z \xrightarrow{\sim}
\fS$ whose inverse is explicitly given by
\[
  m\in \fS \imp
  \mu\ii(m) = \big(\ue(m),\ell(m) \big) \ens \text{with}\ens
  %
  %
  %
  \ell(m) :=  \frac{m - (r-2)P - \sN(\ue)}{2P}
\]
(with reference to~\eqref{eqdefueonfS} for the first component).
Thus
\beglab{eqidentitysofar}
\sum_{m\ge m_0} \ti\chi(m) z^m =
\sum_{m\in\fS \, \text{s.t.}\,\ell(m)\ge0} \chi(m) z^m.
\edla
Given an arbitrary $m = \mu(\ue,\ell)\in\fS$ such that $\chi(m)\neq0$,
we note that, if $r\ge4$, then necessarily
$\ell\notin \{-(r-3),-(r-4),\ldots,-1\}$ (so as to have
$(\ell+1)\cdots(\ell+r-3)\neq0$), hence in all cases
we have the alternative
\begin{itemize}
\item either $\ell\ge0$ and $m\ge \sN(\ue)+(r-2)P \ge m_0$,
  \item or $\ell\le-(r-2)$ and $m\le \sN(\ue)-(r-2)P \le - m_0$
  \end{itemize}
  (we have used that $\min\sN +(r-2)P =\sN(-\uo)+(r-2)P = m_0$
  and $\max\sN -(r-2)P =\sN(\uo)-(r-2)P = -m_0$,
by~\eqref{eqmzeroinfS}).

If $(p_1,\ldots,p_r)\neq(2,3,5)$, then $m_0\ge1$ and~\eqref{eqidentitysofar} immediately
yields~\eqref{eqtichisummjchijti}.

The case $(p_1,\ldots,p_r)=(2,3,5)$ requires special treatment.
We then have $m_0=-1$, $P=30$ and $\sN(\ue) =
15\eps_1+10\eps_2+6\eps_3$.
When comparing the \rhs\ of~\eqref{eqidentitysofar} and
$\sum\limits_{m\in\fS \, \text{s.t.}\,m\ge1} \chi(m) z^m$,
we see that
\begin{itemize}
\item all the terms in the former series are found in the latter
except the one with $m=\mu(-\uo,0)=-1$, i.e.\ the term $-z\ii$
(because $\ell\ge0$ and $\sN(\ue) +30+60\ell<1$ entail
$(\ue,\ell)=(-\uo,0)$),
\item all the terms in the latter series are found in the former
except the one with $m=\mu(\uo,-1)=1$, i.e.\ the term~$z$
(because $\sN(\ue) +30+60\ell\ge1$ and $\ell<0$ entail
$(\ue,\ell)=(\uo,-1)$).
  \end{itemize}
  We thus get $\sum\limits_{m\ge m_0} \ti\chi(m) z^m +z\ii=
  \sum\limits_{m\in\fS \, \text{s.t.}\,m\ge1} \chi(m) z^m-z$,
which amounts to~\eqref{eqtichisummjchijEXCEPT}.
\end{proof}

\begin{proof}[Proof of Proposition~\ref{propPsiassumpts}]
It only remains to prove that the function~$\chi$ of Lemma~\ref{lemchitichi} can be
put in the form of the \rhs\ of~\eqref{eqtichisummjchij}.

Note that the set~$\fS$, which contains the support of~$\chi$, is invariant mod~$2P\Z$.
The definition of $\chi(m)$ 
involves the function $\ue(m)$ of~\eqref{eqdefueonfS},
which clearly is an odd $2P$-periodic function of~$m$, and $\ell(m):=\ell$ given
by the rightmost equation in~\eqref{eqdefchimell}.

Let us introduce the positive integer coefficients $\sig(n,k)$ as the
coefficients of the polynomial
\beglab{eqdefsigsk}
(\ell+1)\cdots(\ell+n) = \sum_{k=0}^n \sig(n,k) \ell^k
\quad \text{for any $n\in\Z_{\ge1}$.}
\edla
For $m\in \fS$, we now compute $\chi(m)$:
with the notation $\ell:=\ell(m)$, $\ue:=\ue(m)$ and
\begla
\pi(\ue) := \eps_1\cdots \eps_r,
\edla
plugging the rightmost equation of~\eqref{eqdefchimell}
into~\eqref{eqdefsigsk} and then the result into the leftmost equation
of~\eqref{eqdefchimell}, we get
\begin{align}
\notag  \chi(m) &=
\frac{\pi(\ue)}{r'!} \sum_{0\le k \le r'} \sig(r',k)\ell^k
 &&\hspace{-17em}
     \text{with notation $r':=r-3$}\\[1ex]
\notag          &= \frac{\pi(\ue)}{r'!} \sum_{j+s+t\le r'} \sig(r',j+s+t)
            \frac{(j+s+t)!}{j!s!t!(2P)^{j+s+t}} (-1)^{s+t} m^j((r-2)P)^t \big(\sN(\ue)\big)^s &&\\[2ex]
  \label{eqdefchij} &= \sum_{0\le j \le r-3} m^j \chi_j(m)
                      \qquad \text{with} \ens
\chi_j(m) := -\sum_{0\le s \le r-3-j} C_{j,s}\, \pi\big(\ue(m)\big) \big(\sN(\ue(m))\big)^s,\\[1ex]
\text{where} &
\label{eqdefCjs}
\quad C_{j,s} := \sum_{0\le t \le r-3-j-s} (-1)^{s+t+1} \sig(r-3,j+s+t)
\frac{(j+s+t)! (r-2)^t}{2^{j+s+t}(r-3)!j!s!t! P^{j+s}}.
\end{align}
We thus define the functions  $\chi_0,\ldots,\chi_{r-3}$ by
\begla
\chi_j(m):=0 \ens\text{if}\ens m\notin\fS, \quad
\text{ $\chi_j(m)$ as in~\eqref{eqdefchij}} \ens\text{if}\ens m\in\fS.
\edla
%

To conclude the proof, we just need to check that the~$\chi_j$'s are all odd or even, of same parity
as $j+1$.
We observe that $\chi=\chi_{1,1}$, where $\chi_{x,y}(m) \in \Q[x,y]$
is defined by $m\notin\fS\iimp \chi_{x,y}(m):=0$ and
\begin{multline*}
m\in\fS \imp
\chi_{x,y}(m) :=
\frac{(\ell_{x,y}(m)+1)\cdots(\ell_{x,y}(m)+r-3)}{(r-3)!} \pi(\ue(m)) \\[1ex]
    %
\text{with}\ens
    \ell_{x,y}(m) := \frac{ x m - y \sN(\ue(m))}{2P} - \frac{r-2}{2}.
%
    %
\end{multline*}
The very same computation as above gives
\beglab{eqexpandchixy}
m\in\fS \imp
\chi_{x,y}(m) = -\sum_{0\le j \le r-3} \, \sum_{0\le s \le r-3-j}
C_{j,s}\, m^j \pi(\ue(m)) \big(\sN(\ue(m))\big)^s
x^j y^s.
\edla
Given $m\in\fS$, we have $\ue(-m)=-\ue(m)$, whence
$\pi(\ue(-m))=(-1)^r\pi(\ue(m))$ and
$\ell_{x,y}(-m) = -\frac{ x m - y \sN(\ue(m))}{2P} -
\frac{r-2}{2} = -\ell_{x,y}(m) - (r-2)$,
and this implies that~$\chi_{x,y}$ is an odd function of~$m$ for each $(x,y)$.
Therefore, by~\eqref{eqexpandchixy}, each function
$m \mapsto C_{j,s} \, m^j \,\pi(\ue(m)) \big(\sN(\ue(m))\big)^s$ must
be odd.
Since $m\mapsto m^j \pi(\ue(m)) \big(\sN(\ue(m))\big)^s$ has the
same parity as $j+r+s$, this implies that
\begla
C_{j,s}=0 \quad \text{whenever} \quad \text{$j+s+r $ is even.}
\edla
The result now follows from~\eqref{eqdefchij}.
\end{proof}

\begin{remark}   \label{remdefcT}
  We find it convenient to express $\ue=\ue(m)$ of~\eqref{eqdefueonfS}
  in terms of
\beglab{eqdefCNuo}
\cN := P + \sN \col E\to\Z,
\quad \text{for which\; $\tp{\cN} \col
  E\to\Z/2P\Z$ is odd and injective}
\edla
(note that~$\cN$ is not odd).
  Defining $\cT \col \Z\to\Z$ by
  \beglab{eqdefcTm}
  \cT(m) := m \ens\;\text{if $r$ is odd}, \qquad
  \cT(m) := m-P \ens\;\text{if $r$ is even}
  \edla
  for all $m\in \Z$,
  what we have found is equivalent to\footnote{%
    Indeed, one easily checks that
\begin{gather*}
  \label{equerodd}
\text{$r$ odd} \imp \fS = \cN(E)+2P\Z, \ens\text{and}\ens
\ue(m) = \tp{\cN}\ii\big(\tp{m}\big)
\ens \text{for $m\in\fS$}
\\[1ex]
 %
 %
 %
  \label{equereven}
\text{$r$ even} \imp \fS = P+\cN(E)+2P\Z, \ens\text{and}\ens
\ue(m) = \tp{\cN}\ii\big(\tp{m-P}\big)
\ens \text{for $m\in\fS$.}
%
 %
 %
\end{gather*}
}
  \begin{gather}
\label{eqchijmsfuo}
      \chi_j = \sum_{0\le s \le r-3-j\atop s\equiv r-1-j \ [2]}
      C_{j,s} \, m^s f^\uo \circ\cT,\\[-1.5ex]
      \intertext{with $C_{j,s}\in\Q$ defined by~\eqref{eqdefCjs}
          and $m^s f^\uo$ $2P$-periodic function of support}
\label{eqsupportsSuo}
\fS^\uo := \cN(E) + 2P\Z \, \subset\, \Z\\[.5ex]
  \intertext{and same parity as $r+s$ defined by}
\label{eqdefmsfuo}
  m\in\fS^\uo \imp
  m^s f^\uo(m) := -\pi(\ue)\big(\sN(\ue)\big)^s
\quad \text{with}\ens
\ue = \tp{\cN}\ii\big(\tp{m}\big).
\end{gather}
The function $f^\uo:=m^0 f^\uo$ is one of the ``Hikami
functions'', as we call them with reference to~\cite{K05a},
and the functions~$m^s f^\uo$ are examples of what we call
$s$-Hikami functions and study in greater generality in Appendix~\ref{secsHikfcns}.
%
%
\end{remark}


\subsection{Proof of Theorem~\ref{thm:QM}}
%

We now apply the theory of Section~\ref{sec:QMresurgence} to
$\Theta(\cdot\,;j,\chi_j,2P)$ for each~$j$.

\begin{proposition} \label{prop:vanishingmeanvalue}
  Let $j \in \{0,\ldots,r-3\}$.
  The quantum set $\sQ_{\chi_j,2P}$ as defined in~\eqref{eqdefsQfM} is
  all of~$\Q$,
  i.e.\ the periodic function
  $m\in\Z\mapsto 
  \chi_j(m) e^{i\pi m^2\al/(2P)}$
  has mean value zero
  and the non-tangential limit
    $\lim\limits_{\tau\to\al} \Theta(\tau;j,\chi_j,2P)$ thus exists
  for each $\alpha \in \Q$.
  Moreover,
  \beglab{eqThetaalplusTmedsum}
  \Theta(\al+T;j,\chi_j,2P) = \frac12 \big( \bS^{\frac\pi2-\epsilon} + \bS^{\frac\pi2+\epsilon} \big)
  \wt\Theta_{j,\chi_j,\alpha,2P}(T)
  \quad \text{for all $T\in\HH$,}
  \edla
where the resurgent formal series
$\wt\Theta_{j,\chi_j,\alpha,2P}(T)$ is defined as in~\eqref{eq:ResurgentSeries}.
\end{proposition}

\begin{proof}
  It is proved in Corollary~\ref{corcVs} in the Appendix that,
    for each $s\in\{0,\ldots,r-3\}$, the function $m^s f^\uo\circ\cT$
    belongs to a vector space~$\cVs$, all of whose elements~$g$
    have the property $\sQ_{g,2P}=\Q$ according to
    Proposition~\ref{propConseqTh}.
%
%
By~\eqref{eqchijmsfuo}, this implies $\sQ_{\chi_j,2P} = \Q$.
  %
  %
The rest of the statement follows by the results of Sections~\ref{secptsLt}--\ref{sec:resurgence}.
\end{proof}
    
Here,~$T$, which plays the role of an indeterminate in the formal
series $\wt\Theta_{j,\chi_j,\alpha,2P}(T)$ and of the
resurgence-summability variable in~\eqref{eqThetaalplusTmedsum}, can
be interpreted as a new variable
\begla
T=\tau-\al,
\edla
in accordance with the
notations of Theorem~\ref{thm:QM}(iii).
In view of~\eqref{eqPsicQsumThchij}, we obtain
\begin{corollary} \label{cor:wtPsial}
  For each $\al\in\Q$, the modified GPPV invariant satisfies
  \beglab{eqPsialTSwtPsial}
  \Psi(\al+T)= \frac12 \big( \bS^{\frac\pi2-\epsilon} + \bS^{\frac\pi2+\epsilon} \big)
  \wt\Psi_\al(\al+T)
  \quad \text{for all $T\in\HH$,}
  \edla
  with a resurgent formal series
\begin{equation}  \label{def:Zalpha}  
  \wt\Psi_\al(\al+T) := \sum_{j=0}^{r-3} \wt\Theta_{j,\chi_j,\alpha,2P}(T)
  + \sum_{m\ge0} \frac1{m!}\cQ^{(m)}(\al)T^m
  \in \C[[T]]
\end{equation}
that is Borel summable in all directions except $\pi/2$ and has a
Borel transform all of whose singular points
are of the form $\frac{i\pi m^2}{2P}$, $m\in\Z_{\ge1}$.
%
%
(Recall that $\cQ=0$ except when
    $(p_1,\ldots,p_r)=(2,3,5)$, in which case $\cQ(\tau) = -e^{-i\pi\tau/60}-e^{i\pi\tau/60}$.)
\end{corollary}

The formal series~$\wt\Psi_\al$ that we just defined are those
mentioned in Point~(iii) of Theorem~\ref{thm:QM}.
We go from them to the formal series $\wt Z^*_{\ze}(q)$
of~\eqref{eq:HatAlpha} by a further change of variable
\begla
Q := q-\ze = e^{2\pi i\tau}-\ze =
\ze ( e^{2\pi i T}-1 )
\edla
and multiplication by an analytic function:

\begin{proposition}
Let~$\ze$ be a root of unity and pick any $\al\in\Q$ such that
$\ze=e^{2\pi i\al}$.
The formula
\begla
\wt Z^*_\ze(\ze+Q) := e^{-\frac{i\pi m_0^2\al}{2P}}
\, ( 1 + \ze\ii Q)^{-\frac{m_0^2}{4P}} \,
\wt\Psi_\al\big( \al + \frac{1}{2\pi i} \log ( 1 + \ze\ii Q) \big)
\in \C[[Q]]
\edla
defines a formal series that does not depend on~$\al$ but only
on~$\ze$,
is resurgent in~$Q$, has a Borel transform all of whose singular
points are of the form $-\pi^2 m^2\ze/P$, $m\in\Z_{\ge1}$,
and is Borel summable in all directions except
$\tht_\al:=2\pi\al+\pi$ mod $2\pi$.
The normalized GPPV invariant satisfies
  \begla
  Z^*(\ze+Q)= \frac12 \big( \bS^{\tht_\al-\epsilon} + \bS^{\tht_\al+\epsilon} \big)
  \wt Z^*_\ze(\ze+Q)
  \quad \text{for all $Q\in\C$ such that $\ze+Q\in\DD$.}
  \edla
\end{proposition}

\begin{proof}
  Let us use lightened notation
  \beglab{eqlightPsialT}
  \Psi(\al+T) = \bS^{\frac\pi2}\median \ti\psi(T),
  \qquad
  \ti\psi(T) := \wt\Psi_\al(\al+T)
  \edla
  for the identity~\eqref{eqPsialTSwtPsial}.
  On the one hand, according
  to~\eqref{partialthetaseriesPsist}--\eqref{partialthetaseries}, the
  normalized GPPV invariant at $q=\ze+Q$ can be retrieved from the modified GPPV
  invariant by the formula
  \begla
Z^*(\ze+Q) = e^{-\frac{i\pi m_0^2\al}{2P}} e^{-\frac{i\pi m_0^2 T}{2P}} \,\Psi(\al+T),
\edla
where $2\pi i(\al+T)$ is the branch of $\log(\ze+Q)$ that is close to
$2\pi i\al$ when~$Q$ is close to~$0$.
This means that $2\pi i T$ is the principal branch of $\log(1+\ze\ii
Q)$ and we get
\beglab{eqZstzepQ}
Z^*(\ze+Q) = e^{-\frac{i\pi m_0^2\al}{2P}}
(1+\ze\ii Q)^{-\frac{m_0^2}{4P}}
\,\Psi\big( \al + \tfrac1{2\pi i}\log(1+\ze\ii Q)\big)
\edla
(using the principal branches of the analytic functions
$(1+x)^{-\frac{m_0^2}{4P}}$ and $\log(1+x)$).

On the other hand, if we define a formal series
\beglab{eqdeftiphQ}
\ti\ph(Q) := e^{-\frac{i\pi m_0^2\al}{2P}} (1+\ze\ii Q)^{-\frac{m_0^2}{4P}}
\,\ti\psi\big( \tfrac1{2\pi i}\log(1+\ze\ii Q)\big)
\in\C[[Q]]
\edla
(using the formal series
$(1+x)^{-\frac{m_0^2}{4P}},\log(1+x)\in\C[[x]]$),
then
\begin{quotation}
\noindent $\ti\ph(Q)$ is resurgent, the singular
points of its Borel transform are of the form
$-\pi^2 m^2\ze/P$, $m\in\Z_{\ge1}$,
and it is Borel summable in all directions except~$\tht_\al$ with
\beglab{eqbSthtmedtiph}
\bS^{\tht_\al}\median\ti\ph(Q) = e^{-\frac{i\pi m_0^2\al}{2P}} (1+\ze\ii Q)^{-\frac{m_0^2}{4P}}
\,(\bS^{\pi/2}\median\ti\psi)\big( \al + \log(1+\ze\ii Q)\big).
\edla
\end{quotation}

\noindent
Indeed, we can write
\begla
\ti\psi\big( \tfrac1{2\pi i}\log(1+\ze\ii Q)\big) = 
\ti\psi_1\big(h(Q)\big) \ens \text{with}\ens
h(Q) := \ze \log(1+\ze\ii Q),
\ens
  \ti\psi_1(Q) := \ti\psi\big( \tfrac1{2\pi i\ze} Q\big),
  \edla
  the properties of~$\ti\psi(T)$ indicated in
Corollary~\ref{cor:wtPsial} then trivially entail that
$\ti\psi_1(Q)$
is resurgent, with all singular
points of its Borel transform of the form
$2\pi i\ze\cdot \frac{i\pi m^2}{2P} = -\pi^2m^2\ze/P$, $m\in\Z_{\ge1}$
(because $\wh\psi_1(\xi) = \tfrac1{2\pi i\ze}\wh\psi\big( \tfrac1{2\pi i\ze} \xi\big)$),
and Borel summable in all directions except~$\tht_\al$ with
$(\bS^{\tht}\median\ti\psi_1)(Q) = (\bS^{\pi/2}\median\ti\psi)\big( \tfrac1{2\pi i\ze}Q\big)$;
the point is that these properties of~$\ti\psi_1$ carry over through the
composition with $h(Q)$ and the
multiplication by $e^{-\frac{i\pi m_0^2\al}{2P}} (1+\ze\ii
Q)^{-\frac{m_0^2}{4P}}$ because these are convergent formal series
in~$Q$ and~$h$ is tangent to identity\footnote{%
  See e.g.\ \cite[Thm~5.55]{MS16} for the summability statement and
  \cite[Sec.~6.2 and proof of Thm~6.32]{MS16} for the resurgence statement,
  with the caveat that, in the standard terminology, the
  resurgence-summability variable is $z:=Q\ii$, hence it is
  composition with
  $\frac{1}{h(z\ii)} = z + \text{convergent power series in~$z\ii$}$
    that must be considered (cf.\ \cite[Sec.~5.15]{MS16}).}.

  Now, the \rhs\ of~\eqref{eqbSthtmedtiph} coincides with that
  of~\eqref{eqZstzepQ} thanks to~\eqref{eqlightPsialT}.
  Since~\eqref{eqdeftiphQ} is equivalent
  to~\eqref{eqdefwtPsifromwtZst} with $\ti\ph(Q)=\wt Z^*_\ze(\ze+Q)$, the proof is complete.
  \end{proof}





This yields Point~(i) and the beginning of of Point~(iii) of
Theorem~\ref{thm:QM} until property~\eqref{eq:resurgentexpansion},
with the reinforcement indicated in Remark~\ref{remstrongerf}.

Let us now prove the rest of Point~(iii) of Theorem~\ref{thm:QM}. We
need
\begin{proposition}   \label{propsuppchijmz}
All the functions~$\chi_j$, $j\in\{0,\ldots,r-3\}$, in the
decomposition~\eqref{eqtichisummjchij} fulfill
condition~\eqref{eqsupportcondition} of Theorem~\ref{thmqmpts} and
Corollary~\ref{corollarypartialthetaserieshigherdepth}
with $M=2P$ and
$n_0=m_0$, where~$m_0$ is as
  in~\eqref{eqformulamzero}.
  %
%
\end{proposition}

\begin{proof}
  In view of \eqref{eqdefcTm}--\eqref{eqdefmsfuo},
  %
  %
  it is sufficient to prove that the condition is fulfilled by all the
  functions
$m^sf^\uo\circ\cT$, $s\ge0$. 
They all have the same support, $\cT\ii\big(\cN(E)\big) + 2P\Z$.
It is thus sufficient to check that the function
$\fourp{(\cT\ii\circ\cN)^2} \col E \to \Z/4P\Z$ is constant, taking on
the value $\fourp{m_0^2}$.

Observe that the value of $\fourp{(\cT\ii\circ\cN)^2}$ at $\ue=-\uo$
is $\fourp{m_0^2}$. This is because~\eqref{eqmzeroinfS} yields
$m_0 = \cN(-\uo)+(r-3)P$,
whence
\begla
\cT\ii\circ\cN(-\uo) =  \begin{cases}
  m_0 - (r-3)P  & \text{ if $r$ is odd} \\[.8ex]
  m_0 - (r-4)P  & \text{ if $r$ is even}
		 \end{cases} 
%
%
\edla
and  in all cases $\cT\ii\circ\cN(-\uo) = m_0 \!\mod 2P$, thus
$(\cT\ii\circ\cN)^2(-\uo) = m_0^2 \!\mod 4P$.

We conclude the proof by showing that $\fourp{(\cT\ii\circ\cN)^2}$ is constant.
For any $\ue,\ue'\in E$,
taking the square of $\sN$ as defined by~\eqref{eqdefwhpj}, we get
\begla
\big( \sN(\ue') \big)^2 - \big( \sN(\ue) \big)^2 =
2 \sum_{i<j} (\eps'_i\eps'_j-\eps_i\eps_j) \whp_i \whp_j,
\edla
whence $\big( \sN(\ue') \big)^2 - \big( \sN(\ue) \big)^2 \in 4P\Z$
because $\eps'_i\eps'_j-\eps_i\eps_j \in \{-2,0,2\}$ and $\whp_i \wh
p_j\in P\Z$
(since~$\whp_j$ is a multiple of~$p_i$),
i.e.\ $\fourp{(\sN)^2} \col E \to \Z/4P\Z$ is constant.
The function $\fourp{(\cN+P)^2}=\fourp{(\sN+2P)^2}$ is thus
constant too, and this is nothing but $\fourp{(\cT\ii\circ\cN)^2}$
when~$r$ is even.

Another implication of our previous computation is that
$\fourp{\cN^2+2P\cN}$ is constant.
But $\ttwo{\cN}$ is constant too, since
$\cN(\ue')-\cN(\ue) = \sum(\ue'_j-\ue_j)\whp_j \in 2\Z$ (because
$\ue'_j-\ue_j$ is always even),
therefore $\fourp{\cN^2}$ is constant, and this is $\fourp{(\cT\ii\circ\cN)^2}$
when~$r$ is odd.
\end{proof}

We can thus apply Theorem~\ref{thmqmpts} and
Corollary~\ref{corollarypartialthetaserieshigherdepth}
to all the functions~$\chi_j$;
in view of~\eqref{eqPsicQsumThchij}, this yields quantum
modularity for~$\Psi(\tau)$ on the congruence subgroup $\Ga_1(4P)$.
As for the part of Point~(iii)
%
%
relative to a vector-valued quantum modular form on the full modular
group $\SL(2,\Z)$, this follows from~\eqref{eqPsicQsumThchij}, \eqref{eqchijmsfuo},
and Corollary~\ref{corcVs} and Proposition~\ref{propositionvectorvalue} in the appendix.





\smallskip

Finally, to conclude the proof of Theorem~\ref{thm:QM}, we just need
to prove Point~(ii), which is equivalent to


\begin{proposition}   \label{propradiallimitPsial}
\noindent \emph{\textbf{(i)}}
The \lhs\ of~\eqref{eqdefnstar} is an odd integer.

\noindent \emph{\textbf{(ii)}}
  For each $\al\in\Q$, the constant term $\Psi_{\al,0}$ of the formal
  series $\wt\Psi_\al(\al+T)$ of~\eqref{def:Zalpha} is
  \beglab{eqradiallimitPsial}
  \Psi_{\al,0} = 2 (-1)^r
  e^{\pi i\al\phi/2} (e^{\pi i\al}-e^{-\pi i\al}) \WRT(X,e^{2\pi i\al})
  \edla
  with~$\phi$ as in~\eqref{eqdefphi}.
\end{proposition} 


\begin{proof}[Proof that Proposition~\ref{propradiallimitPsial} is
  equivalent to Theorem~\ref{thm:QM}(ii)]
Point~(i) of Proposition~\ref{propradiallimitPsial} makes it possible
to define the integer~$n_*$ that is involved in the
formula~\eqref{eq:Zhatradiallimit} that we now prove.

Let $\ze\in\gR$ and pick any $\al\in \Q$ such that $\ze:=e^{2\pi i\al}$.
  The relation~\eqref{eqdefwtPsifromwtZst} between the formal series
  $\wt\Psi_\al(\al+T)$ and $\wt Z^*_\ze(\ze+Q)$ implies that their
  constant terms are related by
  \begla
  \Psi_{\al,0} = e^{\frac{i\pi m_0^2\al}{2P}} Z^*_{\ze, 0}.
  \edla
  Formula~\eqref{eqradiallimitPsial} is thus equivalent to
  \begla
  Z^*_{\ze, 0} = 2 (-1)^r e^{2\pi i\al\Om} (\ze-1) \WRT(X,\ze)
\quad \text{with}\ens
\Om := - \frac{m_0^2}{4P} +\frac{\phi}4 - \frac12.
  \edla
  On the one hand, \eqref{eqdefphi} yields
  \begla
  -\phi - 24\la = P \Big(r-2 - \sum_{1\le j\le r} \frac1{p_j^2} \Big).
  \edla
  On the other hand, \eqref{eqformulamzero} yields
  \begla
  \frac{m_0^2}{P^2} = \Big( r-2 - \sum_{1\le i\le r}\frac1{p_i} \Big)^2
  = (r-2)^2 - 2(r-2)\sum_{1\le i\le r}\frac1{p_i}
  + \sum_{1\le i\le r}\frac1{p_i^2}
  + 2\sum_{1\le i<j\le r}\frac1{p_i p_j}.
  \edla
  Recall the notation $\whp_i:= \frac P{p_i}$ from~\eqref{eqdefwhpj}.
  We also introduce the notation $\whp_{i,j} := \frac P{p_i p_j}$ for
  $i<j$.
  We thus find
  \beglab{eqmphipmzsqP}
  -\phi - 24\la + \frac{m_0^2}P = (r-2)(r-1)P
  - 2(r-2)\sum_i\whp_i + 2\!\sum_{i<j}\whp_{i,j}
  = -2(2n_*+1)
  \edla
  by~\eqref{eqdefnstar},
  whence $\Om = n_*-6\la \in\Z$
  and the conclusion follows from the identity $Z^*_{\ze, 0} = 2 (-1)^r \ze^{\Om} (\ze-1) \WRT(X,\ze)$.  
\end{proof}

\begin{proof}[Proof of Proposition~\ref{propradiallimitPsial}]
\noindent {\textbf{(i)}}
With the previous notation the \lhs\ of~\eqref{eqdefnstar} can be
written
\begla
-\frac{(r-1)(r-2)}2 P + (r-2) \sum_{1\le i \le r} \whp_i
- \sum_{1\le i < j \le r} \whp_{i,j}.
\edla
Recall that, since the beginning, we have assumed $p_2,\ldots,p_r$ odd,
i.e.\ $p_i\equiv 1\!\mod 2$ for $i>1$,
and we can go on computing $\!\!\mod 2$:
\begin{gather*}
  \whp_1 \equiv 1, \qquad
  \whp_{1,j} \equiv 1 \ens\text{for $j>1$,}
  \\[.8ex]
  P\equiv p_1, \qquad
  \whp_i \equiv p_1 \ens\text{for $i>1$,} \qquad
  \whp_{i,j} \equiv p_1 \ens\text{for $1<i<j$,}
  \end{gather*}
  whence it follows that the \lhs\ of~\eqref{eqdefnstar} $\!\!\mod 2$
  is
  \begin{multline*}
  -\frac{(r-1)(r-2)}2 p_1 + (r-2) \big( 1 + (r-1)p_1 \big)
  - \sum_{2\le j\le r} 1
  - \sum_{2\le i < j \le r} p_1 \\[.8ex]
  \equiv \frac{(r-1)(r-2)}2 p_1 + r-2 - (r-1) - \frac{(r-1)(r-2)}2 p_1
  \equiv1.  
\end{multline*}

\smallskip

\noindent {\textbf{(ii)}}
We now prove~\eqref{eqradiallimitPsial} for an arbitrary $\al\in\Q$.

The case $\al=\ell\in\Z$ has been taken care of in Remark~\ref{remcasekis1}:
we then have $\Psi_{\ell,0}=0$ as desired, by~\eqref{eqlimellPsizero}.

Suppose now that $\al = \ell/k$ with coprime integers~$\ell$ and~$k$
such that $k\ge2$.
%
%
  Recall the notation $\zeta_k=e^{2\pi i/k}$. We have
  $\Gal(\Q(\zeta_k):\Q)\cong (\Z/ k\Z)^{\times}$, where the Galois
  transformation $\sigma_{k,u}$ associated with
  $u\in (\Z/ k\Z)^{\times}$ is the field automorphism defined by
$\sigma_{k,u}\cdot\zeta_k=\zeta_k^u$ and by $\Q$-linearity.
  
    By~\eqref{def:Zalpha},
    %
    %
 the property~\eqref{eqsupportcondition} of the support of the
    functions~$\chi_j$ stated in Proposition~\ref{propsuppchijmz}
    %
%
    and~\eqref{eq:Lseriesevaluation},
  we have that
  \begin{equation} \label{eq:integrality}
    e^{-2\pi i m_0^2 \ell / (4Pk)}\Psi_{\al,0} = Z^*_{e^{2\pi i\al}, 0} \in \Q(\zeta_k).
    %
    %
    %
\end{equation}
%
%
%

Assume first that $\ell$ is odd. Then $\gcd(\ell,4k)=1$ by assumption, and
we can consider the associated Galois transformation
$\sigma_\ell :=\sigma_{4k,\ell}$. We note that the action of $\sigma_{4k,\ell}$
on the image of $\Q(\zeta_k)$ under the natural inclusion
$\Q(\zeta_k) \hookrightarrow \Q(\zeta_{4k})$
%
%
is equal to the action of $\sigma_{k,\ell}$ on $\Q(\zeta_k)$, and we see that
%
%
\begin{equation} \label{eq:ZhatGaloistransformation}
  e^{-2\pi i m_0^2 \ell / (4Pk)}\Psi_{\al,0} =
  \sigma_{\ell}\cdot\!\big( e^{-2\pi i  m_0^2 / (4Pk)} \Psi_{1/k,0}\big).
  \end{equation}
  In view of Corollary~\ref{cor:1/klimit}, we thus get
  \beglab{eqgoonwithcomput}
  e^{-2\pi i m_0^2 \ell / (4Pk)}\Psi_{\al,0} =
  \sigma_\ell\cdot\!\left(2 (-1)^r e^{2\pi i(P\phi- m_0^2) / (4Pk)}
    (\zeta_k^{1/2}-\zeta_k^{-1/2}) \WRTk(X)\right).
  \edla
  Earlier, we have encountered the quantity $\phi - \frac{m_0^2}P$
  and proved~\eqref{eqmphipmzsqP}, which amounts to
\begin{equation} \label{eq:Pdivisibility}
  -P\phi+ m_0^2 = 2P\big(2(6\la-n_*)-1\big)
  \quad \text{(odd multiple of $2P$).}
    \end{equation} 
    This implies that
    $e^{2\pi i(P\Phi- m_0^2) / (4Pk)} \in \Q(\zeta_{4k})$
    is mapped to
    $e^{2\pi i\ell(P\Phi- m_0^2) / (4Pk)}$
    by~$\sig_\ell$, and as also
    $2(\zeta_k^{1/2}-\zeta_k^{-1/2}) \WRTk(X) \in \Q(\zeta_{4k})$, we
    may continue the computation in~\eqref{eqgoonwithcomput} as
    follows (using the fact that~$\sig_\ell$ is a field automorphism):
    \begin{align}
      \begin{split} \label{eq:almostthere}
        e^{-2\pi i m_0^2 \ell /
          (4Pk)}\Psi_{\al,0}&=\sigma_\ell\cdot\!\big(2\,e^{2\pi i(P\phi-
          m_0^2) / (4Pk)} (\zeta_k^{1/2}-\zeta_k^{-1/2}) \big)\,
        \sigma_\ell\cdot\!\big(\WRTk(X)\big) \\&= 2\, e^{-2\pi i m_0^2 \ell /
          (4Pk)} e^{\pi i\al\phi/2}(e^{\pi i\al}-e^{-\pi
          i\al})\WRT(X,e^{2\pi i\al}),
\end{split}
    \end{align}
    where,
    %
    %
    for the last equality, we used the Galois
    equivariance~\eqref{eq:GaloisEquivariance} of WRT invariants.
    Multiplying both sides of~\eqref{eq:almostthere} by
    $e^{2\pi i m_0^2 \ell / (4Pk)}$ gives~\eqref{eqradiallimitPsial}.

    If~$\ell$ is even, then~$k$ must be odd, thus
    $\zeta_{2k} = -(\zeta_k)^{\frac{k+1}2} \in \Q(\zeta_k)$ and $4 \in (\Z/k\Z)^{\times}$. Therefore
    \eqref{eq:Pdivisibility} implies
    \begin{equation}
   e^{2\pi i(P\Phi-  m_0^2) / (4Pk)} \in \Q(\zeta_k)
    \end{equation}
    and the proof goes through as before, except that we apply the
    Galois automorphism $\sigma_{k,\ell}$ directly (instead of applying
    $\sigma_{4Pk,\ell}$ under the embeddding
    $\Q(\zeta_k) \hookrightarrow \Q(\zeta_{4k})$).
    %
    %
  \end{proof}

  This ends the proof of Theorem~\ref{thm:QM}.


  \subsection{The WRT invariant of~$X$ as limit of a median sum}
  \label{secWRTmedsum}

  In the previous subsection, Theorem~\ref{thm:QM}(ii) was proved in
  the form of formula~\eqref{eqradiallimitPsial}, which we saw is
  equivalent to~\eqref{eq:Zhatradiallimit}, and which stems
  from~\eqref{eq1/klimit} in Corollary~\ref{cor:1/klimit}.
  This is the link between the GPPV invariant and the WRT invariant.
  As a preparation for the proof of Theorem~\ref{Thm:main} (to be
  found at the end of the next section), we now
  put together~\eqref{eq1/klimit} and the $\al=0$ case of Corollary~\ref{cor:wtPsial}:


\begin{proposition}   \label{propradiallimitWz}
  Consider the resurgent-summable formal series
  \beglab{eqdefwtWz}
  \wt W_0(\tau) := \cE(\tau) \wt\Psi_0(\tau) / \tau \in \C[[\tau]]
  \quad\text{with}\;\, \cE(\tau)\;\, \text{as in~\eqref{eqdefcE}.}
    \edla
Then $\WRT_k(X)$ can be recovered as a
non-tangential limit at~$1/k$ of the function $\bS^{\pi/2}\median \wt
W_0$ that is holomorphic in~$\HH$:
\beglab{eqWRTklimbSmedWz}
\WRT_k(X) = \lim \bS^{\pi/2}\median \wt W_0(\tau)
\quad \text{as}
 \; \tau\to 1/k \; \text{non-tangentially from within~$\HH$} 
 \edla
for every $k\in\Z_{\ge2}$.
\end{proposition}

Note that the series~$\wt W_0(\tau)$ is nothing but the \rhs\
of~\eqref{eqWzintermsofcEPsiz}.
 Later, at the end of the next section, we will show that~$\wt
 W_0(\tau)$ coincides with the series~$W_0(\tau)$ of~\eqref{eqdefWzerolaX}.

\begin{proof}[Proof of Proposition~\ref{propradiallimitWz}]
The formal series defined by~\eqref{eqdefwtWz} is summable in
the same directions as~$\wt\Psi_0$ and resurgent with the same
location of singularities in the Borel plane,
because~$\wt\Psi_0(\tau)$ is divisible by~$\tau$
(cf.~\eqref{eqlimellPsizero}) and the above properties are preserved
by division by~$\tau$ and mutiplication by~$\cE(\tau)$ (since the
latter is a convergent series).

Now,~\eqref{eq1/klimit} gives the non-tangential limit of~$\Psi(\tau)$
at~$1/k$ in the form
\begla
4i(-1)^r \Big(\!\sin\frac\pi k\Big) e^{\frac{i\phi\pi}{2k}} \WRT_k(X) =
k\ii \cE(1/k)\ii \WRT_k(X).
\edla
Here, we identify the convergent formal series~$\cE(\tau)$
of~\eqref{eqdefcE} with its sum, which is a meromorphic function
in~$\C$ regular in $\C\setminus\Z^*$, because~$1/k$ belongs to its
disc of convergence.

We thus find
\begin{align*}
\WRT_k(X) &= k\cE(1/k) \lim_{\tau\to1/k} \Psi(\tau)
= \lim_{\tau\to1/k} \tau\ii\cE(\tau) \Psi(\tau) \\[.8ex]
&= \lim_{\tau\to1/k} \tau\ii\cE(\tau) \bS^{\pi/2}\median\wt\Psi_0(\tau)
=  \lim_{\tau\to1/k} \bS^{\pi/2}\median\big(\tau\ii\cE(\tau)\wt\Psi_0(\tau)\big)
\end{align*}
(where the last-but-one step is justified by the $\al=0$ case
of~\eqref{eqPsialTSwtPsial}
and the last step results from the compatibility of Borel-Laplace
summation with multiplication) and we are done.
\end{proof}

\begin{remark}   \label{remonwtPsiz}
The formal series $\wt\Psi_0$ was defined in~\eqref{def:Zalpha}.
It can also be written in terms of the the formal series
$\wt\cB_0(\tau):=\sum_{p\ge1} \cB_{0,p} \tau^{p+\frac12}$ 
that appears in~\eqref{eqasymptcBzero} and whose Borel transform is
the explicit meromorphic function~$\cB_0(\xi)$ of~\eqref{eq:Borel}:
indeed, at the beginning of Section~\ref{sec:proofs}, we
proved~\eqref{eqPsiasmedBzero}, which amounts to
\beglab{eqwtPsizaltern}
\wt\Psi_0(\tau) = \frac{(-1)^r e^{\pi i/4}}{\sqrt{2P}} \tau^{-1/2} \wt\cB_0(\tau).
\edla
\end{remark}



  \section{Witten's asymptotic expansion conjecture for Seifert fibered
  homology spheres}
\label{secWAECSFHS}

This section aims at stating and proving Theorem~\ref{thmMflat}, which
has been alluded to in the introduction of this paper, and then proving Theorem~\ref{Thm:main}.

\subsection{The moduli space of flat connections with compact gauge group} \label{sec:Mflat}

The orbifold surface of the Seifert fibered $3$-manifold~$X$ is the
two-sphere with $r$ marked points. Removing from~$X$ a tubular
neighbourhood of
the exceptional fibres 
results in a $3$-manifold naturally
homeomorphic to $\Sigma_{0,r} \times S^1$, where $\Sigma_{0,r}$ is a
two-sphere with $r$ boundaries.
Let $G:=\SU(2)$, and let $C_G$ be the set of conjugacy
classes of $G$.

For each tuple
$C=(C_1,\ldots,C_r) \in (C_G)^r$, 
denote by $\M(\Sigma_{0,r},C)$ 
the moduli space of flat $G$-connections on $\Sigma_{0,r}$ with
holonomy around the $j^{\text{th}}$ boundary component contained in~$C_j$ for
each $j\in \{1,\ldots,r\}$.
It is well-known that the moduli space
$\M(\Sigma_{0,r},C)$ 
is connected (when 
non-empty), and that the subspace given by the moduli space
$\M^{\Irr}(\Sigma_{0,r},C)$ 
of flat irreducible connections is a smooth manifold whose dimension
is known \cite[Sec.~4]{Freed}---see~\eqref{eqdimformula} below.

Denote by $\M^{\Irr}(X)$ the moduli space of irreducible flat
$G$-connections on~$X$. Denote by $T \in \M(X)$ the gauge
equivalence class of the trivial flat $G$-connection. As~$X$ is an
integral homology sphere, we have that
\begin{equation}
\M(X) = \{T\} \sqcup \M^{\Irr}(X).
%
%
\end{equation}
%


For each $\ul= (\ell_1,\ldots,\ell_r)\in \Z^r$, define $C^{(\ul)}=(C_1^{(\ul)},\ldots,C_r^{(\ul)}) \in
(C_G)^r$ by
\begin{equation}
  C_j^{(\ul)} := \;\text{conjugacy class of}\;
  \begin{pmatrix}  e^{\pi i \ell_j/ p_j} & 0 \\ 0 &  e^{-\pi i \ell_j/ p_j} \end{pmatrix}
  %
  %
  \ens\text{for}\ens j=1,\ldots,r.
\end{equation}
%
%
Recall that in \cite[Proposition~6]{AM22} \label{secrefAMprop6}
it is established that there is a one-to-one correspondence between the components of the moduli space
of irreducible flat $\SL(2,\mathbb{C})$-connections on~$X$ and the
elements of the set
\begin{multline}   \label{eqmathfrakL}
  \fL(p_1,\ldots,p_r) := \Big\{\ul\in\Z^r\ \Big| \ 
      0\leq \ell_{j}\leq p_{j}\ \text{for all~$j$,} \\
        \qquad\qquad \frac{\ell_j}{p_j}\notin\Z  \ \text{for at least three
          values of~$j$,}\\
     \ell_j \text{ is even for all $j\ge2$}  \Big\}.
\end{multline}
%
%
We now introduce a subset of~\eqref{eqmathfrakL}, 
which we will prove parametrizes the components of~$\M^{\Irr}(X)$.


\begin{definition} \label{defR}
  %
%
%
%
We set
\begin{align} 
\fR(p_1,\ldots,p_r) &:= \Big\{\ul\in \fL(p_1,\ldots,p_r) \Big| \ 
                      \text{for each subset $J \subset \{1,\ldots,r\}$ of odd cardinality,}
  \notag \\
&\hspace{17.35em}      \sum_{j \in J} \frac{p_j-\ell_j}{p_j}  \, + \,
      \sum_{j \in \{1,\ldots,r\}\setminus J}\, \frac{\ell_j}{p_j} >1. \Big\}
      \label{eq:JeffreyTypeInequality}
\end{align}
\end{definition}

We are now ready to state and prove

\begin{theorem} \label{thmMflat}
  For each tuple $\ul \in \fR(p_1,\ldots,p_r)$ we have that
  $\M(\Sigma_{0,r}, C^{(\ul)})=\M^{\Irr}(\Sigma_{0,r}, C^{(\ul)})$, and
  this moduli space is non-empty.
  Pullback with respect to the
embedding $\iota:\Sigma_{0,r} \hookrightarrow X$ induces a
  homeomorphism
  \begin{equation}
        \M^{\Irr} (X)  \cong \bigsqcup_{\ul \in
          \fR(p_1,\ldots,p_r)} \M^{\Irr}(\Sigma_{0,r}, C^{(\ul)}).
      \end{equation}
    In particular, the set $\pi_0(\M^{\Irr}(X))$ is in bijection with
    $\fR(p_1,\ldots,p_r)$.
\end{theorem}

Let us introduce the notation
\begin{equation}   \label{eqdeftul}
t_\ul := \text{ number of
  indices $j\in \{1,\ldots,r\}$ such that $\ell_j$ is multiple of~$p_{j}$}
\end{equation}
for any $\ul\in\Z^r$;
thus $t_\ul\leq r-3$ for $\ul\in \fR(p_1,\ldots,p_r)$ or
$\fL(p_1,\ldots,p_r)$.
The aforementioned dimension formula from \cite[Sec.~4]{Freed} is
\begin{equation}   \label{eqdimformula}
  \dim\M^{\Irr}(\Sigma_{0,r}, C^{(\ul)}) = 2(r-3-t_\ul)
  \qquad \text{for each} \ens \ul\in\fR(p_1,\ldots,p_n).
\end{equation}



  \begin{remark} Theorem~\ref{thmMflat} builds on the works
    \cite{KK91,FS90}, in which the component labelled by
    $\ul\in \fR(p_1,\ldots,p_r)$ in our notation, was described as
    a so-called admissable linkage, shown to be a closed manifold of
    dimension $2(r-3-t_\ul)$ in \cite{FS90}. The novelty of
    Theorem~\ref{thmMflat} is to use the work \cite{JM05} to describe
    the components in terms of moduli spaces of flat $G$-connections
    on~$\Sigma_{0,r},$ which is a deformation retract of the punctured
    orbifold surface of~$X$, with punctures at the exceptional orbits.
  %
  %
    The utility of Theorem~\ref{thmMflat} is that the condition
    indicated in Definition~\ref{defR} 
    will allow us to parametrize the only contributions to the GPPV
    invariant that may not vanish in the limit
    $q \rightarrow e^{2\pi i/k}$, as proven below.
%
  %
  This will be used 
  in our proof of Theorem~\ref{Thm:main}.
\end{remark}

\begin{proof}[Proof of Theorem~\ref{thmMflat}]
  We begin by recalling the character variety presentations of the
  relevant moduli spaces.
  For each $j=1,\ldots,r$, let $x_j \in \pi_1(X)$ be the
  homotopy class of a small circle in $\Sigma_{0,r} \times \{1\}$
  encircling the $j^{\text{th}}$ boundary component of $\Sigma_{0,r}$,
  these~$r$ circles being connected to a common base point by a star-shaped set of arcs.
  We have the following finite presentations
\begin{align} 
\begin{split} \label{eq:pi_1}
 & \pi_1(\Sigma_{0,r}) \cong \langle x_1,\ldots, x_r \rangle / \langle x_1 \cdots x_r \rangle,
    \\ &\pi_1(X) \cong \langle x_1,\ldots,x_r, h \rangle /R,
    \end{split}
\end{align}
where $R$ is the normal subgroup of $\langle x_1,\ldots,x_r, h \rangle$
generated by $\prod_{j=1}^r x_j$ and the elements $ [x_j,h],$ and
$x_j^{p_j} h^{-q_j}$ for $ j=1,\ldots,r$.
Let $C=(C_1,\ldots,C_r) 
\in (C_G)^r$. Let $I\in G$ be the identity matrix, and let
$Z=\langle -I \rangle$ denote the center of $G$.
Regard $U(1)$ as a subgroup of $G$ through the standard embedding,
defined for all $\zeta \in U(1)$ by
$\zeta \mapsto
\begin{psmallmatrix} \zeta & \ 0 \\ 0 & \ \overline{\zeta}\end{psmallmatrix}$.
%
Recall that a $G$-representation~$\rho$ is
irreducible if and only if the image of $\rho$ is not conjugate to a
subgroup of $U(1)$. By \cite[Lemma~2.1]{FS90} we have that any
representation $\rho:\pi_1(X)\rightarrow G$ must satisfy
$\rho(h) \in Z$. For $[\rho]=T$ this is clear, and for $\rho$
irreducible, we note that the image of $\rho$ is contained in the
centralizer of $\rho(h)$, and if $h$ is not central, this implies that
the image of $\rho$ is conjugate to a subgroup of $U(1)$, and
therefore $\rho$ is reducible. Associating to a flat $G$-connection
the holonomy representation of the first fundamental group induces
bijections
\begin{align} \begin{split} \label{eq:modulispaces}
    \M(\Sigma_{0,r}, C) 
    &\cong \big\{Y \in C_1\times \cdots \times C_r \mid 
    Y_1\cdots Y_r = I \big\}/G, 
    \\[1ex]
    \M(X) &\cong \big\{(H,Y) \in Z \times G^{r} \mid 
    Y_j^{p_j}=H^{q_j} \ \text{for each $ j\in \{1,\ldots,r\}$}
    \ \text{and}\ 
    Y_1\cdots Y_r = I 
     \big\} /G.
     \end{split}
\end{align}

Define $\M(\Sigma_{0,r})=\bigsqcup_{C \in (C_G)^r}
\M(\Sigma_{0,r},C)$. We now analyze the image of $\iota^*: \M(X)
\setminus \{T\} \rightarrow \M(\Sigma_{0,r})$.
Towards that end, let a non-trivial flat $G$-connection on~$X$ be
represented by an irreducible $G$-representation
$\rho: \pi_1(X)\rightarrow G$.
As explained in \cite[Sec.~2]{AM22},\footnote{in the context of
  \emph{complex} Chern-Simons theory, i.e.\ with $G = \SL(2,\C)$
  instead of $\SU(2)$}
we can and will assume that
$q_1$ is odd and $q_j$ is even for $j \in \{2,..,r\}$. For each
$j\in \{1,\ldots,r\}$ set $Y_j=\rho(x_j)$ and set $H=\rho(h)$. For
$j\geq 2$ the fact that $H=\pm I$ together with the relation
$x^{p_j}h^{-q_j}$ implies that $Y_j^{p_j}=I$, as $q_j$ is even, and
therefore the eigenvalues of $Y_j$ are two mutually inverse
$p_j$'th roots of unity.
%
%
Thus there is a uniquely determined even number
$\ell_j \in \{0,\ldots,p_j-1\}$, which is an invariant of the gauge
equivalence class $[\rho]$, such that
 \begin{equation} \label{eq:t1}
   \Tr(Y_j)=e^{\pi i\ell_j/p_j}+e^{-\pi i\ell_j/p_j}
   \quad \text{for $j\ge2$.}
   \end{equation} Similarly, from the equation $Y_1^{p_1}=H^{q_1}=H$, we deduce that there exists a uniquely determined $\ell_1 \in \{0,\ldots,p_1\}$ (not necessarily even) such that
 \begin{equation} \label{eq:t2}
     \Tr(Y_1)= e^{\pi i\ell_1/p_1}+e^{-\pi i\ell_1/p_1}.
 \end{equation}
 Recall that two elements $A,B \in G$ are conjugate if and only if
 $\Tr(A)=\Tr(B)$. Therefore, if we set
 $\ul=(\ell_1,..,\ell_r) \in \Z^r$, then we obtain from~\eqref{eq:t1}--\eqref{eq:t2}
 %
 %
 that $Y_j$ is contained in $C^{(\ul)}_j$ for each $j\in \{1,\ldots,r\}$.
 Therefore we have that
 \begin{equation} \iota^*([\rho]) \in \M(\Sigma_{0,r},C^{(\ul)}).
 \end{equation} By analyzing the presentations of moduli spaces given in~\eqref{eq:modulispaces}, it is straightforward to see that pullback induces a homeomorphism \begin{equation} \label{eq:homeomomorphism} (\iota^*)^{-1}(\M(\Sigma_{0,r},C^{(\ul)})) \rightarrow \M(\Sigma_{0,r},C^{(\ul)}),
\end{equation} 
where the inverse of a flat $G$-connection on $\Sigma_{0,r}$ represented by a homomorphism $\rho': \pi_1(\Sigma_{0,r})\rightarrow G$ is represented by the homomorphism $\rho: \pi_1(X) \rightarrow G$ given by
$\rho(x_j)=\rho'(x_j)$ for $j \in \{1,\ldots,r\}$ and
$\rho(h)=\rho'(x_1)^{p_1}$.
Further, we note that by \cite[Lemma~2.2]{FS90} at most $r-3$ of the
matrices $Y_j$ are equal to $\pm I$, and therefore
$\ul\in\fL(p_1,\ldots,p_r)$
%
%
Indeed, if this was not so, the equation $\prod_{j=1}^r Y_j=I$ would simplify to $Y_{j_1}Y_{j_2}=\pm I$ for some $1\leq j_1< j_2 \leq r$, and by coprimality considerations, this would imply that $Y_{j_1}, Y_{j_2} \in \{\pm I\}$. In particular, we would have that $Y_j \in \{\pm I\}$ for all $j\in \{1,\ldots,r\}$, and this implies the image of $\rho$ is conjugate to a subgroup of $U(1)$, and in particular $\rho$ is reducible. 

We now argue that $\M(\Sigma_{0,r},C^{(\ul)})$ contain only
irreducible connections. Recall that a $G$-representation $\rho$ is
irreducible if and only if the image of $\rho$ is not conjugate to a
subgroup of $U(1)$. From the presentations of the first fundamental
groups given in~\eqref{eq:pi_1} we deduce that for every
$\rho:\pi_1(X)\rightarrow G$ the image of
$\iota^*(\rho):\pi_1(\Sigma_{0,r})\rightarrow G$ is equal to the image
of $\rho$. Since
$(\iota^*)^{-1}(\M(\Sigma_{0,r},C^{(\ul)})) \subset \M^{\Irr}(X)=\M(X)
\setminus \{T\}$, and since \eqref{eq:homeomomorphism} is a
homeomorphism (and in particular surjective), we
obtain \begin{equation}
  \M(\Sigma_{0,r},C^{(\ul)})=\M^{\Irr}(\Sigma_{0,r},C^{(\ul)}).
\end{equation}

Thus it only remains to show that for each $\ul$ as above, the moduli
space $\M(\Sigma_{0,r},C^{(\ul)})$ is non-empty if and only if
$\ul\in \fR(p_1,\ldots,p_r)$. We already noted that
$\ul\in\fL(p_1,\ldots,p_r)$,
%
%
and thus it remains to prove that $\M(\Sigma_{0,r},C^{(\ul)})$ is
non-empty if and only if~$\ul$ satisfies~\eqref{eq:JeffreyTypeInequality}.
%
%
Towards that end, we recall the content of
\cite[Theorem~2.2]{JM05}. For any
$\lambda=(\lambda_j)_{j=1}^r \in [0,\pi]^r$, let
$C^{\lambda} \in (C_G)^r$ be the tuple such that for each
$j\in \{1,\ldots,r\}$ the class $C^{\lambda}_j$ contains the
matrix $\begin{psmallmatrix} e^{i\lambda_j} & 0 \\ 0 & e^{-i\lambda_j}\end{psmallmatrix}$.
%
%
From the character variety presentation~\eqref{eq:modulispaces} we see
that $\M(\Sigma_{0,r},C^{\lambda})$ is non-empty if and only if
$I \in C^{\lambda}_1 \cdots C^{\lambda}_r$. Thus, by
\cite[Remark~1]{JM05}, we see that \cite[Theorem~2.2]{JM05} is
equivalent to the assertion that $\M(\Sigma_{0,r},C^{\lambda})$ is
non-empty if and only if for any non-negative $d \leq (r-1)/2$ and any
subset $W \subset \{1,\ldots,r\}$ of cardinality $(r-1)-2d$ we have that
\begin{equation} \label{eq:JeffreyInequality}
S_W:=\sum_{j \in \{1,\ldots,r\} \setminus W} \lambda_j-\sum_{j \in W}     \lambda_j \leq 2d \pi.
\end{equation}
We will finish the proof by showing that this condition is equivalent
to~\eqref{eq:JeffreyTypeInequality}.
%
Given~$\ul$, we define $\lambda^{(\ul)}\in [0,\pi]^r$ by
$\lambda^{(\ul)}_j=\pi \ell_j/p_j$. Then
$C^{(\ul)}=C^{\lambda^{(\ul)}}$. Let $W \subset \{1,\ldots,r\}$ be a
subset of cardinality $r-1-2d$ for some non-negative integer
$d\leq (r-1)/2$. Multiplying both sides
of~\eqref{eq:JeffreyInequality} by $-\pi^{-1}$, we see
that~\eqref{eq:JeffreyInequality} is equivalent to
\begin{equation} \label{eq:JeffreyInequality2}
    -2d \leq \sum_{j \in W} \ell_j/p_j- \sum_{j \in \{1,\ldots,r\} \setminus W} \ell_j/p_j.
\end{equation}
Let $W^c$ denote the complement of $W\subset \{1,\ldots,r\}$. We can rewrite the right hand side as follows
\begin{align}
    \sum_{j \in W} \ell_j/p_j- \sum_{j \in \{1,\ldots,r\} \setminus W} \ell_j/p_j&=\sum_{j \in W} \ell_j/p_j+ \sum_{j \in \{1,\ldots,r\} \setminus W} (p_j-\ell_j-p_j)/p_j
    \\ \label{eq:secondline} & =\sum_{j \in W} \ell_j/p_j+ \sum_{j \in \{1,\ldots,r\} \setminus W} (p_j-\ell_j)/p_j- \lvert W^c \rvert.
\end{align}
Thus, by adding $\lvert W^c \rvert=r-\lvert W \rvert=1+2d$ to both
sides of~\eqref{eq:JeffreyInequality2}, we see
from~\eqref{eq:secondline} that~\eqref{eq:JeffreyInequality} is
equivalent to
\begin{equation} 
    1\leq \sum_{j \in W} \ell_j/p_j+ \sum_{j \in \{1,\ldots,r\} \setminus W} (p_j-\ell_j)/p_j.
\end{equation}
By coprimality considerations, we see that this is equivalent to
\eqref{eq:JeffreyTypeInequality} with $J=W^c$
(note that every subset~$J$ of odd cardinality is of that form).
%
This finishes the proof.
\end{proof}


\begin{corollary} \label{Cor:CS}
  For each $\ul=(\ell_1,\ldots,\ell_r) \in
      \fR(p_1,\ldots,p_r)$, the Chern-Simons action
      functional~$\cS$ of~\eqref{eqdefaction} is constant on the
      component of~$\M(X)$ isomorphic to $\M^{\Irr}(\Sigma_{0,r}, C^{(\ul)})$;
its value there is
\begin{equation}   \label{eqdefSl}
S_\ul := -\frac{1}{4P} \bigg ( \sum\limits_{j=1}^r  \ell_j\wh
        p_j\bigg)^{\!2} 
\end{equation}
\!\!$\mod\Z$ (with the notation~$\wh p_j$ of~\eqref{eqdefwhpj}).
Consequently,
  \beglab{eqformulCSX}
    \CS(X)=\{0\} \sqcup \{ S_\ul \!\!\!\mod\Z
    %
    %
    \mid \ul \in 
      \fR(p_1,\ldots,p_r) \} \subset \Q / \Z.
  \edla
\end{corollary}

\begin{proof}
  This follows directly from Theorem~\ref{thmMflat} 
  %
  %
  together with \cite[Proposition~8]{AM22} (which of course builds on \cite{KK91,FS90}).
\end{proof}

\begin{remark}
This is to be compared with Theorem~1 of \cite{AM22} for the
$\SL(2,\C)$ Chern-Simons actions, which the results of
Appendix~\ref{secHikamisets} allow to rephrase as
\[
{\CS}_{\C}(X) = \{0\} \sqcup \{ S_\ul \!\!\!\mod\Z
    \mid \ul \in  \fL(p_1,\ldots,p_r) \}
  \]
with the natural extension of the explicit
  definition~\eqref{eqdefSl} of~$S_\ul$ to the case of
  $\ul \in \fL(p_1,\ldots,p_r)$. 
\end{remark}

\begin{example}
%
%
%
  The triple $(p_1-1,p_2-1,\ldots,p_r-1)$ always belongs to
  $\fL(p_1,\ldots,p_r)$. When $r=3$, it belongs to $\fR(p_1,p_2,p_3)$ if and only if
  $(p_1,p_2,p_3)=(2,3,5)$.
  In fact, $\fL(2,3,5)=\fR(2,3,5)$ consists of this triple, $(1,2,4)$,
  and one more triple: $(1,2,2)$, the corresponding Chern-Simons
  actions being $S_{(1,2,4)}=-1/120$ and $S_{(1,2,2)}=-49/120$ $\!\!\mod \Z$.
  In the case $(p_1,p_2,p_3)=(2,3,7)$, we find
$\fR = \{(1,2,2), (1,2,4)\} \subsetneq \fL = \{(1,2,2), (1,2,4),
(1,2,6)\}$, and the corresponding $\SU(2)$ Chern-Simons actions are
$-25/168$ and $-121/168$ $\!\!\mod\Z$, while ${\CS}_{\C}(X)$ has one
more element, $S_{(1,2,6)}=-1/168$ $\!\!\mod\Z$.
An example with $r=4$ is $(p_1,\ldots,p_4)=(2,3,5,7)$, for which
$(1,2,4,6)\in\fR$ and the cardinalities are $\#\fR=22$ and $\#\fL=29$.
\end{example}


\subsection{Proof of Theorem~\ref{Thm:main}}  \label{secpfthmmain}



The case $\al=0$ of Corollary~\ref{cor:wtPsial} says that the modified
GPPV invariant can be written as a median sum of the
resurgent-summable formal series~$\wt\Psi_0(\tau)$ defined
by~\eqref{def:Zalpha} or~\eqref{eqwtPsizaltern},
\beglab{eqPsiasamedsumofwtPisz}
\Psi(\tau) = \bS^{\frac\pi2}\median \wt\Psi_0(\tau)
\quad \text{for} \ens \tau\in\HH,
\edla
median sum meaning the half-sum of lateral Borel-Laplace sums in our
case (cf.\ footnote~\ref{ftnLaplMedSum}).
For any $k\in\Z_{\ge1}$, the case $\al=1/k$ of
Corollary~\ref{cor:wtPsial} entails that the non-tangential limit
\beglab{eqlimPsioneovk}
\lim_{\tau\to1/k}\Psi(\tau) = \Psi_{1/k,0}
\edla
exists.
Theorem~\ref{Thm:main} is about $\WRT_k(X)$, but
Proposition~\ref{propradiallimitWz} shows that it is sufficient to
study the numbers~\eqref{eqlimPsioneovk}
(compare~\eqref{eqdefwtWz}--\eqref{eqWRTklimbSmedWz} with
  \eqref{eqPsiasamedsumofwtPisz}--\eqref{eqlimPsioneovk}).

We will compare~$\Psi(\tau)$ written as the median sum of $\wt\Psi_0(\tau)$ and one of its two lateral
Borel-Laplace sums, namely
\begla
\bS^{\frac\pi2-\eps}\wt\Psi_0(\tau) = \bS^0\wt\Psi_0(\tau).
\edla
Clearly, $\Psi(\tau)- \bS^0\wt\Psi_0(\tau)$ is half the difference of
the two lateral Borel-Laplace sums, which is a particular case of
``Stokes phenomenon'',\footnote{The terminology ``Stokes phenomenon''
  classically pertains to the theory of linear meromorphic systems of ODEs, but in
  the context of resurgence it is often used to refer to the
  difference of two Borel-Laplace sums computed by means of resurgent
  analysis in the Borel plane.}
and Proposition~\ref{propUnivPol} will give us the tools to compute
it.

Note that the function $\bS^0\wt\Psi_0$ analytically extends to much
more than the upper half-plane $\HH=\{0<\arg\tau<\pi\}$: the Borel
summability statement in Corollary~\ref{cor:wtPsial} allows us to 
follow its analytic continuation up
to $\{-2\pi<\arg\tau<\pi\}$;
in particular it is analytic on $\R_{>0}=\{\arg\tau=0\}$.
By way of contrast, the difference of the two lateral
  Borel-Laplace sums of~$\wt\Psi_0(\tau)$ is a priori defined in~$\HH$
  only, but we will see that it can be expressed as a sum of partial
  theta series evaluated at~$-\tau\ii$ that have non-tangential limits
  at any rational number.
Letting~$\tau$ tend to the positive rational number~$1/k$, the upshot will be

\begin{proposition}   \label{propequivmainthm}
  For each $k\in\Z_{\ge2}$, we have
  \beglab{eqlimPsioneovksum}
\lim_{\tau\to k\ii}\Psi(\tau) = (\bS^0 \wt\Psi_0)(k\ii) +
    \sum_{\ul \in \fR(p_1,\ldots,p_r)} e^{2\pi ikS_\ul}\,
   k^{1/2} H^\ul(k)
  \edla
with $\fR(p_1,\ldots,p_r)\subset\Z^r$ and~$S_\ul$ as in~\eqref{eq:JeffreyTypeInequality} and~\eqref{eqdefSl},
and where $H^\ul$ is a polynomial, defined by~\eqref{eqdefHul} below, satisfying
\beglab{ineqdegHul}
\deg H^\ul \le r-3-t_\ul
\quad\text{with}\;\, t_\ul \;\, \text{as in~\eqref{eqdeftul}.}
\edla
  \end{proposition}

  \begin{proof}[Proof that Proposition~\ref{propequivmainthm} implies
    Theorem~\ref{Thm:main}]
    Let $k\in\Z_{\ge2}$.
    As in the proof of Proposition~\ref{propradiallimitWz},
    since~$\cE(\tau)$ is convergent in the unit disc
    and~$\wt\Psi_0(\tau)$ is divisible by~$\tau$,
    formula~\eqref{eqWRTklimbSmedWz} implies that
    \begla
    \WRT_k(X) = k\cE(k\ii)\lim_{\tau\to k\ii}\bS^{\pi/2}\median\wt\Psi_0(\tau), \qquad
    (\bS^0\wt W_0)(k\ii) = k\cE(k\ii) (\bS^0 \wt\Psi_0)(k\ii).
    \edla
    Taking~\eqref{eqlimPsioneovksum} for granted and using~\eqref{eqdefSl}--\eqref{eqformulCSX}, we get
    \begin{align}
    \WRT_k(X) &= (\bS^0 \wt W_0)(k\ii) +
    \sum_{\ul \in \fR(p_1,\ldots,p_r)} e^{2\pi ikS_\ul}\,
                k^{3/2} \cE(k\ii) H^\ul(k) \notag \\[.8ex]
      &= (\bS^0 \wt W_0)(k\ii) +
      \sum_{S \in \CS(X) \setminus \{0\}} e^{2\pi ikS}\,
   k^{3/2} \cE(k\ii) H_S(k)   \label{eqWRTkwtWz}
  \\[-2.5ex]
      \intertext{with}
    & \quad H_S(k) := \sum_{\ul \in \fR(p_1,\ldots,p_r) \;\text{s.t.}\; S_\ul = S}
    H^\ul(k).   \label{eqdefHSsumHul}
    \end{align}
    Since the functions~$H_S$ are polynomials in~$k$, 
    formula~\eqref{eqWRTkwtWz} gives rise to an asymptotic expansion of $\WRT_k(X)$ for
    $k\to\infty$, with finitely many different exponentials
    modulated by Laurent formal series in~$k^{-1/2}$, and these formal
    Laurent series are uniquely determined.
In particular, the formal series~$\wt W_0(k\ii)$ is uniquely determined and
must coincide with the series~$W_0(k\ii)$ already found by Lawrence
and Rozansky in~\cite{LR99} in terms of the Ohtsuki series (cf.~\eqref{eqdefWzerolaX}).
This proves that $\wt W_0=W_0$ (alternatively, this identity can be inferred
from \cite{AM22}).

Now, for each $\ul\in\fR(p_1,\ldots,p_r)$, the upper bound~\eqref{ineqdegHul} for the degree of the
polynomial $H^\ul(k)$ is nothing but
$\frac12 \dim\M^{
  \Irr}(\Sigma_{0,r}, C^{(\ul)})$,
by~\eqref{eqdimformula}.
Therefore, for $S \in \CS(X) \setminus \{0\}$, formula~\eqref{eqdefHSsumHul} shows
that $\deg H_S$ is at most half the dimension~$d_S$ referred to in~\eqref{eqWSkZskdS}.
\end{proof}


\begin{proof}[Proof of Proposition~\ref{propequivmainthm}]
  Putting together~\eqref{eqPsicQsumThchij}
  and~\eqref{eqchijmsfuo}, we have
  \beglab{eqPsicQsumjsCjsTh}
    \Psi(\tau) = \cQ(\tau) + \sum_{\substack{j,s\ge0, \; j+s \le r-3 \\[.3ex] j+s\equiv r-1 \ [2]}}
    C_{j,s} \, \Theta(\tau; j,m^s f^\uo \circ\cT,2P)
    \quad \text{for}\ens \tau\in\HH,
    \edla
with~$C_{j,s}$ as in~\eqref{eqdefCjs} and $m^s f^\uo$ defined
by~\eqref{eqsupportsSuo}--\eqref{eqdefmsfuo}.
Recall that $\cT=\ID_{\Z}$ or $\ID_{\Z}-P$ according
as $r$ is odd or even.

Our strategy is to make use of formula~\eqref{equationShigherdepth} in
Proposition~\ref{propUnivPol}.
This is possible because, for each $(j,s)$ involved
in~\eqref{eqPsicQsumjsCjsTh}, the $2P$-periodic function
$m^s f^\uo\circ\cT$ fulfills the hypotheses of
Proposition~\ref{propUnivPol}:
its parity is $r+s \equiv j+1 \ [2]$,
and its mean value is~$0$,
as proved in Corollary~\ref{corcVs} in the appendix.
%
%
Applying~\eqref{equationShigherdepth} and recalling the
definition~\eqref{def:Zalpha} of $\wt\Psi_0$, we find
\begin{align}
\Psi(\tau) &= \bS^{\frac\pi2- \epsilon}\wt\Psi_0(\tau)
\notag \\[.8ex]
& \qquad - \sum_{\substack{j,s\ge0, \; j+s \le r-3 \\[.3ex] j+s\equiv r-1 \ [2]}}
      2^{-[\frac{j}{2}]} \, i^{\frac12} C_{j,s}
      \sum\limits_{\substack{0 \leq \nu \leq j \\[.3ex] \nu \equiv j \ [2]}}
      \Big( \frac{2P}{\pi i}\Big)^{\!\frac {j-\nu}2} P_{j,\nu} \,
      \tau^{-\frac{j+\nu+1}{2}}
  \Theta(-\tau^{-1};\nu,\rwh{m^s f^\uo \circ\cT},2P)
  \\[1.3ex]
           &= \bS^0\wt\Psi_0(\tau) +
             \sum_{\substack{\nu,s\ge0, \; \nu+s \le r-3 \\[.3ex]
  \nu+s\equiv r-1 \ [2]}}
  \tau^{-1/2} Q_{\nu,s}(\tau\ii) \Theta(-\tau^{-1};\nu,\rwh{m^s f^\uo \circ\cT},2P)
\label{eqtherhs}
\end{align}
with polynomials
\beglab{eqdefQnusx}
    Q_{\nu,s}(x) := - \sum_{\substack{\nu\le j\le r-3-s \\[.3ex] j\equiv \nu [2]}}
    2^{-[\frac{j}2]} \, i^{\frac12}
    \Big(\frac{2P}{\pi i}\Big)^{\frac{j-\nu}{2}}
    C_{j,s} P_{j,\nu} \, x^{\frac{j+\nu}2}
    \in \C[x].
    \edla
    We thus need to inquire about the non-tangential limit as
    $\tau\to1/k$ of the \rhs\ of~\eqref{eqtherhs}, which amounts to
    asking whether, for each pair $(\nu,s)$ in the finite sum,
    \beglab{eqlimmkThrwhmsfuo}
    \lim_{\tau\to-k}\Theta(\tau;\nu,\rwh{m^s f^\uo \circ\cT},2P)
    \edla    
    exists and what it is.
    To that end, we need information about the DFT $\rwh{m^s f^\uo
      \circ\cT}$, which is provided by
    
    \begin{lemma}   \label{lemDFTmsfuo}
      For each $0\le s \le r-3$, the support of $\rwh{m^s f^\uo
      \circ\cT}$ is contained in the union over
    $\{\ul\in\fL(p_1,\ldots,p_r)\mid t_\ul\le s\}$ of the sets
    \beglab{eqdefSul}
    \fS^\ul := \cN^\ul(E) + 2P\Z \, \subset\, \Z,
    \edla
    where the function $\cN^\ul \col E = \{+1,-1\}^r \to \Z$ is defined by
    \beglab{eqdefcNul}
    \cN^\ul(\ue) :=
    P + \sum_{j=1}^r \eps_j \ell_j \whp_j.
    \edla
    \end{lemma}

Lemma~\ref{lemDFTmsfuo} is a direct consequence of
    Proposition~\ref{propositionDFTcomputation} in
    Appendix~\ref{apppflemDFTmsfuo}.
    Note that~$\cN^\ul$ and~$\fS^\ul$ are generalizations of~$\cN$
    and~$\fS^\uo$ defined in~\eqref{eqdefCNuo} and~\eqref{eqsupportsSuo}.
    
    Another information that we need is provided by
    Lemma~\ref{lemmaHikamisets2}(iii)--(iv)
    to be found in Appendix~\ref{secHikamisets}: it shows that the
    sets $\fS^\ul$ are pairwise disjoint and, on each of them, the
    function $m\mapsto\fourp{m^2}$ is constant.
    Equivalently, the function $m\in\fS^\ul\mapsto\oone{-\frac{m^2}{4P}}\in\Q/\Z$ is constant;
    evaluating at $\ue=(-1,1,\ldots,1)$ and comparing with~\eqref{eqdefSl}, we find
    \beglab{eqfSulSsigoul}
    m\in\fS^\ul \imp
    \oone{-\frac{m^2}{4P}} = \oone{S_{\sig_1(\ul)}}
    \quad \text{with} \quad
    \sig_1(\ell_1,\cdots,\ell_r):=(p_1-\ell_1,\ell_2,\cdots,\ell_r)
    \edla
    (note that~$\sig_1$ is an involution of $\fL(p_1,\ldots,p_r)$ that
    leaves~$t_\ul$ invariant).
Finally, another useful consequence of
Proposition~\ref{propositionDFTcomputation} is that, for each $\ul\in\fL(p_1,\ldots,p_r)$,
\beglab{eqprodindic}
\text{the product function}\ens
\rwh{m^s f^\uo\circ\cT} \cdot \indul
\ens \text{has mean value~$0$}
\edla
(note that the indicator function~$\indul$ is
$2P$-periodic too).

    Taking these few facts for granted, to study~\eqref{eqlimmkThrwhmsfuo},
    we can write
    \begla
    \rwh{m^s f^\uo\circ\cT}  = \sum_{\ul\in\fL(p_1,\ldots,p_r)\,\text{s.t.}\, t_\ul\le s}
    \rwh{m^s f^\uo\circ\cT} \cdot \indul
    \edla
    (note that, for each $\ul\in\fL$, considering the product function
    of~\eqref{eqprodindic} amounts to considering the restriction of
    $\rwh{m^s f^\uo\circ\cT}$ to~$\fS^\ul$).
    We now apply to the corresponding sum of partial theta
      series two elementary observations 
    (obvious consequence of the definition~\eqref{eqdefThetafjM} for
    the first one, and the $\al=0$ case of~\eqref{eqexistlimal} for
    the second one):
    \begin{lemma}
      Let $\nu\in\Z_{\ge0}$.
      %

      \noindent \emph{\textbf{(i)}}
      If~$f$ is an $M$-periodic function on~$\Z$ and there exists
      $\tht\in\Q/\Z$ such that, for any $m\in\Z$,
      \begla
      f(m)\neq0 \iimp \oone{-\frac{m^2}{4P}} = \tht,
      \edla
      then
      \begla
      \Theta(\tau-k;\nu,f,M) = e^{2\pi ik\tht} \Theta(\tau;\nu,f,M)
      \quad
      \text{for all}\ens \tau\in\HH
      \ens \text{and}\ens k\in\Z.
      \edla
      %

      \noindent \emph{\textbf{(ii)}}
      If moreover~$f$ has zero mean value, then
      $\Z\subset \sQ_{f,M}$
      and
      \begla
      \lim_{\tau\to-k} \Theta(\tau;\nu,f,M) = e^{2\pi ik\tht} \lim_{\tau\to0} \Theta(\tau;\nu,f,M)
      \quad \text{for all}\ens k\in\Z.
      \edla
      \end{lemma}

      We thus obtain
          \begin{multline}
          \lim_{\tau\to-k}\Theta(\tau;\nu,\rwh{m^s f^\uo \circ\cT},2P) =
          \sum_{\ul\in\fL(p_1,\ldots,p_r)\,\text{s.t.}\, t_\ul\le s}
          e^{2\pi ikS_{\sig_1(\ul)}} \La_{\nu,s,\ul} \\[.8ex]
          \text{with}\quad \La_{\nu,s,\ul} := \lim_{\tau\to0}
          \Theta(\tau;\nu,\rwh{m^s f^\uo \circ\cT}\cdot \indul,2P).
    \end{multline}
    Plugging that into~\eqref{eqtherhs}, we get
    \begin{align}
      \notag
    \lim_{\tau\to k\ii}\Psi(\tau) &= (\bS^0 \wt\Psi_0)(k\ii) +
    \sum_{\substack{\nu,s\ge0, \; \nu+s \le r-3 \\[.3ex]
  \nu+s\equiv r-1 \ [2]}}
  k^{1/2} Q_{\nu,s}(k) \sum_{\ul\in\fL(p_1,\ldots,p_r)\,\text{s.t.}\, t_\ul\le s}
      e^{2\pi ikS_{\sig_1(\ul)}} \La_{\nu,s,\ul} \\[.8ex]
      \notag
      &= (\bS^0 \wt\Psi_0)(k\ii) +
      \sum_{\ul\in\fL(p_1,\ldots,p_r)} e^{2\pi ikS_{\sig_1(\ul)}}\, k^{1/2}
        H^{\sig_1(\ul)}(k)
      \\[-2.5ex]
      \intertext{with}
      & \qquad H^\ul(k) := \sum_{\substack{\nu\ge0,\; s\ge t_\ul, \; \nu+s \le r-3 \\[.3ex]
  \nu+s\equiv r-1 \ [2]}}
  \La_{\nu,s,{\sig_1(\ul)}} \, Q_{\nu,s}(k) \in \C[k].
      \label{eqdefHul}
      \end{align}
      Note that, by~\eqref{eqdefQnusx}, $\deg Q_{\nu,s} \le
      \frac12(\nu+r-3-s)$.
      We thus have $\deg Q_{\nu,s} \le r-3-s$ for each term
      in~\eqref{eqdefHul}, whence~\eqref{ineqdegHul} follows.

      There only remains to be proved that
      \begla
      \ul \notin \fR(p_1,\ldots,p_r) \imp
      H^\ul=0.
      \edla
      This follows from Corollary~\ref{corvanish},
      which says that, if $\ul \notin \fR(p_1,\ldots,p_r)$ while
      $\nu+s\le r-3$ and $\nu+s\equiv r-1\ [2]$, then
      $\lim\limits_{\tau\to0}
      \Theta(\tau;\nu,\rwh{m^sf^{\uo}\circ\cT}\cdot \indusl,2P)
      = 0$, i.e.\
      $\La_{\nu,s,\sig_1(\ul)}=0$.
\end{proof}

\appendix

\section{Appendix on normalizations} \label{Appendix:Normalizations}

Let $L \subset S^3$ be a framed oriented link and let
$M=S^3_L$. Recall the definition of $\WRTk(M)$ as given
in~\eqref{def:WRT} (Definition~\ref{ddef:WRT}).
The topological invariant $\mathcal{F}_k(M)\in\C$
introduced in \cite{RT91} is given by
$\mathcal{F}_k(M)=(C_0)^{b_1(M)}\WRTk(M)$, where $C_0 \in \C^\times$
is a $k$-dependent constant discussed in \cite[Appendix~A]{Hansen01}. In \cite{GPPV} they use the following notations
\[\tau_k(M)= (G_{k,0})^{b_1(M)}\WRTk(M), \qquad Z_{\SU(2)_k}(M):=\frac{\tau_k(M)}{\tau_k(S^1 \times S^2)}.
\]
The authors of \cite{GPPV} refer to the $(S^1 \times S^2)$-normalized invariant $Z_{\SU(2)_k}(M)$ as the Witten-Reshetikhin-Turaev invariant, or the quantum Chern-Simons partition function. We note that for a rational homology sphere $M$, we have $b_1(M)=0$ and therefore $\tau_{k}(M)=\WRTk(M)=\tau_{\zeta_k}(M),$ where $\tau_{\zeta_k}(M)$ is the invariant considered in \cite{Habiro08}. In particular, for a rational homology sphere $M$, we have $Z_{\SU(2)_k}(M)= \WRTk(M)/\WRTk(S^1\times S^2)$.

\subsection{Rational surgeries}

In \cite{Hansen01} the rational surgery formula from \cite{Jeffrey92}
is generalized to Reshetikhin-Turaev invariants defined for more
general modular tensor categories. The main result is \cite[Theorem~5.3]{Hansen01}. Consider the modular tensor category $\mathcal{V}_k$
(denoted by $\mathcal{V}_t$ in \cite{Hansen01}) associated with the
quantum group $U_q(\mathfrak{sl}(2,\C))$, where $q=\zeta_k$. Let
$\mathcal{D}_k=G_{k,0}=\WRT(S^1 \times S^2)$. This is a so-called rank
of $\mathcal{V}_k$, and it satisfies
$\mathcal{D}^2= \sum_{j=1}^{k-1} [j]^2$. The invariant
$\tau_{\mathcal{V}_k,\mathcal{D}_k}(M) \in \C$ considered in
\cite{Hansen01} is given by
\begin{equation}\label{eq:MTCnormalization}
    \tau_{\mathcal{V}_k,\mathcal{D}_k}(M)= \frac{\WRTk(M)}{\WRTk(S^1 \times S^2)}.
\end{equation}
This identity follows from the material presented in \cite[Appendix~A]{Hansen01}. The invariant \eqref{eq:MTCnormalization} extends to triples $(M,L',\lambda')$, where $L'\subset M$ is a framed oriented link and $\lambda' \in \Lambda_k^{\pi_0(L')}$ is a coloring. For $M=S^3$, we have that $\tau_{\mathcal{V}_k,\mathcal{D}_k}(S^3,L',\lambda')=\mathcal{D}_k^{-1}J_{\lambda'}(L',\zeta_k)$, where, as above, $J_{\lambda'}(L',\zeta_k)$ is the colored Jones polynomial. Given rational surgery data $(L,B)$ \cite[Corollary~8.3]{Hansen01} gives 
\begin{equation} \label{eq:RationalSurgery2}
\tau_{\mathcal{V}_k,\mathcal{D}_k}(S^3_{L,B})=\frac{\exp\left({\frac{\pi i}{4}\left(\frac{k-2}{k}\right) \Phi(L,B)}\right)}{\WRTk(S^1 \times S^2)} \sum_{\lambda \in \Lambda_k^{\pi_0(L)}} J_{\lambda}(L,\zeta_k) \prod_{j \in \pi_0(L)} \rho_k(B_j)_{\lambda_j,1},
\end{equation}
where we used the notation from Section \ref{sec:RationalSurgery} (and
substituted the identity
$\tau_{\mathcal{V}_k,\mathcal{D}_k}(S^3,L',\lambda')=\mathcal{D}_k^{-1}J_{\lambda'}(L',\zeta_k)$
into the right hand side of the central equation in
\cite[Corollary~8.3]{Hansen01}). Note that
\eqref{eq:rationalsurgeryformula} is consistent with
\eqref{eq:MTCnormalization} and \eqref{eq:RationalSurgery2}.

\subsection{The normalization used in the work of Lawrence and
  Rozansky} Consider again the Seifert fibered integral homology
sphere $X$. As described above: in \cite{LR99} the invariant
$\WRTk(X)$ is computed by implementing the rational surgery formula
\eqref{eq:rationalsurgeryformula} to a specific surgery
presentation. They work with a $S^3$-normalized invariant which they
denote by $Z_k(X)$, and they state a rational surgery factor for $Z_k$
in \cite[eqn~(3.2)]{LR99}. The surgery formula is equal to the one
given in this article in~\eqref{eq:rationalsurgeryformula} times
$G_{k,0}$. However, we observe that in their computation of $Z_k(X)$
in \cite[Sec.~4]{LR99}, they actually omit this factor $G_{k,0}$ and
implement the formula given in~\eqref{eq:rationalsurgeryformula}.

\begin{remark}
  As a heed of caution we remark that the normalization coefficients
  $G_{k,\pm}$ introduced in this article in Section~\ref{sec:defWRT}
  are standard in the literature, but they differ from the
  normalization coefficients denoted by $G_\pm$ in \cite{LR99} and
  used in their surgery formula for WRT invariants
  \cite[eqn~(3.1)]{LR99}. However, the coefficients $G_{\pm}$ are not
  used directly in the computation of $Z_k(X)$ in \cite[Sec.~4]{LR99},
  where they use instead the rational surgery formula
  \eqref{eq:rationalsurgeryformula}, which does not involve $G_{\pm}$
  directly, but agree (up to an overall factor of $G_{k,0}$ as
  explained above) with the standard formula for WRT invariants in
  terms of rational surgeries, as can be found in
  \cite{Jeffrey92,Hansen01}. Therefore, in spite the fact that there
  seems to be a minor inconsistency between \cite[eqn~(3.1)]{LR99} and
  \cite[eqn~(3.2)]{LR99}, the results from \cite[Sec.~4]{LR99}
  applies to the normalized invariant which we denote by $\WRTk(M)$ in
  this article.
\end{remark}



\section{Appendix on Hikami sets, $s$-Hikami functions and their Discrete
  Fourier transforms} \label{sectiontoolbox}

We recall that $r\ge3$ and $p_1,\ldots,p_r$ are positive and pairwise
coprime, with~$p_j$ odd for $j\ge2$.
%
%
%
Recall also the notation~\eqref{eqdefnotabrN} for the canonical projection
$\brN{\ \cdot\ } \col \Z \to \Z/N\Z$.
From now on, we will simply denote by~$\fL$ and~$\fR$ the subsets
of~$\Z^r$ introduced in~\eqref{eqmathfrakL}
and~\eqref{eq:JeffreyTypeInequality}.
Since we will need to deal with various subsets of $\{1,\ldots,r\}$
and their complements, we use the notation
\[
J \subset \{1,\ldots,r\} \Imp \compl J := \{1,\ldots,r\} \setminus J.
\]
\bblack

\subsection{Hikami sets}   \label{secHikamisets}

We define a subset $\fH=\fH(p_1,\ldots,p_r)$ of~$\Z^r$ by
  \begin{multline}   \label{equationHikamisets}
  \qquad \fH := \Big\{ \uh=(h_1,\ldots,h_r)\ \Big| \ 
      0\leq h_{j}\leq p_{j}\ \text{for all~$j$,} \\
        \qquad\qquad \frac{h_j}{p_j}\notin\Z  \ \text{for at least three
          values of~$j$}  \Big\}. \qquad
\end{multline}
Note that
%
  %
$ \fR \subset \fL \subset \fH$.
  %
%
For any $\uh \in \fH$, we define
\beglab{eqdefJuh}
J^\uh := \big\{ j\in\{1,\ldots,r\} \;\text{such that}\; h_j \equiv0\ [p_j] \big\}.
\edla
%
%
Thus, with reference to~\eqref{eqdeftul}, %
\begla
0\le t_\uh=|J^\uh|\le r-3.
\edla
We also generalize the sets~$\fS^\ul$ and the functions~$\cN^\ul$
defined in~\eqref{eqdefSul}--\eqref{eqdefcNul} for $\ul\in\fL$ to the
case of an arbitrary $\uh\in\fH$:
\begin{equation}
\fS^\uh := \cN^\uh(E)+2P\Z \, \subset\, \Z, \qquad
    \cN^{\uh} \col \ue \in E=\{ +1,-1\}^r \mapsto P+\sum\limits_{j=1}^r \eps_j h_j \whp_j \in \Z
  \end{equation}
 (recall that $P=p_1\cdots p_r$ and $\wh p_j=P/p_j$).
  %
%
%
 Finally, we define
 \beglab{eqdeftiJn}
 \ti J_n := \big\{ j\in\{1,\ldots,r\} \;\text{such that}\; n \equiv0\ [p_j] \big\}
 \quad \text{for any}\ens n\in\Z.
 \edla
%
 %
Writing $\cN^\uh(\ue) = P + \sum\limits_{j\in J^\uh} \eps_j h_j \wh
p_j + \sum\limits_{j\in \ccompl J^\uh} \eps_j h_j \wh p_j$,
one easily checks

\begin{lemma}\label{lemmaHikamisets1}
  Let $\uh \in \mathfrak{H}$. 
  \smallskip

      \noindent \emph{\textbf{(i)}}
      For all $\ue\in E$, $\ti J_{\cN^\uh(\ue)} = J^\uh$.
      \medskip
      
      \noindent \emph{\textbf{(ii)}}
      Consider the map $\tp{\cN^\uh} \col E \to \Z/2P\Z$. 
      Each element of its range is of the form $\tp n$
      with $n\in\fS^\uh$, and it has exactly $2^{t_\uh}$ preimages~$\ue$, all
      of which have the same restriction to~$\compl J^\uh$:
      the restriction $\ue_{\mid\ccompl J^\uh}$ is determined by
      $\tp{n}$ but the restriction $\ue_{\mid J^\uh}$ is free.
%
      \medskip
      
      
      \noindent \emph{\textbf{(iii)}}
 The map $\ue \in E \mapsto \fourp{\big(\cN^{\uh}(\ue)\big)^2}$ is
 constant.
  \end{lemma}


We now define a equivalence relation in~$\fH$ by declaring that
$\uh\sim\uh'$ if
\begin{equation} \label{equhsimuhp}
\exists J \subset \{1,\ldots,r\} \;\, \text{with}\; |J|\; \text{even such that}\;
    \begin{cases}
         h_j = p_j- h_j^\prime  & \text{for}\ens j\in J
        \\[.8ex]
         h_j =  h_j^\prime & \text{for}\ens j\in \compl J.
      \end{cases}
\end{equation}
%
%
The reader may check

\begin{lemma}\label{lemmaHikamisets2} The following properties hold:
      \medskip
      
      \noindent \emph{\textbf{(i)}}
The set~$\fL$ of~\eqref{eqmathfrakL} is a system of representatives
of~$\fH/\sim$.
      \medskip
      
      \noindent \emph{\textbf{(ii)}}
      Given $\uh,\uh' \in \fH$, \quad
      $
        \uh \sim \uh' \imp J^\uh=J^{\uh'} \ens\text{and}\ens \fS^{\uh}= \fS^{\uh'}.
        $

         \medskip
      
      \noindent \emph{\textbf{(iii)}}
      Given $\ul,\ul' \in \fL$, \quad
      $
        \ul\neq\ul' \imp \fS^\ul \cap \fS^{\ul'} = \emptyset.
      $
      
      \medskip
      
      \noindent \emph{\textbf{(iv)}}
      For each $\ul\in\fL$, \; the map $m\in\fS^\ul \mapsto \fourp{m^2}
      \in \Z/4P\Z$ is constant.
    \end{lemma}


We call ``Hikami sets'' the sets $\fS^\uh$, $\uh\in\fH$ (or, without
loss of generality, $\uh\in\fL$).
Note that, in view of Lemma~\ref{lemmaHikamisets1}(i), no multiple
of~$P$ belongs to any of these sets:
\begla
\fS^\uh \, \cap\,  P\Z = \emptyset
  \quad \text{for all}\ens \uh\in\fH.
\edla
In particular, we cannot have $\cN^\uh(\uo)=2P$, whence
  \beglab{eqneverone}
  \sum\limits_{j=1}^{r} \frac{h_j}{p_j}\neq1
  \quad \text{for all}\ens \uh\in\fH.
  \edla

\begin{remark}
The set~$\fR$ of~\eqref{eq:JeffreyTypeInequality} can be written as 
\begin{equation}   \label{equationRp}
  \fR = \Big\{\ul\in \fL \Big| \ 
\sum\limits_{j=1}^{r} \frac{h_j}{p_j}>1
  \; \text{ for all } \uh\sim \sig_1(\ul) \text{ in } \fH \Big\},
\end{equation}
  where~$\sig_1$ is the involution 
$(\ell_1,\ldots,\ell_r)\mapsto(p_1-\ell_1,\ell_2,\ldots,\ell_r)$.
\end{remark}

\begin{remark}   \label{remCSHikamisets}
  Putting together Lemma~\ref{lemmaHikamisets2} and
  \cite[Proposition~6]{AM22} recalled in Section~\ref{secrefAMprop6},
we obtain a one-to-one correspondence between the components of the
  moduli space of irreducible flat $\SL(2,\mathbb{C})$-connections
  (labelled by $\fL$) and Hikami sets.
The Chern-Simons action associated with the component labelled by
$\ul\in\fL$ has been computed in~\eqref{eqdefSl} and~\eqref{eqfSulSsigoul}:
it is $\oone{S_{\sig_1(\ul)}}$.
\end{remark}

We will be interested in subsets of~$\Z$ obtained as disjoint
unions of certain Hikami sets.
%

\begin{lemma}\label{lemmaresiduesets}
  Let $s\in\{0,\ldots,r-3\}$ and
  \begla
  \Mss := \{\, n\in \Z \ \text{such that}\ \, |\ti J_n| > s \,\}
  \edla
  (with reference to~\eqref{eqdeftiJn} for the notation~$\ti J_n$).
  Then
  \beglab{eqdecZMssfL}
  \Z = \Mss \, \sqcup \bigsqcup_{\ul\in\fL \,\text{s.t.}\, t_\ul\le s}\fS^\ul.
  \edla
\end{lemma}

\begin{proof}
  We will prove the following more precise statement:
  for every subset $J\subset \{1,\ldots,r\}$ such that $|J|\le r-3$,
  \begin{equation}\label{equationdecompositioninlemmamultiplesetisHikamiset}
\{\, n\in\Z\mid \ti J_n = J \,\}
    = \bigsqcup_{ \ul \in \fL \, \text{s.t.} \,  J^\ul = J } \fS^\ul,
      %
      %
 \end{equation}
with reference to~\eqref{eqdefJuh} for the notation~$J^\ul$.
The decomposition~\eqref{eqdecZMssfL} will then follow
from~\eqref{equationdecompositioninlemmamultiplesetisHikamiset}
by writing~$\Z$ as the disjoint union of 
$\{\, n\in \Z \ \text{such that}\ \, \ti J_n = J \,\}$
over all subsets~$J$ of $\{1,\ldots,r\}$.

  The \rhs\
  of~\eqref{equationdecompositioninlemmamultiplesetisHikamiset} is a
  disjoint union by Lemma \ref{lemmaHikamisets2}(iii) and the inclusion
  ``$\supset$'' directly follows from Lemma~\ref{lemmaHikamisets1}(i).

  Let us prove the converse inclusion. Let $n\in\Z$ satisfy
  $\ti J_n=J$.  We just need to find $\ul\in\fL$ and $\ue\in E$ such
  that $n\equiv\cN^\ul(\ue)\ [2P]$ and
  $J^\ul=J$.
  %
  
   According to the Chinese Remainder Theorem, since $2P = 2 p_1
   p_2\cdots p_r$, the congruence equation
     $n\equiv\cN^\ul(\ue)\ [2P]$ is equivalent to the system of equations
   \begin{align}
     \label{eqcongrjeq1}
     \eps_1 \ell_1\wh p_1 +\eps_2 \ell_2\wh p_2 +\cdots+\eps_r \ell_r\wh p_r
     &\equiv n -P \ [2p_1], \\[.8ex]
     \eps_1 \ell_1\wh p_1 +\eps_2 \ell_2\wh p_2 +\cdots+\eps_r \ell_r\wh p_r
     &\equiv n -P \ [p_j] \quad \text{for $2\le  j\le r$}.
       \label{eqcongrjge2}
  \end{align} 
  We will first check that the congruence equations~\eqref{eqcongrjge2}
  uniquely determine $\ell_2,\ldots,\ell_r$ as well as $\eps_2 \ell_2\wh p_2
  +\cdots+\eps_r \ell_r\wh p_r$.

  Suppose $j\ge2$.
  Since $\wh p_k$ is divisible by~$p_j$ for each $k\in\{1,\ldots,r\}\setminus\{j\}$, the
  congruence equation~\eqref{eqcongrjge2} is equivalent to
  $\eps_j \ell_j\wh p_j \equiv n\ [p_j]$;
  since $\pj{\wh p_j}$ is invertible in $\Z/p_j\Z$, the latter
  equation is equivalent to
  \beglab{eqpjepsjhj}
\pj{\eps_j \ell_j} = \pj{\wh p_j}\ii \pj{n} \ens \text{in}\ens \Z/p_j\Z.
\edla
The \rhs\ of~\eqref{eqpjepsjhj} can be written in a unique way as
$\pj{m_j}$ with $0\le m_j<p_j$,
and we note that $j\in\{2,\ldots,r\}\setminus J\iimp 0<m_j<p_j$,
whereas $j\in\{2,\ldots,r\}\cap J\iimp m_j=0$
(because $\ti J_n=J$).
Having $\ul\in\fL$ imposes the constraint $0\le \ell_j\le p_j$ and~$\ell_j$
even. Since~$p_j$ is odd, we find a unique solution $(\ell_j,\eps_j)$ for
$j\in\{2,\ldots,r\}\setminus J$, namely
\begla
(\ell_j,\eps_j) = (m_j,+1) \ens\text{if $m_j$ is even},
\quad
(\ell_j,\eps_j) = (p_j-m_j,-1) \ens\text{if $m_j$ is odd},
\edla
whereas~$\eps_j$ is left undetermined if $j\in\{2,\ldots,r\}\cap J$
and $\ell_j=0$ in that case.

Therefore, $\ell_2,\ldots,\ell_r$ are determined, as well as $M:=\eps_2 \ell_2\wh p_2
+\cdots+\eps_r \ell_r\wh p_r$.
Note that, according to our findings, 
\beglab{eqJsofar}
\big\{j\in \{2,\ldots,r\} \mid \ell_j \equiv0\ [p_j] \big\}
= \{2,\ldots,r\}\cap J.
\edla

We can now solve the first congruence equation:
since $\po{\wh p_1}$ is invertible in $\Z/2p_1\Z$,~\eqref{eqcongrjeq1} is equivalent to
  \beglab{eqtpo}
\po{\eps_1 \ell_1} = \po{\wh p_1}\ii \po{n - P - M} \ens \text{in}\ens \Z/2p_1\Z.
\edla
The \rhs\ of~\eqref{eqtpo} can be written in a unique way as
$\po{m_1}$ with $0\le m_1<2p_1$.
There are two cases, and in each of them we will determine~$\ell_1$ taking into account the
constraint $0\le \ell_1\le p_1$ due to $\ul\in\fL$:
\begin{itemize}
\item
either $1\notin J$: we then have $0<m_1<p_1$ or $p_1<m_1<2p_1$, and
we must take $(\ell_1,\eps_1)=(m_1,+1)$ in the former subcase and
$(\ell_1,\eps_1)=(2p_1-m_1,-1)$ in the latter one;
\item
or $1\in J$: we then have $m_1=0$ or $m_1=p_1$, and we must take
$\ell_1=m_1$ in both subcases (with~$\eps_1$ left undetermined).
\end{itemize}
The unique~$\ell_1$ that we just found is multiple of~$p_1$ if and
only if $1\in J$;
%
%
together with~\eqref{eqJsofar}, this yields $J^\ul=J$ and we can
confirm that $\ul\in\fL$. The proof is thus complete.
%
%
%
\end{proof}


\subsection{Generalized Hikami functions}   \label{secsHikfcns}

The $s$-Hikami functions~$m^sf^\uo$ were defined
in~\eqref{eqsupportsSuo}--\eqref{eqdefmsfuo}, based
on the definition~\eqref{eqdefwhpj} of the function~$\sN$ and the
definition~\eqref{eqdefCNuo} of~$\cN$.
Here,~$s$ can be any non-negative integer, but only the case $s\le r-3$
is relevant to this paper.
We now define functions~$m^sf^\uh$ for any $\uh\in\fH$ such that
$t_\uh=0$. Since $J^\uh=\emptyset$, by virtue of
Lemma~\ref{lemmaHikamisets1}(ii) there is a well-defined map
\begla
\tp{\cN^\uh}\ii \col \fS^\uh \to E.
\edla
Note that $\cN^\uh=P+\sN^\uh$ with
$\sN^\uh(\ue) := \eps_1 h_1 \wh p_1+\cdots +\eps_r h_r \wh p_r$
for any $\ue\in E$.

\begin{definition}
  For any $\uh\in\fH$ with $t_\uh=0$, we define the $s$-Hikami
  function $m^s f^{\uh}\col \Z\to \Z$  by
  \begin{equation}   \label{eqdefmsfuh}
m^sf^{\uh}(n) := \begin{cases}
      -\pi(\ue)\big(\sN^{\uh}(\ue)\big)^s
      \;\, &\text{if $n\in\fS^\uh$,\;
        with $\ue =\tp{\cN^\uh}\ii\big(\tp{n}\big)$} \\[.8ex]
      \ens 0    &\text{if $n\in\Z\setminus\fS^\uh$,}
    \end{cases}
  \end{equation}
  with the notation $\pi(\ue)=\eps_1\cdots\eps_r$.
  %
  %
\end{definition}

As a particular case, we may take $\uh=\uo$: one always have
$\uo\in\fH$ and $t_\uo=0$, and one then
recovers the function~$m^sf^\uo$ of~\eqref{eqsupportsSuo}--\eqref{eqdefmsfuo}.

\begin{lemma}
  Let $h\in\fH$ have $t_\uh=0$.
      \medskip
      
      \noindent \emph{\textbf{(i)}}
      The function $m^sf^\uh$ is $2P$-periodic and even or odd, of same parity as $r-s$.
      The set $\cN^\uh(E)\subset \Z$ is a system of representatives of
      its support $\!\!\!\mod 2P\Z$,
      and there is an identity between Laurent polynomials of $\Z[z,z\ii]$:
\begin{equation}\label{equationgenratingfunctionHikami}
  \sum_{n\in\cN^{\uh}(E)} m^sf^{\uh}(n) z^n =
  -  z^P 
  \Big(z\frac{d}{dz}\Big)^{\!s}
  \Big(\prod_{j=1}^r (z^{h_j\whp_j}-z^{-h_j\whp_j}) \Big).
\end{equation}
  %
      
\noindent \emph{\textbf{(ii)}}
The \rhs\ of~\eqref{equationgenratingfunctionHikami} can also be
written as
\beglab{eqresultofLeibniz}
\sum_{ \substack{ s_1,\ldots,s_r\ge0\ \text{s.t.} \\[.3ex] s_1+\cdots+s_r=s } } \,
\frac{ - s! \, z^P }{ s_1! \cdots s_r! } \, \prod_{j=1}^r
\Big(z\frac d{dz}\Big)^{\!s_j} (z^{h_j\whp_j}-z^{-h_j\whp_j}).
\edla
      
\noindent \emph{\textbf{(iii)}}
      Suppose $0\le s < r$. Then the mean value of $m^sf^\uh$ is zero.
      \medskip
      
      \noindent \emph{\textbf{(iv)}}
      Suppose $s=0$ and let $\uh' \in\fH$ be such that $t_{\uh'}=0$.
      If $\uh\sim\uh'$, then $m^0f^{\uh} = m^0f^{\uh'}$.
    \end{lemma}

\begin{proof}
  \textbf{(i):} Parity is obvious, since~$\sN$ is odd.
  Then, use Lemma~\ref{lemmaHikamisets1}(ii) and
  get~\eqref{equationgenratingfunctionHikami} by mimicking
  the passage from~\eqref{eqoriggenser} to~\eqref{eqcomputgenpol}.
  \textbf{(ii):} Leibniz rule.
  \textbf{(iii):} Evaluate~\eqref{eqresultofLeibniz} at
  $z=1$: if $r>s$, then at least one of~$s_j$'s is~$0$ and the corresponding factor vanishes.
  \textbf{(iv):} The generating function of~$m^0f^\uh$ is
  \begla
  \sum_{n\in\Z} m^0f^\uh(n) z^n = \sum_{k\in\Z} \sum_{n\in\cN^h(E)}
  f^\uh(n) z^{n+2kP}
  = \sum_{k\in\Z} z^{2kP} \gP{\uh,0}(z) \in\Z[[z,z\ii]],
  \edla
  where $\gP{\uh,s}(z)\in\Z[z,z\ii]$ is the Laurent polynomial~\eqref{equationgenratingfunctionHikami}.
  If $\uh$ and~$\uh'$ satisfy~\eqref{equhsimuhp}, then
  $z^{h'_j\whp_j}-z^{-h'_j\whp_j} =
  -z^P(z^{h_j\whp_j}-z^{-h_j\whp_j})$ for $j\in J$, whereas these two
  factors are identical for $j\in\{1,\ldots, r\}\setminus J$,
  thus $\gP{\uh',0}(z) = (-z^P)^{|J|}\gP{\uh,0}(z)$.
Since~$|J|$ is even, this implies that $m^0f^\uh$ and~$m^0f^{\uh'}$ have the
same generating function.
\end{proof}


What about the case when we do not assume $t_\uh=0$? Next section will
require ``generalized Hikami functions'' associated with $s=0$ and arbitrary
$\uh\in\fH$, but Lemma~\ref{lemmaHikamisets1}(ii) shows that in general
the map $\tp{\cN^\uh}$ is no longer injective, so we need to modify
the definition~\eqref{eqdefmsfuh}.

\begin{definition}   \label{defgJuh}
  For any $\uh\in\fH$ and any subset~$J$ of $\{1,\ldots,r\}$ such
  that $J^\uh\cap J=\emptyset$, we define the generalized Hikami
  function $g_J^\uh\col \Z\to \{-1,0,1\}$  by
  \begin{equation}   \label{eqdefgJuh}
g_J^\uh(n) := \begin{cases}
      \prod\limits_{j\in J} \eps_j
      \;\, &\text{if $n\in\fS^\uh$,\;
        with any $\ue\in E$ such that $\tp{\cN^\uh(\ue)}=\tp{n}$} \\[.8ex]
      \ens 0    &\text{if $n\in\Z\setminus\fS^\uh$.}
    \end{cases}
  \end{equation}
\end{definition}
Note that the definition~\eqref{eqdefgJuh} makes sense because we have
assumed $J\subset \compl J^\uh$, thus Lemma~\ref{lemmaHikamisets1}(ii)
implies that the restriction $\ue_{\mid J}$ is determined for any
$n\in \fS^\uh$.
In the particular case $J=\{1,\ldots,r\}$, we recover $m^0 f^\uh$ as
$-g_{\{1,\ldots,r\}}^{\,\uh}$.

\begin{lemma}   \label{lemmagsmeanvalue}
  Suppose $h\in\fH$, $J\subset\{1,\ldots,r\}$ and $J^\uh\cap J=\emptyset$.
      \medskip
      
      \noindent \emph{\textbf{(i)}}
      The function $g_J^\uh$ is $2P$-periodic and even or odd, of same parity as $|J|$.
      \medskip
      
\noindent \emph{\textbf{(ii)}}
    There is an identity in the quotient ring $\Z[z]/(z^{2P}-1)$:
\beglab{equationgeneratinggs}
  2^{t_\uh} \cdot\! \sum_{n \!\!\mod 2P} g_J^\uh(n) z^n \equiv
z^P \cdot
  \prod_{j\in \ccompl J} (z^{h_j\whp_j}+z^{-h_j\whp_j})
    \cdot \prod_{j\in J} (z^{h_j\whp_j}-z^{-h_j\whp_j})
  \! \mod (z^{2P}-1).
\end{equation}
  %
      
\noindent \emph{\textbf{(iii)}}
      Suppose $J\neq\emptyset$. Then the mean value of $g_J^\uh$ is zero.
    \end{lemma}

\begin{proof}
\textbf{(i)} is obvious and \textbf{(iii)} follows from~\textbf{(ii)}
by evaluation at $z=1$.
Let us prove~\textbf{(ii)}: we mimic the passage
from~\eqref{eqoriggenser} to~\eqref{eqcomputgenpol} and write
 the \rhs\ of~\eqref{equationgeneratinggs} as
 \beglab{eqnewnumber}  
  z^P \cdot
  \prod_{j\in \ccompl J} \Big( \sum_{\eps\in\{\pm1\}}
  z^{\eps_jh_j\whp_j} \Big)
    \cdot \prod_{j\in J} \Big( \sum_{\eps\in\{\pm1\}} \eps_j
    z^{\eps_jh_j\whp_j} \Big)
    = \sum_{\ue\in E} \Big( \prod_{j\in J} \eps_j \Big)
    z^{P + \sum\limits_{j=1}^r \eps_jh_j\whp_j}
    %
    %
= \sum_{\ue \in E} g_J^\uh\big(\cN^\uh(\ue)\big) z^{\cN^\uh(\ue)}.
  \edla \bblack
This is a polynomial of $\Z[z,z\ii]$ that we can project to the quotient
    ring
    $\Z[z,z\ii]/(z^{2P}-1) = \Z[z]/(z^{2P}-1)$:
    this amounts to replacing the power $\cN(\ue)$ by $\tp{\cN(\ue)}$
    and Lemma~\ref{lemmaHikamisets1}(ii) thus yields
    $2^{t_\uh}\cdot \!\sum\limits_{n \in \cN^\uh(E) \!\!\mod 2P}
    g_J^\uh(n) z^n \! \mod (z^{2P}-1)$,
    which is the \lhs\ of~\eqref{equationgeneratinggs}.
\end{proof}

\subsection{Discrete Fourier transforms of $m^sf^\uo$ and
  generalized Hikami functions}
\label{secDFTsHik}

Recall that, if $f$ is a
$2P$-periodic function from~$\Z$ to~$\C$, then according to footnote~\ref{ftnUM} its DFT is the $2P$-periodic
function defined by
  \beglab{eqDFTwhfell}
  n\in \Z \mapsto \wh{f}(n) := \frac{1}{\sqrt{2P}}
  \sum_{\ell \!\!\mod 2P} e\Big[\!-\frac{\ell n}{2P}\Big] f(\ell),
  \quad\text{where}\ens e[x]:=e^{2\pi ix}.
  \edla
  In other words,
  \beglab{eqDFTwhfRGP}
\sqrt{2P}\wh{f}(n) = \; \text{evaluation of}\;
\sum\limits_{\ell\!\!\mod 2P} f(\ell) z^\ell \!\!\mod(z^{2P}-1)
\; \text{at the root of unity}\; e\Big[\!-\frac n{2P}\Big].
\edla
Note that we are using here what may be called the ``Reduced
Generating Polynomial'' of~$f$, an element of the quotient ring
$\C[z]/(z^{2P}-1)$.
  
  The first part of this section aims at computing the DFT of the $s$-Hikami function
  $m^s f^\uo$.
  More precisely, we need $\rwh{m^sf^\uo\circ\cT}$ with~$\cT$ as in~\eqref{eqdefcTm},
  i.e.\ $\cT=\ID_{\Z}$ or $\ID_{\Z}-P$ according as $r$ is odd or
  even.
  The first step is

\begin{lemma}\label{lemmageneratingfunctionDFT}
  For any $s\ge0$, the DFT of $m^s f^\uo \circ\cT$ is given by
  \begin{multline}   \label{eqDFTmsfuofirst}
\rwh{m^s f^\uo \circ\cT}(n) = 
    \ka_n \cdot \Big[
    \Big(z\frac{d}{dz}\Big)^{\!s}
    \Big( \prod_{j=1}^r(z^{\whp_j}-z^{-\whp_j})\Big)
    \Big]_{z=e[-\frac{n}{2P}]}
    \\[.8ex]
    = \sum_{ \substack{ s_1,\ldots,s_r\ge0\ \text{s.t.} \\[.3ex] s_1+\cdots+s_r=s } } \,
\frac{ s! \, \ka_n }{ s_1! \cdots s_r! } \, \prod_{j=1}^r \Big[
\Big(z\frac d{dz}\Big)^{\!s_j} (z^{\whp_j}-z^{-\whp_j})
    \Big]_{z=e[-\frac{n}{2P}]}
    \end{multline}
    for all $n\in\Z$, with the notation $\ka_n := \frac{(-1)^{rn+1}}{\sqrt{2P}}$.
  %
%
%
%
\end{lemma}

\begin{proof}
  For any $2P$-periodic function $f\col \Z \to \C$, the DFT of
  $g:=f\circ(\ID_{\Z}-P)$ is $\wh g(n)=(-1)^n\wh f(n)$,
  whence $\rwh{f\circ\cT}(n) = (-1)^{(r+1)n}\wh f(n)$.
  %
  %
  The result thus follows from~\eqref{equationgenratingfunctionHikami}--\eqref{eqresultofLeibniz}.
  %
\end{proof}

We now assume $s\in\{0,\ldots,r-3\}$ and set out to compute the \rhs\
of~\eqref{eqDFTmsfuofirst}, first when~$n$ belongs to~$\Mss$, and
then when
$n\in\bigsqcup\limits_{\ul\in\fL \,\text{s.t.}\, t_\ul\le s}\fS^\ul$,
with reference to the decomposition of~$\Z$ given by Lemma~\ref{lemmaresiduesets}.

\begin{lemma}\label{lemmamsf=0}
Let $s\in\{0,\ldots,r-3\}$.  Then the function $\rwh{m^sf^{\uo}
    \circ \cT}$ vanishes on $\Mss$.
%
%
\end{lemma}

\begin{proof}
Suppose $n\in\Mss$, i.e. the subset~$\ti J_n$ of~$\{1,\ldots,r\}$
has cardinality $>s$.
Pick any term labelled by $\us=(s_1,\ldots,s_r)$ in the \rhs\
of~\eqref{eqDFTmsfuofirst}; the condition $s_1+\cdots+s_r=s$ implies
that $L:=\big\{ j \in \{1,\ldots,r\} \mid s_j=0 \big\}$ has cardinality
$\ge r-s$,
whence $L\cap \ti J_n \neq\emptyset$.
Now pick any $j\in L\cap \ti J_n$: because~$n$ is multiple of~$p_j$,
$e\big[\!-\frac n{2P}\big]$ is a root of the corresponding factor
%
in the term associated with~$\us$.
Therefore all terms are~$0$.
%
  %
\end{proof}

%
%


In view of the decomposition of~$\Z$ given by
Lemma~\ref{lemmaresiduesets}, we thus have
\beglab{eqrwhmsfuoasasumovul}
\rwh{m^sf^{\uo} \circ \cT} = \sum_{\ul\in\fL \,\text{s.t.}\, t_\ul\le s} \rwh{m^sf^{\uo} \circ \cT}\cdot \indul
\edla
(with the notation $\indfS$ for the indicator function of a
subset~$\fS$ of~$\Z$).
We now give ourselves
\begla
\ul\in\fL \ens\; \text{such that}\ens
0 \le t_\ul = |J^\ul| \le s
\edla
and focus, for the rest of the computation, on the values of
$\rwh{m^sf^{\uo} \circ \cT}$ on~$\fS^\ul$.
%

Let us consider arbitrary $n\in\fS^\ul$ and $\uh\in\fH$ such that $\uh\sim\ul$.
Recall that, thanks to Lemmas~\ref{lemmaHikamisets1}(ii) and~\ref{lemmaHikamisets2}(ii),
%
we have $\fS^\ul=\fS^\uh$ and $\tp{n}$ can be written as
\beglab{eqellcNulue}
\tp n= \tp{\cN^\uh(\ue)}
\quad \text{for some} \ens \ue\in E,
\edla
where the restriction of $\ue\in E$ to $\compl J^\uh$ is uniquely
determined (and is thus a $2P$-periodic function of~$n$), whereas
its restriction to~$J^\uh$ is free (there are $2^{t_\uh}$ possibilities for~$\ue$).
By Lemma~\ref{lemmageneratingfunctionDFT},
\begin{align*}
  \rwh{m^sf^\uo\circ\cT}(n)
&= \sum_{ \substack{ s_1,\ldots,s_r\ge0\ \text{s.t.} \\[.3ex] s_1+\cdots+s_r=s } } \,
\frac{ s! \, \ka_n}{ s_1! \cdots s_r! } \, \prod_{j=1}^r \Big[
\whp_j^{\hspace{.3ex}s_j} (z^{\whp_j}-(-1)^{s_j} z^{-\whp_j})
    \Big]_{z=e[-\frac{n}{2P}]}
\\[.8ex]
&= \sum_{ \substack{ s_1,\ldots,s_r\ge0\ \text{s.t.} \\[.3ex] s_1+\cdots+s_r=s } } \,
\frac{ s! \, \ka_n}{ s_1! \cdots s_r! } \, \prod_{j=1}^r
\whp_j^{\hspace{.3ex}s_j} (e^{-i\pi n/p_j}-(-1)^{s_j} e^{i\pi n/p_j}).
  \end{align*}
  Clearly, the factors associated with $j\in \ti J_n$ such that $s_j$ is even
  vanish (because, for each of these,~$n$ is multiple of~$p_j$).
  Now, by Lemma~\ref{lemmaHikamisets1}(i), we have $J^\uh = \ti
  J_{\cN^\uh(\ue)} = \ti J_n$
  (in view of~\eqref{eqellcNulue}),
  thus we can restrict the summation to
  \begin{gather}
    \label{eqdefuSul}
  \uS_s^\uh :=
  \{\, \us=(s_1,\ldots,s_r) \in \Z_{\ge0}^r \mid
  s_1+\cdots+s_r = s
  \;\, \text{and}\;\, Ev_\us \cap J^\uh=\emptyset \,\} 
  \\
  \intertext{with the notation}
  \label{eqdefOdusEvus}
Ev_\us := \big\{ j \in \{1,\ldots, r\} \big|\ \text{$s_j$ is even}
\big\}
\quad \text{for any}\ens \us \in \Z_{\ge0}^r.
\\[-2ex]
\intertext{We get}
    \notag
    \rwh{m^sf^\uo\circ\cT}(n) =
  \sum_{ \us \in \uS_s^\uh }
\frac{ s! \, \ka_n}{ s_1! \cdots s_r! } \, \prod_{j=1}^r \whp_j^{\hspace{.3ex}s_j} 
\, \prod_{ j\in \cEv_\us }
(e^{-i\pi n/p_j}+e^{i\pi n/p_j})
\, \prod_{ j\in Ev_\us }
(e^{-i\pi n/p_j}-e^{i\pi n/p_j}).
%
%
\end{gather}
Notice that, for each $\us\in\uS_s^\uh$, 
%
%
\beglab{ineqEvusencadr}
%
%
 Ev_\us\subset\compl J^\uh
  \quad \text{and} \quad |Ev_\us| \ge r-s \ge 3
\edla
(because $j\in \CEv_\us \iimp s_j\ge1$, thus
$s=s_1+\cdots+s_r \ge |\CEv_\us|$). 
In particular $Ev_\us$ is never empty, and
   \beglab{eqparitycardEvus}
   s \equiv \sum_{j\in Ev_\us} 0 + \sum_{j\in \cEv_\us} 1
   \equiv r-|Ev_\us| \! \mod 2, 
   \edla
   i.e.\ $|Ev_s|$ and $r-s$ have same parity.

We pursue the computation by observing that, in view of its definition
in Lemma~\ref{lemmageneratingfunctionDFT}, $\ka_n$ only depends on~$\uh$, not on
$n\in\fS^\uh$.
Indeed, it only depends on the parity of~$n$, and
$n\equiv\cN^\uh(\ue)\, [2]$
with $\cN^\uh(\ue')-\cN^\uh(\ue)= \sum (\eps'_j-\eps_j)h_j\wh p_j
\equiv 0 \, [2]$, because $\eps'_j-\eps_j$ is always even.
Thus,
\beglab{eqDFTasprodcossin}
    \rwh{m^sf^\uo\circ\cT}(n) =
  \sum_{ \us \in \uS_s^\uh } \wt K_\uh^\us
\, \prod_{ j\in \cEv_\us }
\cos(\pi n/p_j)
\, \prod_{ j\in Ev_\us }
\sin(\pi n/p_j)
\edla
with
$\displaystyle \wt K_\uh^\us := \frac{(-1)^{|Ev_\us|} \, 2^r s! \, \ka_n}{ s_1!
  \cdots s_r! } \, \prod_{j=1}^r \whp_j^{\hspace{.3ex}s_j}$ for any $n\in\fS^\uh$.

We will now show that~\eqref{eqDFTasprodcossin} can be rewritten as
\beglab{eqDFTmsfuoKp}
    \rwh{m^sf^\uo\circ\cT}(n) =
  \sum_{ \us \in \uS_s^\uh } K_\uh^\us
\, \prod_{ j\in Ev_\us }
\eps_j
\edla
with coefficients $K_\uh^\us$ independent of~$n$
(whereas the \rhs\ depends on~$n$ through the restriction of~$\ue$ to $Ev_\us$,
which is determined by~\eqref{eqellcNulue}).

For a given~$j$,
$\cos(\pi n/p_j)$ and $\sin(\pi n/p_j)$ depend on~$n$ only through
$\tpj{n}$ and
\begla
n \equiv \cN^\uh(\ue)
\equiv \pm
P + \sum_{i=1}^r \eps_i h_i \whp_i \!\mod 2p_j,
\edla
hence we just need to deal with
%
%
$\cos\big(\pi\cN^\uh(\ue)/p_j\big)$
or 
%
%
$\sin\big(\pi\cN^\uh(\ue)/p_j\big)$
according as $j\in \CEv_\us$ or $j\in Ev_\us$.
This quantity does not depend
on~$\eps_i$ for $i\in\{1,\ldots,r\}\setminus\{j\}$ (because switching
the sign of~$\eps_i$ amounts changing $\cN^\uh(\ue)$ by adding to it $\pm 2h_i\wh p_i$, which is
a multiple of $2p_j$), thus it is a function of~$\eps_j$ only; now,
that function is even or odd in~$\eps_j$:
$\cos\big(\pi\cN^\uh(\ue)/p_j\big) =
\cos\big(\pi\cN^\uh(-\ue)/p_j\big)$ is even in~$\eps_j$ and thus does not depend on~$\ue$ at all,
whereas $\sin\big(\pi\cN^\uh(\ue)/p_j\big) =
-\sin\big(\pi\cN^\uh(-\ue)/p_j\big)$ is odd in~$\eps_j$ and is thus a
multiple of~$\eps_j$. This yields~\eqref{eqDFTmsfuoKp} with
\beglab{eqdefKulus}
K_\uh^\us = \wt K_\uh^\us \, \prod_{ j\in \cEv_\us }
\cos\big(\pi\cN^\uh(\uo)/p_j\big)
\, \prod_{ j\in Ev_\us }
\sin\big(\pi\cN^\uh(\uo)/p_j\big).
\edla

Thus~\eqref{eqDFTmsfuoKp} is proved. We now observe that, in view of~\eqref{eqellcNulue}, $\prod_{ j\in
  Ev_\us }\eps_j$ is nothing but
$g_{Ev_\us}^\uh(n)$, with reference to
Definition~\ref{defgJuh}. Therefore, our result is

\begin{lemma} \label{lemDFTonfSul}
  For every $s\in\{0,\ldots,r-3\}$ and $\ul\in\fL$ such that $t_\ul\le s$, the restriction of
  the DFT of $m^sf^{\uo} \circ \cT$ to~$\fS^\ul$ is
\beglab{eqresultforDFTonfSul}
\rwh{m^sf^{\uo} \circ \cT}\cdot \indul = 
\sum_{\us\in\uS_s^\uh} K_\uh^\us \, g_{Ev_\us}^\uh
\quad \text{for any $\uh\in\fH$ such that $\uh\sim\ul$}
\edla
with $\uS_s^\uh$ as in~\eqref{eqdefuSul} and
$K_\uh^\us$ as in~\eqref{eqdefKulus}.
  \end{lemma}

  Note that for each $\us\in\uS_s^\uh$ we have $Ev_\us\neq\emptyset$,
  thus Lemma~\ref{lemmagsmeanvalue} shows that $g_{Ev_\us}^\uh$ has
  zero mean value.
  One can check that for any $\uh\sim \ul$, $\uS_s^\uh = \uS_s^\ul$,
    but different choicies of~$\uh$ may lead different decompositions of
    $\rwh{m^sf^{\uo} \circ \cT}\cdot \indul$ (because the
    constants~$K_\uh^\us$ and the functions $g_{Ev_\us}^\uh$ depend
    on~$\uh$), and this flexibility will prove useful at the end of
    next section.
    Choosing $\uh=\ul$ we obtain,
    as a direct consequence of~\eqref{eqrwhmsfuoasasumovul} and
    Lemma~\ref{lemDFTonfSul}: 
  
  \begin{proposition}\label{propositionDFTcomputation}
    For every $s\in\{0,\ldots, r-3\}$, we have
    \beglab{eqdecrwhmsfuocTovulus}
    \rwh{m^sf^\uo\circ\cT} = \sum_{\ul\in\fL \,\text{s.t.}\, t_\ul\le
      s}\, \,
    \sum_{\us\in\uS_s^\ul} \, K_\ul^\us \, g_{Ev_\us}^\ul
    \edla
    where each~$g_{Ev_\us}^\ul$ is supported on the Hikami set~$\fS^{\ul}$
    and has zero mean value.
  \end{proposition}

  \label{apppflemDFTmsfuo}

We conclude this section by describing the DFT of the functions
$g_J^\uh$.
Their reduced generating polynomials are given 
by~\eqref{equationgeneratinggs}; 
thanks to~\eqref{eqDFTwhfRGP}
%
%
%
and computations similar to those of this section (but much simpler), one finds
%
%
\begin{lemma}   \label{lemDFTofgJuh}
  For any $J\subset \{1,\ldots,r\}$ such that $|J|\ge3$
  and $\uh\in\fH$  
such that $J^\uh \subset \compl J$,
  the function $g_J^\uh$ has a DFT supported in the disjoint union
  of all the Hikami sets~$\fS^{\ul'}$ with $\ul'\in\fL$ such that
  %
  %
  $J^{\ul'}\subset \compl J$.
  The restriction of this DFT to such a set~$\fS^{\ul'}$ is of the form
    \begin{equation}
      \rwh{\fakeheight{$g_J^\uh$}{$h_J$}} \cdot \indulp
      = \Ga_{\!J}(\uh,\ul')\, g_J^{\ul'}
    \quad \text{for some constant}\;\, \Ga_{\!J}(\uh,\ul')
    \end{equation}
  and $\rwh{\fakeheight{$g_J^\uh$}{$h_J$}}$ is thus a linear
  combination of these functions $g_J^{\ul'}$.
  \end{lemma}

  \begin{proof}[Sketch of proof]
    Given $J\subset \{1,\ldots,r\}$ such that $|J|\ge3$
    and $\uh\in\fH$ such that $J^\uh \subset \compl J$,
    in view of the identity~\eqref{equationgeneratinggs} satisfied by
    the reduced generating polynomial of~$g_J^\uh$,
    \eqref{eqDFTwhfRGP} yields
    \begin{align}
      \notag
    \sqrt{2P} \, \rwh{\fakeheight{$g_J^\uh$}{$h_J$}}(n) &= (-1)^n \, 2^{-t_\uh}
\, \prod_{ j\in \ccompl J }
(e^{-i\pi h_j n/p_j}+e^{i\pi h_j n/p_j})
\, \prod_{ j\in J }
(e^{-i\pi h_j n/p_j}-e^{i\pi h_j n/p_j})
      \\[.8ex]
      \label{eqrwhgJuhprocsossin}
 &= (-1)^{n+|J|} \, 2^{r-t_\uh}
\, \prod_{ j\in \ccompl J }
\cos(\pi h_j n/p_j)
\, \prod_{ j\in J }
\sin(\pi h_j n/p_j).
    \end{align}
    Therefore, the support of $\rwh{\fakeheight{$g_J^\uh$}{$h_J$}}$ is contained in $\{ n\in
    \Z \mid \ti J_n \subset \compl J\}$.
    Since the inclusion $\ti J_n \subset \compl J$ entails
    $|\ti J_n| \le r-|J|$, we can use Lemma~\ref{lemmaresiduesets}
    with $s=r-|J|$: we obtain that the support of $\rwh{\fakeheight{$g_J^\uh$}{$h_J$}}$ is
    contained in the disjoint union of all the Hikami
    sets~$\fS^{\ul'}$ with $\ul'\in\fL$ such that $|J^{\ul'}| \le
    r-|J|$, 
    and, thanks to Lemma~\ref{lemmaHikamisets1}(i),
    we can even restrict to those such that $J^{\ul'}\subset \compl J$.

  Take~$n$ in one of these sets~$\fS^{\ul'}$ and write $n\equiv
  \cN^{\ul'}(\ue) \ [2P]$ with some $\ue\in E$:
  the restriction $\ue_{\mid J^{\ul'}}$ is free but
  $\ue_{\mid \ccompl J^{\ul'}}$ is determined;
  in particular, $\ue_{\mid J}$ is determined.
  In formula~\eqref{eqrwhgJuhprocsossin}, each of the $\cos$ or $\sin$ factors
  depends only on $\tpj{n}$; arguing exactly as in the proof
  of~\eqref{eqDFTmsfuoKp}, we find that each $\cos$ is proportional
  to~$1$ and each $\sin$ is proportional to~$\eps_j$, with
  proportionality constants depending only on~$\ul'$ (not on~$\ue$,
  i.e.\ not on~$n$),
and the products of the $\eps_j$'s with $j\in J$ is precisely
$g_J^{\ul'}(n)$.
    \end{proof}
  
Since $m^sf^\uo\circ\cT$ is always even or odd of
same parity as $r-s$ (because that is the case for $m^sf^\uo$ itself),
the DFT of $\rwh{m^sf^\uo\circ\cT}$
is none other than $(-1)^{r-s} m^sf^\uo\circ\cT$;
 putting together Proposition~\ref{propositionDFTcomputation} and
 Lemma~\ref{lemDFTofgJuh} we thus obtain

 \begin{corollary}   \label{corcVs}
   For each $s\in\{0,\ldots, r-3\}$, both $m^s f^\uo\circ\cT$ and its DFT
   belong to the $\C$-vector space
\begla
\cVs := \Span\big\{ g_J^\ul \ \big|\ 
J\subset \{1,\ldots,r\}, \;\, |J|\ge r-s, \;\, |J| \equiv r-s \ [2],\;\,
 \ul\in\fL, \;\, J^\ul\subset \compl J \big\}.
\edla
All the elements of~$\cVs$ are zero mean value $2P$-periodic functions, even or odd of
same parity as $r-s$.
Moreover, the space~$\cVs$ is stable under DFT.
 \end{corollary}

 \begin{proof}
   We show that $\rwh{m^sf^\uo\circ\cT} \in \cVs$ by
   rewriting~\eqref{eqdecrwhmsfuocTovulus} as
   \begla
    \rwh{m^sf^\uo\circ\cT} = \sum_{\ul\in\fL \,\text{s.t.}\, t_\ul\le
      s}\, \,
    \sum_{J\subset\{1,\ldots,r\}} \, \uK_\ul^J \, g_J^\ul
    \quad \text{with}\ens
    \uK_\ul^J := \sum\limits_{\us\in\uS_s^\ul \,\text{s.t.}\, Ev_\us=J} K_\ul^\us.
    \edla
    Each constant~$\uK_\ul^J$ vanishes unless
    \eqref{ineqEvusencadr}--\eqref{eqparitycardEvus} hold, which
    implies
    $J \subset \compl J^\ul$, $|J|\ge r-s$ and $s\equiv r-|J| \ [2]$.
    We thus find
   \begla
   \rwh{m^sf^\uo\circ\cT} =
   \sum_{J \,\text{such that}\atop |J|\ge r-s,\  |J|\equiv r-s \, [2]}\, \,
    \sum_{\ul\in\fL \,\text{such that}\atop t_\ul\le s,\ J^\ul\subset\ccompl J}
    \, \uK_\ul^J \, g_J^\ul
    \, \in \, \cVs
    \edla
    (note that the condition $t_\ul\le s$ in the latter summation can
    be omitted, since $\uK_\ul^J \neq0 \iimp t_\ul = |J^\ul| \le r-|J|$
    and we need $r-|J|\le s$).    

    The functions in~$\cVs$ are all $2P$-periodic and of same parity
    as $r-s$, since this is the case for $g_J^\ul$ when $|J|\equiv
    r-s\ [2]$ by Lemma~\ref{lemmagsmeanvalue}(i); since $|J|\ge r-s\ge3$
    for each $g_J^\ul\in \cVs$, we get zero mean value by Lemma~\ref{lemmagsmeanvalue}(iii).

    We easily obtain that~$\cVs$ is stable under DFT from Lemma~\ref{lemDFTofgJuh}.
    In particular, $m^sf^\uo\circ\cT$, being the DFT of $(-1)^{r-s}
    \rwh{m^sf^\uo\circ\cT}$, is in~$\cVs$ too.
    %
%
%
 \end{proof}

 Note that we also have $g_J^\uh\in \cVs$ for every $J\subset \{1,\ldots,r\}$
 such that $|J|\ge r-s$ and $|J| \equiv r-s \ [2]$ 
 and every $\uh\in\fL$ such that $J^\uh\subset \compl J$,
 by the same argument as for $m^sf^\uo\circ\cT$ (using parity and Lemma~\ref{lemDFTofgJuh}).

\subsection{Consequences for some partial theta series}
\label{secConseqTh}

\begin{proposition}   \label{propConseqTh}
Let $s\in\{0,\ldots, r-3\}$.
For every $g \in \cVs$, the quantum set $\sQ_{g,2P}$ as defined in~\eqref{eqdefsQfM} is
  all of~$\Q$,
  i.e.\ 
  the periodic function
  $m\in\Z\mapsto g(m) e^{i\pi m^2\al/(2P)}$
  has zero mean value and the non-tangential limits
    $\lim\limits_{\tau\to\al} \Theta(\tau;\nu,g,2P)$ thus exist 
  for all $\al \in \Q$ and $\nu\in\Z_{\ge0}$.
\end{proposition}

\begin{proof}
  If $s\equiv r+1\ [2]$, then all functions in~$\cVs$ are odd and
  the conclusion follows from Remark~\ref{remtrivsQ}.

  We now suppose
  that~$s$ and~$r$ have same parity, thus all functions $g\in\cVs$
  are even.
  Define
  \begla
  \oQs := \bigcap\limits_{g\in\cVs}\sQ_{g,2P}.
  \edla
  We will prove that $\oQs = \Q$ by using the following
  characterization (consequence
  of~\eqref{eqexistlimal}--\eqref{eqdefsQfM}):
  \begla
   \oQs = \{\, \al\in\Q \mid
  \text{for each $g\in\cVs$,\; $\Theta(\tau; 1,g,2P)$ has a limit as $\tau\to\al$}
  \,\}.
    \edla
    Since $0\in\oQs$ (by~\eqref{eqdefsQfM}, because
    each $g\in\cVs$ has zero mean value), it is sufficient to
    prove that $\oQs$ is invariant under (i) the unit translation
    $\al\in\Q\mapsto\al+1$ and (ii) the negative inversion
    $\al\in\Q\setminus\{0\}\mapsto-\al\ii$.

    \smallskip

    \noindent (i)
    Suppose $\al\in\oQs$. Every $g\in\cVs$ can be written as a
    linear combination
    of functions~$g_J^\ul$ belonging to~$\cVs$;
    for each of them, \eqref{eqdefThetafjM} and~\eqref{eqfSulSsigoul}
    yield
    \[ \Theta(\tau+1; 1,g_J^\ul,2P) = e^{-2\pi i
      S_{\sig_1(\ul)}}\Theta(\tau; 1,g_J^\ul,2P), \]
      whence the existence of $\lim_{\tau\to\al}\Theta(\tau+1; 1,g,2P)$
      follows.
      Therefore $\al+1\in \oQs$.
    
    \smallskip

    \noindent (ii)
    Suppose that $0\neq\al\in\oQs$. For every $g\in\cVs$, since~$g$ is
    even and has zero mean value, we can
    apply~\eqref{equationShigherdepth} with $j=1$:
    \begla 
      \Theta(\tau;1,g,2P) \mp i^{\frac12} \tau^{-\frac{3}{2}} \Theta(-\tau^{-1};1,\wh g,M)
= \bS^{\frac\pi2\mp \epsilon}\wt\Theta_{1,g,0,2P}(\tau).
\edla
Since $\wh g\in\cVs$ and $\al\in\oQs$, the second term of the
\lhs\ has a limit as $\tau\to-\al\ii$.
So does the \rhs\ if $-\al\ii>0$ and we consider the lateral sum
$\bS^{\frac\pi2- \epsilon}$, or if $-\al\ii<0$ and we consider the
lateral sum $\bS^{\frac\pi2+ \epsilon}$.
Thus, in all cases, $\Theta(\tau;1,g,2P)$ itself has a limit as
$\tau\to-\al\ii$. Therefore $-\al\ii\in\oQs$.
\end{proof}


We now give a result that is crucial to our proof
of Witten's conjecture: the point is that, in our decomposition of the
DFT of $m^sf^\uo\circ\cT$, some pieces do not contribute of the
non-tangential limits we are interested in.
   
\begin{proposition}\label{propositionlimitvanish}
Let $s\in \{0,\ldots,r-3\}$. Let $\nu\in\{0,\ldots,r-s-2\}$ satisfy $\nu\equiv r-s-1\ [2]$.
  Then, for any $\uh\in\mathfrak{H}$,
  \begin{equation}   \label{implicvanish}
    0<\sum\limits_{j=1}^r \frac{h_j}{p_j}<1
    \Imp
    \lim\limits_{\tau\rightarrow 0}\Theta(\tau;\nu,\rwh{m^sf^{\uo}\circ\cT}\cdot \induh,2P)
    = 0.
  \end{equation}
  \end{proposition}

  Note that the conclusion in~\eqref{implicvanish} depends only on the
  class of~$\uh$ modulo the equivalence relation~$\sim$ that we have
  introduced before the statement of Lemma~\ref{lemmaHikamisets2};
  indeed, we can write $\fS^\uh=\fS^\ul$ with a uniquely determined
  $\ul\in\fL$. However, the premise of~\eqref{implicvanish} does
  depend on~$\uh$ itself and not only on its equivalence class.
  It is here that we use the flexibility provided by Lemma~\ref{lemDFTonfSul}.
  
  \begin{proof}[Proof of Proposition~\ref{propositionlimitvanish}]
    Let~$s$, $\nu$ and~$\uh$ be as in the statement, with~$\uh$
    satisfying the premise of~\eqref{implicvanish}, which we rewrite
    as
    \beglab{ineqcrucialhj}
    0 < \sum_{j=1}^r h_j \, \wh p_j < P.
    \edla
    In view of Proposition~\ref{propositionDFTcomputation}, there is
    no loss of generality in assuming $t_\uh\le s$.
    Lemma~\ref{lemDFTonfSul} together with~\eqref{ineqEvusencadr} show that it is enough to prove
    \begin{equation}
      J\subset \compl J^\uh \ens\text{and}\ens
      |J|\ge r-s \Imp
      \lim\limits_{\tau\rightarrow0} \Theta(\tau;\nu,g_J^\uh,2P) =0.
\end{equation} 
Equation~\eqref{equationlimitvalueTheta} with $\alpha=0$ gives
\begla
\lim\limits_{\tau\rightarrow0} \Theta(\tau;\nu,g_J^\uh,2P) =
-\frac{(2P)^\nu}{\nu+1} \sum_{m=1}^{2P}
B_{\nu+1}\!\Big(\frac{m}{M_{\alpha}}\Big)g_J^\uh(m),
\edla
where the $(\nu+1)^{\text{th}}$ Bernoulli polynomial has degree
$\nu+1\le r-s-1$.
The desired result is thus implied by 
\begin{equation}   \label{eqvanisha}
      J\subset \compl J^\uh \ens\text{and}\ens
      |J|\ge r-s \Imp
  \sum\limits_{m=1}^{2P} m^a g_J^\uh(m) = 0
  \ens \text{for each} \ens
  a\in\{0,\ldots, r-s-1\}.
\end{equation}

We prove~\eqref{eqvanisha} by exploiting~\eqref{ineqcrucialhj} as follows.
According to~\eqref{eqnewnumber}, 
we have
\begla
z^P \cdot
  \prod_{j\in \ccompl J} (z^{h_j\whp_j}+z^{-h_j\whp_j})
    \cdot \prod_{j\in J} (z^{h_j\whp_j}-z^{-h_j\whp_j})
    =
\sum_{\ue \in E} g_J^\uh\big(\cN^\uh(\ue)\big) z^{\cN^\uh(\ue)}
\edla
but we observe that, due to~\eqref{ineqcrucialhj},
$-P < \sum \eps_j h_j \wh p_j < P$ for each $\ue\in E$,
whence $0 < \cN^\uh(\ue) < 2P$. Lemma~\ref{lemmaHikamisets1}(ii) thus
yields
\beglab{eqnonreducedgP}
2^{t_\uh} \cdot\! \sum_{n=1}^{2P} g_J^\uh(n) z^n
= z^P \cdot
  \prod_{j\in \ccompl J} (z^{h_j\whp_j}+z^{-h_j\whp_j})
    \cdot \prod_{j\in J} (z^{h_j\whp_j}-z^{-h_j\whp_j}).
\edla
This is a reinforcement of~\eqref{equationgeneratinggs} inasmuch as we
just computed the \emph{non-reduced} generating polynomial
$\gP{{\uh,J}}(z) = \sum\limits_{n=1}^{2P} g_J^\uh(n) z^n$.
The evaluation at $z=1$ of $\big(z \frac d{dz}\big)^a\gP{{\uh,J}}(z)$
will give the sum in the \rhs\ of~\eqref{eqvanisha},
but~\eqref{eqnonreducedgP} shows that $z=1$ is a root of
multiplicity~$|J|$ for the polynomial~$\gP{{\uh,J}}$
and we thus get~$0$ for $a \le r-s-1 < |J|$.
\end{proof}

\begin{corollary}   \label{corvanish}
Let $s\in \{0,\ldots,r-3\}$. Let $\nu\in\{0,\ldots,r-s-2\}$
  satisfy $\nu\equiv r-s-1\ [2]$.
  Then the restrictions of $\rwh{m^sf^{\uo}\circ\cT}$ to the sets
  $\fS^{\sig_1(\ul)}$, $\ul\in\fL$, satisfy
  \begla
  \ul \notin \fR \Imp \lim_{\tau\to0}
  \Theta(\tau;\nu,\rwh{m^sf^{\uo}\circ\cT}\cdot \indusl,2P)
  = 0.
  \edla
\end{corollary}

\begin{proof}
  Suppose $\ul \notin \fR$.
  By~\eqref{eqneverone}--\eqref{equationRp}, there exists $\uh\in\fH$ such that
  $\uh\sim \sigma_1(\ul)$ and $\sum\limits_{j=1}^{r} \frac{h_j}{p_j}<1$.
  Proposition~\ref{propositionlimitvanish} shows that
  $\lim_{\tau\to0} \Theta(\tau;\nu,\rwh{m^sf^{\uo}
    \circ \cT}\cdot \induh,2P) =0$,
  but $\fS^\uh=\fS^{\sig_1(\ul)}$ by Lemma~\ref{lemmaHikamisets2}(ii).
  \end{proof}




\subsection{Vector-valued strong quantum modular forms arising from
  partial theta series}   \label{secVVsqmf}
  

We conclude this appendix with quantum modularity properties of
the partial theta series associated with the elements of the vector
space~$\cVs$ introduced in Corollary~\ref{corcVs}.
The aim of Section~\ref{secVVsqmf} is to explain the proof of

\begin{proposition} \label{propositionvectorvalue}
  Let $s,\nu\in\Z_{\ge0}$ satisfy $\nu+s\le r-3$ and $\nu+s\equiv r-3\ [2]$.
Then, for every $g\in\cVs$,
the function $\Theta(\cdot\,;\nu,g,2P)$ is a component of
a vector-valued depth $[\nu/2]+1$ strong quantum modular form on the
full modular group $\SL(2,\Z)$
  with quantum set~$\Q$ and weight $\nu+\frac{1}{2}$.
\end{proposition}

Note that, thanks to Remark~\ref{remparity}, \eqref{eqbSmed} and
Proposition~\ref{propConseqTh}, the function
$\Theta(\cdot\,;\nu,g,2P)$ has a 
resurgent-summable asymptotic expansion at each $\al\in\Q$ and, in
agreement with Remark~\ref{remsqmfupper}, it is more precisely the collection of
these asymptotic expansions that is a strong quantum modular form.

Proposition~\ref{propositionvectorvalue} will follow from a more
precise result, Corollary~\ref{corVVbasisWnu} below.
Statements and computations in this section will be eased by the use
of the
\emph{metaplectic double cover~$\TGA$ of $\GA:=\SL(2,\Z)$} (\cite{Weil64},
\cite{Shi73}, \cite{transseriescompletions}).
With reference to~\eqref{eqdefJga}, one may define this
group as
\begla
\TGA :=
\big\{\, (\ga,j) \in \GA\times\cO(\HH) \mid
j^2 = J_\ga \,\big\}
\ens \text{with product}\ens
(\ga_1,j_1) (\ga_2,j_2) \defeq \big(\ga_1\ga_2, \big(j_1\circ\ga_2) j_2\big).
\edla
A few particular elements of~$\wt\Ga$ are
\beglab{eqdefuoneTS}
\uone \defeq \Big( \matr{1&0}{0&1}, 1 \Big), 
\qquad
\uT \defeq \Big( \matr{1&1}{0&1}, 1 \Big),
\qquad
\uS \defeq \Big( \matr{0&-1}{1&0}, \tau^{1/2} \Big),
\qquad
\edla
where we use the principal branch in the latter case, \ie $\tau^{1/2}$
takes values in the first quadrant.
Note that $\uS^4$ is a nontrivial central element; multiplication
by~$\uS^4$ is the involution $(\ga,j)\mapsto (\ga,-j)$.
The group~$\TGA$ is generated by~$\uT$ and~$\uS$.

We call parabolic the elements $\uga=(\ga,j)\in\TGA$ for which $c=0$
when~$\ga$ is written as in~\eqref{eqdefleftactionGA};
we then have $a=d\in\{1,-1\}$, $\ga\tau=\tau+bd$ and $J_\ga(\tau)=d$,
whence the function~$j$ is constant with values in $\{1,-1\}$ or $\{i,-i\}$.

The advantage of~$\TGA$ over~$\GA$ is that
the weight~$w$ action from
the right of~$\GA$ on the space of all functions on~$\HH$,
$(\ga,\phi)\mapsto J_\ga^{-w} \, (\phi\circ\ga)$,
was defined for integer~$w$ but not for half-integer~$w$ in general,
whereas
\begla
\big(\phi,\uga\big) \mapsto j^{-2w} \cdot (\phi\circ\ga)
\ens\text{defines a right action of~$\TGA$ for any $w\in \tfrac12\Z$.}
\edla

Correspondingly, elaborating on \cite[Theorem~6]{han2022resurgence},
one can define an action of~$\TGA$ from the right, 
$(f,\uga) \mapsto f\bulM\uga$,
on the space
\begla
\cV := \{\, f\col \Z\to\C \mid
\text{$f$ is $2P$-periodic and $\sQ_{f,2P}=\Q$}\,\}
\edla
with the following properties:\footnote{We
  skip some details here and refer the interested reader to
  \cite{transseriescompletions}. One finds, for any $\uga\in\TGA$ and $f\in\cV^\pm$,
$\uga$ parabolic $\iimp f\bulM\uga = j^{-3} \La_{2P}^{db} f$ or
$j\ii \La_{2P}^{db} f$ according as~$f$ is odd or even, and
in the non-parabolic case:
\[
f\bulM\uga(n) =
 j(1-\tfrac d c)\ii (2P)^{-1/2} 
    e^{-i\pi/4} \La_{2P}^{bd}(n)
    \sum_{r \modM} f(r+dn)\, e^{i\pi{bnr}/{P}}
    \sum_{\substack{\ell\modcM \\ \text{s.t.}\;\ell=r\modM}}\La_{2cP}^a(\ell)
  \]
  for all $n\in\Z$,
with the notation $\La_M(n):=e^{i\pi n^2/M}$ for any positive even integer~$M$.
}

\noindent \textbf{(i)}
This action of~$\TGA$ is parity-preserving, i.e.\ it leaves invariant both subspaces
\begla
\cV^- := \{\, f\in\cV \mid \text{$f$ is odd}\,\}
\quad\text{and}\quad
\cV^+ := \{\, f\in\cV \mid \text{$f$ even}\,\}.
\edla

\noindent \textbf{(ii)}
For any $f\in\cV$,
\beglab{eqactionuTuS}
\big(f\bulM\uT\big)(n) = e^{\frac{i\pi n^2}{2P}} f(n),
\quad
\big(f\bulM\uS\big)(n) = e^{-i\pi/4} \wh f(n)
\qquad
\text{for all $n\in\Z$.}
\edla

 \noindent \textbf{(iii)}
 If $\uga=(\ga,j)\in\TGA$ is parabolic and $\nu\in\Z_{\ge0}$, then
 \begin{gather}
   \label{eqqmodparabod}
   f\in\cV^- \Imp
 \Theta(\cdot\,;\nu,f)
 - j^{-3}\, \Theta(\cdot\,;\nu,f\bulM\uga\ii)\circ\ga
 = 0,
   \\[.8ex]
   \label{eqqmodparabev}
   f\in\cV^+ \Imp
 \Theta(\cdot\,;\nu,f)
 - j\ii\, \Theta(\cdot\,;\nu,f\bulM\uga\ii)\circ\ga
 = 0.
 \end{gather}

\noindent \textbf{(iv)}
If $\uga=(\ga,j)\in\TGA$ is non-parabolic and~$\ga$ is written as
in~\eqref{eqdefleftactionGA}, then
 \begin{align}
   \label{eqqmodod}
   f\in\cV^- &\Imp
 \Theta(\cdot\,;0,f)
 \mp j\ii\, \Theta(\cdot\,;0,f\bulM\uga\ii)\circ\ga
 = \bS^{\frac\pi2\mp c\epsilon} \wt\Theta_{0,f,-\frac d c}\circ
 \big(\ID+\tfrac d c\big),
   \\[.8ex]
   \label{eqqmodev}
   f\in\cV^+ &\Imp
 \Theta(\cdot\,;1,f)
 \mp j^{-3}\, \Theta(\cdot\,;1,f\bulM\uga\ii)\circ\ga
 = \bS^{\frac\pi2\mp c\epsilon} \wt\Theta_{1,f,-\frac d c}\circ
 \big(\ID+\tfrac d c\big).
 \end{align}
 Note that the \rhs s of~\eqref{eqqmodod}--\eqref{eqqmodev} are
 independent of~$\epsilon$ provided $\epsilon>0$ is small enough.
 For instance, if $c>0$, then~\eqref{eqqmodod} with the `$-$' sign
 and~$f$ odd says that
 \[
 \Theta(\tau;0,f)
- j(\tau)\ii \, \Theta(\ga\tau;0,f\bulM\uga\ii)
 = \bS^{\frac\pi2-\epsilon} \wt\Theta_{0,f,-\frac d c}
 \big(\tau+\tfrac d c\big),
\]
i.e.\ we use the ``lateral sum to the right'', but if $c<0$ we must
use the ``lateral sum to the left'';
the former Borel-Laplace sum has a holomorphic extension to $\big(\!-\frac
d c,+\infty)$, while the latter one has a holomorphic extension to
$\big(\!-\infty,-\frac d c\big)$.
This is the key to the proof of
 
\begin{lemma}   \label{lemVVbasisW}
  Let $\nu\in\{0,1\}$.
Let~$\cW$ be any linear subspace
of~$\cV^\pm$ invariant under the action of~$\TGA$, where the
sign~$\pm$ is that of $(-1)^{\nu+1}$
(i.e.\ any $f\in\cW$ is odd if $\nu=0$, and is even if $\nu=1$).
Then, for any basis $(g_1,\ldots,g_D)$ of~$\cW$, the functions
  \begla
  \ph_i\pp\nu := \Theta(\cdot\,;\nu,g_i) \col \HH \to \C, \qquad i=1,\ldots,D,
  \edla
  are the components of a resurgent-summable quantum modular
  form on the full modular group $\GA=\SL(2,\Z)$
  with quantum set~$\Q$ and weight $\nu+\frac{1}{2}$: 
  there exists $\eps \col \GA \to \GL(D,\C)$ such that
  \begla
  (\ph_1\pp0,\ldots,\ph_D\pp0) \in \ocQ^1_{\frac{1}{2}}(\Q,\GA,\eps)^{res}_{med}
  \quad\text{if $\nu=0$}, \qquad
  (\ph_1\pp1,\ldots,\ph_D\pp1) \in \ocQ^1_{\frac{3}{2}}(\Q,\GA,\eps)^{res}_{med}
  \quad\text{if $\nu=1$.}  
  \edla
\end{lemma}

Here, we have introduced a reinforcement of Definition~\ref{defstrongqmfvar}:
\begin{definition}
  Given $N\ge0$, $w\in\frac12\Z$ and $\eps\col \GA\to\GL(D,\C)$,
the space $\ocQ^N_{w}(\Q,\GA,\eps)^{res}_{med}$ of depth~$N$ resurgent-summable quantum modular
  forms on~$\GA$
  with quantum set~$\Q$ and weight~$w$
  is defined to be~$\C$ if $N=0$ and, if $N\ge1$, the set of all tuples of holomorphic functions
  $(\ph_1,\ldots,\ph_D) \col \HH \to \C$ such that,
  for each $\al\in\Q$, $\ph_i\pp\nu(\tau)$ can be obtained as the median sum of a
  resurgent-summable formal series of $\C[[\tau-\al]]$
  and, for each $\ga\in\GA$, the modularity defect
         \begin{equation}   \label{eqvvmoddefR}
       \Big( \ph_i - J_\ga^{-w} \sum_{k=1}^D \eps_{i,k}(\ga) \ph_k\circ\ga \Big)_{1\leq i \leq D} 
       \ens \text{belongs to} \ens
       \bigoplus_{m=1}^M \cO(R_\ga) \otimes
       \ocQ_{w_m}^{N_m}(\Q,\GA,\eps^{(m)})^{res}_{med},
     \end{equation}
     %
     %
     where~$R_\ga$ is an open neighborhood of~$\R$ if $c=0$
     and an open neighborhood of $\R\setminus\{-\frac dc\}$ if
     $c\neq0$,
     following the convention~\eqref{eqchoosebranch} to determine
     $J_\ga^{-w}$ on~$R_\ga$,
for some $M\in\Z_{\ge1}$,
 $w_1,\ldots,w_M \in \frac12 \mathbb{Z}$ 
    and matrix-valued multipliers $\eps^{(1)},\ldots,\eps^{(M)}$, 
    and with $0\le N_m <N$ for each~$m$.    
\end{definition}

\begin{proof}[Proof of Lemma~\ref{lemVVbasisW}]
  Let $(g_1,\ldots,g_D)$ be a basis of a $\TGA$-invariant subspace
  $\cW\subset\cV^\pm$.
  For each $\uga\in\TGA$, let $A(\uga)\in \GL(D,\C)$ denote the matrix
  of the linear automorphism $f\in\cW \mapsto f\bulM \uga\in\cW$ in this
  basis, so that
  \begla
  g_i\bulM \uga = \sum_{k=1}^D  A_{i,k}(\uga) g_k,
  \qquad i=1,\ldots,D.
  \edla
Let $\ph_i\pp\nu := \Theta(\cdot\,;\nu,g_i)$ with $\nu=0$ or~$1$
according as the functions in~$\cW$ are odd or even.
We will prove that, for each
$\ga=\begin{psmallmatrix} a&b \\ c&d \end{psmallmatrix}\in\GA$,
%
%
each of the two lifts $\uga:=(\ga,j)\in\TGA$ of~$\ga$ satisfies
       \beglab{eqqumodvvnuzoro}
       \ph_i\pp\nu \mp j^{-2\nu-1} \sum_{k=1}^D \ue_{i,k}(\uga) \,\ph_k\pp\nu\circ\ga
       = \;\left| \begin{split}
         0 \qquad\qquad\qquad\qquad & \text{if $c=0$ and the sign `$\mp$' is `$-$'} \\[.8ex]
         \bS^{\frac\pi2\mp\epsilon} \wt\Theta_{\nu,f,-\frac d c}\circ \big(\ID+\tfrac d c\big)
         \quad & \text{if $c>0$}
       \end{split} \right.
     \edla
     with $\ue_{i,k}(\uga) :=d A_{i,k}(\uga\ii)$ if $c=0$, and $\ue_{i,k}(\uga) := A_{i,k}(\uga\ii)$ if $c>0$.
       
     This is sufficient because, when $c\ge0$, $J_\ga$ takes its
     values in the upper half-plane and we can thus choose the lift
     that has~$j$ taking its values in the first quadrant:
     setting $\eps_{i,k}(\ga):=\ue_{i,k}(\uga)$ with that choice
     of~$\uga$, we get a trivial modular defect~\eqref{eqvvmoddefR} on $R_\ga:=\R$ in the
     parabolic case (because~\eqref{eqchoosebranch} then says that
     $J_\ga(\al)^{1/2}\in\{1,i\}$, i.e.\ $J_\ga(\al)^{1/2}=j$)
     and, in the non-parabolic case, observing that for any
     $\al \in\Q\setminus\{-\tfrac dc\}$
     the non-tangential limit of 
     $j^{-2\nu-1}(\tau)$ as $\tau\to\al$ is
     $J_\ga(\al)^{-\nu-\frac12}$ if $\al>-\frac d c$
     and $-J_\ga(\al)^{-\nu-\frac12}$ if $\al<-\frac d c$
    (due to the convention~\eqref{eqchoosebranch}),
     we see that the modular defect~\eqref{eqvvmoddefR} is the restriction to
     $\Q\setminus\{-\tfrac dc\}$ of a function holomorphic in a
     neighborhood of
     $R_\ga := \R\setminus\{-\tfrac dc\}$.
     Moreover, Remark~\ref{remgamga} allows us to cover the case
     $c<0$ as well (we can compute $\eps_{i,k}(\ga)$ in terms of $\eps_{i,k}(-\ga)$).

     As for the proof of~\eqref{eqqumodvvnuzoro}, the case $c>0$ directly follows from
     \eqref{eqqmodod}--\eqref{eqqmodev};
     for the case $c=0$, use \eqref{eqqmodparabod}--\eqref{eqqmodparabev} noticing that $j^2=d=\pm1$.
\end{proof}


\begin{lemma}
Given $N\ge1$, $w\in\frac12\Z$ and $\eps\col \GA \to \GL(D,\C)$, we
have
\[
  (\ph_1,\ldots,\ph_D) \in \ocQ^N_{w}(\Q,\GA,\eps)^{res}_{med}
  \Imp
  \Big(\frac{d\ph_1}{d\tau},\ldots,\frac{d\ph_D}{d\tau}\Big) \in
  \ocQ^{N+1}_{w+2}(\Q,\GA,\eps)^{res}_{med}.
\]
\end{lemma}

\begin{proof}
Rephrasing the premise in terms of column vectors, we have
\beglab{eqdemoddefvvR}
\Phi - J_\ga^{-w} \, \eps \cdot \Phi\circ\ga =
\sum_{m=1}^M h_m \Phi\pp m,
\edla
 with $h_m\in \cO(R_\ga)$ and $\Phi\pp m \in
 \ocQ_{w_m}^{N_m}(\Q,\GA,\eps^{(m)})^{res}_{med}$,
    where~$R_\ga$ is an open neighborhood of~$\R$ or $\R\setminus\{-\frac dc\}$ and $M\in\Z_{\ge1}$,
    for some weights $w_1,\ldots,w_M \in \frac12 \mathbb{Z}$
    and matrix-valued multipliers $\eps^{(1)},\ldots,\eps^{(M)}$,
    with $0\le N_m <N$ for each~$m$.
    %
    %
    
    Differentiating with respect to~$\tau$, since
    $\frac{d}{d\tau}(\ga\tau) = J_\ga^{-2}$ and $\frac{d J_\ga}{d\tau}
    =c$, we get
    \begla
\frac{d\Phi}{d\tau} - J_\ga^{-w-2} \, \eps \cdot \frac{d\Phi}{d\tau} \circ\ga =
-cw J_\ga^{-w-1} \, \eps \cdot \Phi\circ\ga 
+ \sum_{m=1}^M \Big[ \frac{d h_m}{d\tau} \Phi\pp m + h_m \frac{d\Phi\pp m}{d\tau}\Big],
    \edla
    the desired result thus follows by induction on~$N$.  
\end{proof}


Since $\ph_i\pp{\nu+2} =
\dfrac{2P}{i\pi}\dfrac{d\ph_i\pp\nu}{d\tau\;}$ by the first part
of~\eqref{eqderivThetaj}, we immediately obtain from Lemma~\ref{lemVVbasisW}

\begin{corollary}   \label{corVVbasisWnu}
    Let $\nu\in\Z_{\ge0}$.
Let~$\cW$ be any linear subspace
of~$\cV^\pm$ invariant under the action of~$\TGA$, where the
sign~$\pm$ is that of $(-1)^{\nu+1}$
(i.e.\ any $f\in\cW$ is odd if~$\nu$ is even, and is even if~$\nu$ is odd).
Then, for any basis $(g_1,\ldots,g_D)$ of~$\cW$, the functions
  \begla
  \ph_i\pp\nu := \Theta(\cdot\,;\nu,g_i,2P), \qquad i=1,\ldots,D
  \edla
  are the components of a depth $[\nu/2]+1$ resurgent-summable quantum modular
  form on $\GA=\SL(2,\Z)$
  with quantum set~$\Q$ and weight $\nu+\frac{1}{2}$:
  \begla
  (\ph_1\pp\nu,\ldots,\ph_D\pp\nu) \in \ocQ_{\nu+\frac{1}{2}}^{[\nu/2]+1}(\Q,\GA,\eps)^{res}_{med}.
  \edla
\end{corollary}



We now prove Proposition~\ref{propositionvectorvalue} as follows.
We have seen that $\cVs\subset\cV^\pm$, where the sign is that of
$(-1)^{r-s}$.
Given a nonzero $g\in\cVs$, we consider the orbit $\gTGA$ of~$g$ under the 
action of~$\TGA$ (or rather the group algebra of~$\TGA$),
i.e.\ the minimal linear subspace of~$\cV^\pm$ that contains~$g$ and
is invariant under this action.
Note that
\begla
1 \le D := \dim_{\C} \gTGA \le \dim_{\C} \cV^\pm < P.
\edla
For each of the sequences $g_J^\ul$ that generate~$\cVs$,
we have
\begla
g_J^\ul \bulM \uT = e^{-2\pi iS_{\sig_1(\ul)} } g_J^\ul,
\quad
g_J^\ul \bulM \uS = e^{-i\pi/4} \,\rwh{\fakeheight{$g_J^\ul$}{$L_J$}}
\edla
by~\eqref{eqactionuTuS}, Lemma~\ref{lemmaHikamisets2}(iv) and
Remark~\ref{remCSHikamisets}.
Since~$\cVs$ is invariant under DFT (Corollary~\ref{corcVs}), it is
thus invariant under the action of~~$\uT$ and~$\uS$, and thus under
the action of $\TGA$ because the
group is generated by~$\uT$ and~$\uS$.
Therefore, we can apply Corollary~\ref{corVVbasisWnu} with 
$\cW = \gTGA\subset \cVs$
and any basis $(g_1,\ldots, g_D)$ of $\gTGA$ such that $g_1=g$.


\begin{remark}
  \label{remcongrsubgp}
  Any congruence subgroup $\Ga \subset\GA$ can be lifted to a subgroup
  $\wt\Ga\subset \TGA$.
  The restriction of the action of~$\TGA$ to~$\wt\Ga$ may have many
  more invariant subspaces.
  For instance, with $\Ga = \Ga_1(4P)$, one finds that every nonzero
  $g\in\cV^\pm$ gives rise to an invariant line $\C g$.
  It follows that, if $\nu\ge0$ and~$g$ have opposite parities, then
  $\Theta(\cdot\,;\nu,g,2P)$ is a (scalar) depth $[\nu/2]+1$ resurgent-summable quantum modular
  form on $\Ga_1(4P)$
  with quantum set~$\Q$ and weight $\nu+\frac{1}{2}$.
  This is the mechanism behind the proof of
  Corollary~\ref{corollarypartialthetaserieshigherdepth}.
  \end{remark}

\section*{Acknowledgements}

This paper is partly a result of the ERC-SyG project, Recursive and
Exact New Quantum Theory (ReNewQuantum) which received funding from
the European Research Council (ERC) under the European Union's Horizon
2020 research and innovation programm under grant agreement No.~810573.
This work has been partially supported by the project
CARPLO of the Agence Nationale de la recherche (ANR-20-CE40-0007).
Y.~Li is partially supported by BNSFC No.~JR25001.
S.~Sun is  partially supported by National Key R\&D Program of China
(2020YFA0713300), NSFC (Nos.~11771303, 12171327, 
6111911530092, 12261131498, 11871045).


\printbibliography


\vspace{3ex}

\noindent 
\textbf{J\o rgen Ellegaard Andersen}\\
Center for Quantum Mathematics, and
Danish Institute for Advanced Studies\\
University of Southern Denmark,
DK-5230 Odense, Denmark

\smallskip

\noindent 
\textbf{Li Han}\\
Department of Mathematics\\ Capital Normal University, Beijing 100048, P.~R.~China

\smallskip

\noindent
\textbf{Yong Li}\\
Beijing Institute of Mathematical Sciences and Applications, Beijing 101408, P.~R.~China

\smallskip

\noindent 
\textbf{William Elb\ae k Misteg\aa rd}\\
Center for Quantum Mathematics\\
University of Southern Denmark,
DK-5230 Odense, Denmark

\smallskip

\noindent
\textbf{David Sauzin}\\
Department of Mathematics\\
Capital Normal University, Beijing, P.~R.~China\\
(on leave from CNRS LTE, Paris Observatory, France)

\smallskip

\noindent
\textbf{Shanzhong Sun}\\
Department of Mathematics, and Academy for Multidisciplinary Studies\\
Capital Normal University, Beijing 100048 P.~R.~China

\end{document}